\documentclass[12pt,a4paper]{article}
\usepackage{amsmath,amscd,amssymb,amsthm,amstext,mathrsfs,dcolumn,
booktabs,enumerate}

\usepackage{geometry}
\usepackage{ifpdf}
\ifpdf
\usepackage[all,line,cmtip]{xy}
\else
\usepackage[all,line,cmtip,dvips]{xy}
\fi

\CompileMatrices

\usepackage[final]{showkeys}

\input{xypic}

\theoremstyle{plain}
   \newtheorem{theorem}{Theorem}[section]
   \newtheorem{proposition}[theorem]{Proposition}
   \newtheorem{lemma}[theorem]{Lemma}
   \newtheorem{corollary}[theorem]{Corollary}
   \theoremstyle{definition}
   
   \newtheorem{definition}[theorem]{Definition}
   \newtheorem{example}[theorem]{Example}
   \theoremstyle{remark}
   
   \newtheorem{remark}[theorem]{Remark}

\newcommand{\Th}{{\mathrm{Th}}}

\newcommand{\RR}{{\mathbb R}}
\newcommand{\NN}{{\mathbb N}}
\newcommand{\ZZ}{{\mathbb Z}}
\newcommand{\QQ}{{\mathbb Q}}
\newcommand{\CC}{{\mathbb C}}
\newcommand{\PP}{{\mathbb P}}
\newcommand{\FF}{{\mathbb F}}

\newcommand{\Gr}[2]{{{\mathbf G}_{#1}(\CC^{#2})}}
\newcommand{\St}[2]{{{\mathbf V}_{#1}(\CC^{#2})}}
\newcommand{\Fr}[2]{{{\mathbf{PV}}_{#1}(\CC^{#2})}}
\newcommand{\Fra}[3]{{{\mathbf{PV}}_{#1}^{(#3)}(\CC^{#2})}}
\newcommand{\Fratilde}[3]{{{\widetilde{\mathbf{PV}}}_{#1}^{({#3})}(\CC^{#2})}}
\newcommand{\SG}[1]{{\mathbf{\Delta}_o(\CP^{#1} )}}
\newcommand{\PaG}[2]{{\mathbf{B}_{#2} (\CP^{#1} )}}
\newcommand{\Id}{{\mathbf{Id}}}

\newcommand{\CP}{{\mathbb{C}\mathrm{P}}}

\newcommand{\coker}{{\mathrm {coker}}}
\newcommand{\im}{{\mathrm {im}}}
\newcommand{\FC}{{\mathcal F}}

\newcommand{\TT}{{\mathbb T}}

\newcommand{\F}{{\mathcal F}}
\newcommand{\Alg}{{{\mathcal A}lg}}

\newcommand{\Tor}{{\operatorname{Tor}}}

\newcommand{\simp}{_\bullet}

\newcommand{\Hom}{{\operatorname{Hom}}}
\newcommand{\Ext}{{\operatorname{Ext}}}

\newcommand{\image}{{\operatorname{im}}}

\newcommand{\drf}{{\overline \Omega}}

\newcommand{\PS}{{P}}

\newcommand{\trunc}[2]{{T_{#1}(#2)}}

\newcommand{\ddr}{{\bf d}}

\newcommand{\trcl}[2]{{\tau_{#1}^{(#2)}}}

\newcommand{\Morse}{{\mathcal M}}

\newcommand{\Trans}{{\mathcal T}}

\newcommand{\IT}{{\mathcal IT}}
\newcommand{\IF}{{\mathcal IF}}

\newcommand{\inp}[2]{{\langle #1, #2 \rangle}}

\newcommand{\diag}{{\mathrm{diag}}}

\title{String cohomology groups of complex projective spaces}
\author{Marcel B\" okstedt \& Iver Ottosen \footnote{The main part of
I. Ottosen's work was done at the University of Aarhus. A smaller part 
was done at the Institute Mittag-Leffler, Sweden and at 
the University of Aalborg.}}

\begin{document}
\maketitle

\begin{center}
Department of Mathematical Sciences, University of Aarhus, \\
Ny Munkegade, 8000 Aarhus C, Denmark. \\
marcel@imf.au.dk
\end{center}

\begin{center}
Department of Mathematical Sciences, University of Aalborg, \\
Fredrik Bajers Vej 7 G, 9220 Aalborg East, Denmark. \\
ottosen@math.aau.dk
\end{center}

\begin{abstract}
Let $X$ be a space and write $LX$ for its free loop space equipped
with the action of the circle group $\TT$ given by dilation.
We compute the equivariant cohomology  
$H^*(LX_{h\TT};\ZZ /p)$ as a module over $H^* (B\TT ;\ZZ /p)$ 
when $X=\CP^r$ for any positive integer $r$ and any prime 
number $p$. The computation implies that the associated 
mod $p$ Serre spectral sequence collapses from the $E_3$-page.
\begin{noindent}
\end{noindent}
\end{abstract}

\section{Introduction}

Let $LX$ be the space of maps from the circle to a space
$X$. The circle acts on itself by rotation, and this action 
induces an action of the circle group ${\TT}=S^1=SO(2)$ on $LX$. 
(The action extends to an $O(2)$-action, but we will not consider
the extended action in this paper).
The homotopy orbit under the circle action is
the space  $E{\TT}\times_{\TT} LX$, which we will
also write as $LX_{h\TT }$.
The purpose of this paper is to compute the ${\TT}$ 
Borel cohomology with $\FF_p$ coefficients of the free loop space 
on $\CP^r$, that is $H^*(L \CP^r_{h\TT };\FF_p )$. 

There are several motivations for studying this question.
In differential topology one studies the spectrum $TC(M)$,
which is related to the diffeomorphisms of the manifold $M$.
There is a long exact sequence for the spectrum cohomology 
of $TC(M)$, where the other terms are given by the cohomology
and the Borel cohomology of the free loop space on $M$. 
So this result gets us closer to understanding $TC(\CP^r )$.

The theory of Chas and Sullivan \cite{SullivanChas} constructs 
algebraic operations on groups related to the free loop space. 
Part of this structure has been computed in the case $M=\CP^r$ 
\cite{CohenJones}. It would be interesting to write down as much 
as possible of the Chas-Sullivan structure for the particularly 
simple example $M=\CP^r$. Since one of the groups carrying the 
structure is the Borel cohomology of the free loops space, a first 
step is to compute this group. 

The usual cohomology of the free loop space has been used 
to study the existence of closed geodesics 
\cite{GM}, \cite{Klingenberg}, \cite{VigueSullivan}. The main
problem studied is to decide whether any Riemannian manifold
has infinitely many geometrically distinct closed geodesics.
This method works better the bigger the cohomology of the free loop space is.
For general metrics on $\CP^r$ one has not succeeded. 
It is a far out possibility that the knowledge of the equivariant 
cohomology could possibly be of some help.

We should point out that our computation does not give us much more 
than the cohomology groups as groups. One can ask for product 
structure, cohomology operations and higher torsion. We are working 
on these questions.  

The method we use is convoluted, so we should try to
give an overview of the computation here.

We know several different spectral sequences that converges to 
$H^*(L\CP^r_{h\TT };\FF_p )$. The first is the Serre spectral sequence 
for $L\CP^r \to L\CP^r_{h\TT } \to B{\TT}$. Since the cohomology of 
$L\CP^r$ is known, we know the $E_2$ page of this spectral sequence. 
The $d_2$ differential is given if you can compute the action map 
${\TT} \times L\CP^r \to L\CP^r$ in the ordinary cohomology
of the free loop space. We are able to compute this directly, using
methods from \cite{SpSe}, so we can compute the $E_3$ page of the 
Serre spectral sequence. It is a consequence of our final result that 
there are no further non-trivial differentials, so that the Serre 
spectral sequence collapses from the $E_3$ page. One could argue 
that this is the simplest formulation of our main result. We have no 
clue as to how to prove this directly using only the Serre spectral 
sequence.

The second spectral sequence is derived from infinite dimensional 
Morse theory. For a compact Riemannian manifold $M$, the energy 
function $E:LM\to \RR$ is a ${\TT}$-equivariant map, so we get 
an equivariant filtration of $LM$ by energy levels. 

For $M=\CP^r$ with the symmetric space (Fubini--Study) metric,
the energy function is a Morse--Bott function. The filtration defines 
a spectral sequence. We have enough knowledge of the critical manifolds 
to write down the spectral sequence. This method has previously been 
used by Klingenberg in  \cite{KlingProj} to compute
$H^*(EO(2)\times_{O(2)} L \CP^r;\FF_2 )$. 
In this case, Klingenberg proves that the spectral sequence collapses. 

In our case, the spectral sequence will eventually collapse, but not 
from the $E_1$ page, but only from the $E_p$-page.
There are definitely non-trivial differentials in the Morse
spectral sequence. Not all of these differentials start
at the $\TT$ fixed points, that is at the space of constant curves.
 
The non-triviality of the differentials can be interpreted
as a geometrical statement about $n$-fold iterated geodesics. 
It is a consequence of our calculation that on $\CP^r$
with the standard metric for any $n$ there is 
always a curve close to an $n$-fold geodesic, such that if the curve 
moves according to the dynamics given by the gradient of the energy 
function $E$, then it will eventually approximate an $n-1$-fold 
iterated geodesic. 

We study the Morse filtration spectral sequence using localization 
methods from equivariant homotopy theory. The $p$-fold iteration 
map induces a map 
\[
H^*(LM_{h\TT } ;\FF_p ) 
\to H^*(LM_{h\TT }^{(p)} ;\FF_p)
\]
where the superscript $(p)$ means that we have twisted the 
$\TT$-action. This map is compatible with the Morse filtration,
so it induces a map of spectral sequences.
The localization methods prove that under certain technical
conditions, after inverting the generator $u$ of $H^*(B\TT ;\FF_p )$ 
this map induces an isomorphism of spectral sequences. It turns out 
that the Morse spectral sequence for the twisted action collapses, 
so this isomorphism puts strong conditions on the original spectral 
sequence.

We apply this to our special case. After some brute calculation 
in the Morse spectral sequence we obtain the collapse from the $E_p$ 
page.

So we have to contend with the first differentials.
We know that there are some non-trivial differentials, because the $E_1$
page of the Morse spectral sequence is strictly bigger than 
the $E_3$ term of the Serre spectral sequence. To play these 
two spectral sequences out against each other, we study the 
$\TT$-transfer. This is a degree $-1$ endomorphism on $H^*(X)$ 
defined for any $\TT$-space $X$. What makes it interesting is that we 
can compute it using homotopy theoretical methods, but it seems 
hard to calculate the transfer directly in the Morse theory context,
primarily because the transfer will not 
respect the splitting of $H^*(L\CP^r ;\FF_p )$ into summands 
corresponding to the quotients of the Morse filtration. So you 
can translate information about the $\TT$-transfer into information 
about how the Morse strata fit together, that is information about 
the heteroclinic trajectories of the flow on $L\CP^r$ defined by the 
gradient of the energy function.

A third spectral sequence converging to
$H^*(L\CP^r_{h{\TT}};\FF_p )$ is the spectral sequence that
we developed in \cite{SpSe} and \cite{BO}.
For any space $X$ the $E_2$ term of this spectral 
sequence is algebraically described as non-abelian higher derived 
functors. If $X$ is simply connected, it converges towards
$H^*(LX_{h{\TT}};\FF_p )$. For general $X$ it seems to be a hard 
algebraic problem to compute the $E_2$ page explicitly, but we can handle 
the algebra for $\CP^r$. Preliminary computations show that the 
size of the $E_2$ term of this spectral sequence agrees with the 
size of the Morse spectral sequence, which means that it collapses 
from the $E_2$ page.

The individual sections are organized as follows.
In sections \ref{sec:Geospace} and \ref{sec:CohGeospace} we study 
the space of geodesics on $\CP^r$. Except for the constant ``geodesics'',
this space consists of an infinite number of homeomorphic components.
We compute its cohomology by comparing to the cohomology of Grassmann 
spaces. The $\FF_p$-cohomology of any $\TT$-space is a module over
$H^*(B\TT ;\FF_p )\cong \FF_p [u]$, and we determine these
module structures.
This is very closely related to Klingenberg's work in \cite{KlingProj}. 
The main difference is that we are interested in the circle action, 
while he studies the action of the group $O(2) \supset \TT$. 

In section \ref{sec:Orbits} we compute the equivariant (Borel) 
cohomology of the spaces of geodesics. The action of $\TT$ is different
on the different homeomorphic components.
The equivariant cohomology can distinguish at 
least some of them.

In section \ref{sec:Bundles} we recall some facts about
equivariant cohomology theory, in particular the
localization theorem, which shows that the relative equivariant
cohomology of $X^{C_p}\subset X$ is annihilated by $u$.

In section \ref{sec:Twisting} we study the following abstract 
construction. If $X$ is a $\TT$ space, then the fixed points 
$X^{C_p}$ is a $\TT \cong \TT/C_p$ space. We show that the 
equivariant cohomology of $X^{C_p}$ with coefficients in $\FF_p$
does not depend on the action. It is just the tensor product of 
the ordinary cohomology of $X$ with $\FF_p [u]$. 

In section \ref{sec:Iteration} we apply the general theory of
the two previous section to the free loop space of a manifold. 
The idea is that the iteration map, which maps a loop to the 
same loop run through $p$ times, is really the inclusion of the 
$C_p$ fixed points, so that it induces an isomorphism of localized
equivariant cohomology and by the results of \ref{sec:Twisting} 
it is easy to compute this equivariant cohomology.

The localization theorem needs a finite dimensionality condition. 
This condition is not just technical, but essential. Since the
free loop space is infinite dimensional, we have to reduce the 
situation to a finite dimensional situation. We do this step 
using Morse theory on the Hilbert manifold version of the free 
loop space.
 
In section \ref{sec:MSS} we analyze how the theory from section 
\ref{sec:Iteration} influences the Morse spectral sequence. We show 
that under certain (strong) conditions, the localization of the
Morse spectral sequence for equivariant cohomology is a functor 
in the non-equivariant Morse spectral sequence.

In section \ref{sec:MorseCPr} we specialize to $\CP^r$. We use the 
results from the sections \ref{sec:CohGeospace} and \ref{sec:Orbits}
to explicitly write down the $E_1$ page of the Morse spectral sequence, 
in both the equivariant and the non-equivariant case.

For this space we do know that the non-equivariant spectral sequence 
collapses, and it follows from section \ref{sec:Iteration} that 
the localization of the Morse spectral sequence collapses.
This does not imply that the un-localized Morse spectral sequence
for the equivariant cohomology collapses from the $E_1$ page,
but it is sufficient to prove collapsing from the $E_p$ page.

Now there is a change of scene to homotopy theoretical methods. 
In section \ref{sec:FreeCohomology} 
we compute the $\TT$ transfer map. We do this in the slightly greater 
generality of spaces whose cohomology are truncated polynomial algebras 
on one generator. To do the computation, we use the homological
methods we developed in \cite{SpSe}. The central computation is
a computation in non-abelian homological algebra, which we
do in section \ref{sec:Derived}.

Even if it is not necessary for the proof of our main result, we
study the Serre spectral sequence converging to Borel 
homology in section \ref{sec:SerreSS}. In particular we 
show that our main result implies that the Serre spectral sequence
collapses from the $E_3$ page. We also obtain some information about
the product structure in $H^*(L\CP^r_{h\TT} ;\FF_p )$, which seems
hard to get using only the Morse theory approach. 

In section \ref{se:MorseSerreSS} we show how the transfer map 
forces the existence of non-trivial differentials in the 
Morse spectral sequence. It seems to be hard to  pinpoint these
differentials exactly, but we obtain sufficient information to
compute the $E_\infty$ term as a module over $\FF_p [u]$.  

{\em Notation}: We fix a prime $p$. Cohomology groups will always 
have coefficients in the field $\FF_p$, unless we explicitly 
state otherwise. We state otherwise mainly in section 
\ref{sec:CohGeospace}. The group $\TT$ is the circle group. 
Spaces are supposed to be homotopy equivalent to CW-complexes, and
$\TT$ spaces are supposed homotopy equivalent to $\TT$ CW-complexes.

\section{Geodesics on $\CP^r$}
\label{sec:Geospace}

Let us fix some notation and recall some basic facts. 
Consider $\CC^{r+1}$ with the usual Hermitian inner product  
${\inp \medspace \medspace}_\CC$, and write
$\inp \medspace \medspace$ for its real part. 
Thus $\inp \medspace \medspace$ is the standard inner product on 
$\RR^{2r+2} \cong \CC^{r+1}$.

Let $S^{2r+1}$ be the unit sphere in $\CC^{r+1}$.
Since the left action of $\TT$ on $S^{2r+1}$ given by 
$(z,v)\mapsto zv$ is a free proper action, the Hopf map
\[ Q:S^{2r+1} \to S^{2r+1}/\TT =\CP^r \]
is a smooth submersion. We equip $\CP^r$ with the Fubini Study 
Hermitian metric. This means that for any $x\in S^{2r+1}$ the map
\begin{equation} \label{a_x}
a_x : (\CC x)^\perp \subseteq T_x(S^{2r+1}) \xrightarrow{T_xQ} 
T_{Q(x)}(\CP^r ) 
\end{equation}
is a $\CC$-linear isometry, where 
$(\CC x)^\perp = \{ v \in \CC^{r+1} | {\inp v x}_\CC =0\}$.
Note that the following identity holds
\begin{equation} \label{fseqn}
a_{zx}(zv)=a_x(v) \text{ for } z\in S^1.
\end{equation}
See for example lemma \cite[14.4]{MT} for these facts. 

We equip $S^{2r+1}$ with the Riemannian metric coming from
$\inp \medspace \medspace$ and use the real part 
of the Fubini Study metric as Riemannian metric on $\CP^r$. 
Let $H_x$ be the horizontal subspace ie. the orthogonal 
complement of $\ker (T_xQ)=\RR ix$ in $T_x(S^{2r+1})$. Since 
\[ H_x = \{ v \in \CC^{r+1} | \inp v x = \inp v {ix} =0 \} 
= (\CC x)^\perp \]
we see that $a_x : H_x \to T_{Q(x)}(\CP^r )$ is an 
$\RR$-linear isometry. So the Hopf map $Q$ is a 
Riemannian submersion \cite[Ch. 2]{GHL}. 

We use the Levi--Civita connection corresponding
to the standard metrics on $S^{2r+1}$ and $\CP^r$.
Let $\PaG r 1$ denote the set of primitive closed geodesics
$f:[0,1]\to \CP^r$ and let $\PaG r q$ for integer $q\geq 1$ 
be the set of geodesics of the form $f(qt)$ for $f\in \PaG r 1$.
Note that a geodesic in $\PaG r q$ has length $q\pi$ and energy
$q^2\pi^2$.
The free loop space $L\CP^r$ defined as the set of  
closed curves $g:[0,1] \to \CP^r$ of class $H^1$ is a Hilbert 
manifold modeled on $L(\RR^{2r})$ \cite[\S 1]{KlingSphere}
and by \cite[\S 1.4]{KlingProj} we have that $\PaG r q$ is a (critical) 
submanifold of $L\CP^r$.

The left $\TT$ action on $L\CP^r$ restricts to an action on $\PaG r q$.
We view an element in $\PaG r q$ as a periodic geodesic and then
\[ \TT \times \PaG r q \to \PaG r q ; \quad 
(e^{2\pi i\theta}*f)(t)=f(t-\theta), \quad \theta \in \RR .\] 

We will now give an alternative description of $\PaG r q$.
The description can be found in \cite[\S 1]{KlingProj}, but we will
write down some explicit maps and add information on the $\TT$-action.
Write $\St 2 {r+1}$ for the Stiefel manifold of complex 
{\em orthonormal} 2-frames in $\CC^{r+1}$. We can rotate a 
2-frame $(x,v)$ by an angle $\omega \in \RR$ as follows:
\[ R(\omega )(x,v) = 
(\cos (\omega )x+\sin (\omega )v,-\sin (\omega )x+\cos (\omega )v). \]
Write $\Fr 2 {r+1}$ for the complex projective Stiefel manifold
$\St 2 {r+1} /\diag_2 (U(1))$.

\begin{definition}
For integer $q\geq 1$ we let $\Fra 2 {r+1} q$ denote $\Fr 2 {r+1}$ 
equipped with the well-defined $\TT$-action 
\[ \TT \times \Fr 2 {r+1} \to \Fr 2 {r+1}; \quad
e^{2\pi i \theta}*[x,v] = [R(-q\pi \theta )(x,v)], \quad \theta \in \RR .\]
\end{definition}

\begin{proposition}
There is a diffeomorphism  
$\phi_q : \Fra 2 {r+1} q \to \PaG r q$ defined by 
$\phi_q([x,v])=Q\circ c(q,x,v)$ where
\[ c(q,x,v)(t)= \cos (q\pi t)x+\sin (q\pi t)v, \quad 0\leq t\leq 1. \]
Furthermore, $\phi_q$ is a $\TT$-equivariant map.
\end{proposition}

\begin{proof}
By \cite[2.109 and 2.110]{GHL} the map $\phi_q$ is a bijection.
A direct computation using the trigonometric addition formulas 
shows that $\phi_q$ is $\TT$-equivariant.
\end{proof}

Note that the geodesics $Q\circ c(q,x,v)$ can easily be extended 
from $[0,1]$ to the open interval $(-\epsilon , 1+\epsilon )$ 
for $\epsilon>0$.

\begin{lemma} \label{sigma}
There is a $\TT$-equivariant diffeomorphism 
\[ \sigma_q :\Fra 2 {r+1} q \to \Fratilde 2 {r+1} q ; \quad
[x,v]\mapsto [\frac {x+iv} {\sqrt 2} ,\frac {x-iv} {\sqrt 2}], \]
where $\Fratilde 2 {r+1} q = \Fr 2 {r+1}$ equipped with 
the well-defined $\TT$-action 
\[ e^{2\pi i \theta } \star [a,b]=
[e^{-\pi q i \theta }a,e^{\pi q i \theta }b ].\]
\end{lemma}

\begin{proof}
We have an automorphism $\sigma$ of $\St 2 {r+1}$ as follows: 
\[ \sigma (x,v)= \frac 1 {\sqrt 2} (x+iv,x-iv) \quad , \quad 
\sigma^{-1}(a,b)=\frac 1 {\sqrt 2} (a+b,-i(a-b)). \]
This automorphism respects the diagonal $U(1)$ action. Furthermore,
\[ \sigma (R(\omega )(x,v))= 
\frac 1 {\sqrt 2} (e^{-i\omega }(x+iv),e^{i\omega }(x-iv)). \]
The result follows.
\end{proof}

Note that the $\TT$-action on $\Fratilde 2 {r+1} 1$ is free so by
the lemma the $\TT$-action on $\Fra 2 {r+1} 1$ is also free.

We have yet another interpretation of the space of primitive geodesics.
Let $\eta: S(\tau (\CP^r )) \to \CP^r$ denote the sphere bundle 
of the tangent bundle.

\begin{proposition}
\label{lemma:UnitTangents}
There is a diffeomorphism $\psi :\Fr 2 {r+1} \to S(\tau (\CP^r ))$
defined by
\[ 
\psi ([x,v])=(a_x(v))_{Q(x)}.
\]
The composition of $\eta \circ \psi$ is projection on the first factor 
$[x,v] \mapsto [x]$.
\end{proposition}

\begin{proof}
The map $\psi$ is well-defined by (\ref{fseqn}) and has the inverse
$\psi^{-1}(v_{Q(x)})=[x,a_x^{-1}(v)]$. 
Both $\psi$ and $\psi^{-1}$ are smooth maps.
\end{proof}

We equip $S(\tau (\CP^r ))$ with the $\TT$-action which makes
$\psi$ an equivariant map.
There is a fiber bundle sequence
$\Fra 2 2 1\to \Fra 2 {r+1} 1 \to \Gr 2 {r+1}$ since the Grassmannian 
is the quotient of the Stiefel manifold by the group $U(2)$. 
Equivalently, we have a fiber bundle
\[ S(\tau (\CP^1 )) \to S(\tau (\CP^r )) \to \Gr 2 {r+1} .\]
The action of $\TT$ on the total space is free and preserves the fibers.
After dividing by it, we get another fiber bundle

\begin{equation}
\label{eq:geodesicsandgrassman}
\SG 1 \to \SG r \to \Gr 2 {r+1}.  
\end{equation}
Note that we may view a point in $\SG r$ as the trace of a simple 
geodesics together with an orientation. Klingenberg \cite{KlingProj} does
not consider $\SG r$. He divides out by the entire $O(2)$ action instead
and denote the resulting quotient space by $\mathbf \Delta (\CP^r )$.

\begin{example}
\label{example:S2}  
In the case $r=1$ we have that $\CP^1$ is the standard round 
sphere $S^2$ with radius $\frac 12$. A tangent vector $w\in T_v(S^2)$ is part of a unique 
ordered basis $(v,w,v\times w )$. There is a unique $A\in SO(3)$ 
so that $v=Ae_1$, $w=Ae_2$, $v\times w=Ae_3$. This establishes 
a diffeomorphism between $SO(3)$ and $S(\tau (S^2))$.
The induced action of $e^{i\theta} \in S^1$ on $SO(3)$ 
is given by right multiplication $A\mapsto A\rho (\theta )$ 
where $\rho (\theta )$ is the rotation matrix by angle 
$\theta$ around the axis $e_3$.
The map $A\mapsto Ae_3$ is invariant under the action. 
It defines a map $\SG 1\to S^2$ which is a diffeomorphism.
\end{example}

The example combines with the fibrations to show that we have a 
diagram of fibration sequences
\begin{equation}
\label{eq:fiberdiagram}
\xymatrix@C=1.5 cm{
\TT \ar@{=}[r] \ar[d] & \TT \ar[r] \ar[d] & {*} \ar[d] \\
SO(3) \ar[r] \ar[d] & S(\tau (\CP^r )) \ar[r]^{q} \ar[d]^{q^\prime} 
& \Gr 2 {r+1} \ar@{=}[d] \\
S^2 \ar[r] & \SG r \ar[r]^{\overline q} & \Gr 2 {r+1} 
} \end{equation}

Let $\gamma_2$ be the canonical complex 2-dimensional vector bundle 
over $\Gr 2 {r+1}$ and let $\PP (\gamma_2 )$ be the associated projective 
bundle. A point in the total space of $\PP (\gamma_2 )$ is a flag 
$V_1\subset V_2\subset \CC^{r+1}$ where 
$V_i$ has complex dimension $i$.

\begin{lemma}
\label{lemma:ProjectiveBundle}
The fiber bundle $S^2\to \SG r \to \Gr 2 {r+1}$ 
is isomorphic to the fiber bundle 
$\CP^1 \to \PP (\gamma_2) \to \Gr 2 {n+1}$.
\end{lemma}

\begin{proof}
It is enough to show that there is a bundle map 
$\Fratilde 2 {r+1} 1 /\TT \to \PP (\gamma_2 )$ which is a fiber wise 
diffeomorphism. We have a smooth map 
\[ f: \St 2 {r+1} \to \PP (\gamma_2 )\quad ;\quad
(a,b)\mapsto (\CC \{a \} \subset \CC \{ a, b\} \subset \CC^{r+1}), \]
which factors through $\Fratilde 2 {r+1} 1 /\TT$ since multiplying 
the generators by units does not change a linear span. It suffices 
to see that $\overline f:\Fratilde 2 {r+1} 1 /\TT \to \PP (\gamma_2 )$
is a diffeomorphism when restricted to the fibers over 
$\CC^2 \subset \CC^{r+1}$ since we can then get the result for a 
general fiber by changing basis. The map of fibers is
$\Fratilde 2 2 1 \to \CP^1$ which under the standard identification
$U(2)/(U(1)\times U(1))\cong \CP^1$ corresponds to the quotient map
\[ (\frac {U(2)} {\diag_2 (U(1))})/\TT \to 
\frac {U(2)} {U(1)\times U(1)} .\]
This map is a diffeomorphism because of the following identity of
$2\times 2$ diagonal matrices:
\[ D(e^{2\pi i \alpha }, e^{2\pi i \beta })=
D(e^{\pi i (\alpha + \beta )},e^{\pi i (\alpha + \beta )})
D(e^{\pi i (\alpha - \beta )},e^{\pi i (\beta - \alpha )}); \quad
\alpha , \beta \in \RR .\]
\end{proof}

\section{The cohomology of spaces of geodesics}
\label{sec:CohGeospace}

The purpose of this section is to compute the cohomology
of $\SG r$ and $S(\tau (\CP^r ))$. It turns out to be
convenient to do this computations for cohomology with
coefficients in the integers. We first determine the 
cohomology of the base Grassmann manifold. 

\begin{theorem}
\label{th:GrassmannCohomology}
Let $c_1,c_2$ be the two first Chern classes 
of the canonical bundle $\gamma_2$ over $\Gr 2 {r+1}$. Then one has
\[
H^*(\Gr 2 {r+1} ;\ZZ )\cong \ZZ [c_1,c_2]/(\phi_r ,\phi_{r+1} ),
\]
where $\phi_i \in \ZZ [c_1,c_2]$ is the polynomial defined 
inductively by
\begin{displaymath}
\phi_0 = 1; \quad \phi_1 = -c_1; \quad 
\phi_i = -c_1 \phi_{i-1} - c_2 \phi_{i-2}.
\end{displaymath}
\end{theorem}

\begin{proof}
According to \cite{Borel}, proposition 31.1, we have an isomorphism
\begin{equation} \label{eq:Borel}
H^*(\Gr n {n+m} ;\ZZ ) \cong
\frac {\ZZ [x_1, \dots ,x_n]^{\Sigma_n} \otimes
\ZZ [x_{n+1}, \dots , x_{n+m}]^{\Sigma_m}} 
{(\ZZ [x_1, \dots ,x_{n+m}]^{\Sigma_{n+m}})^+}
\end{equation}
Here the degree of $x_i$ is 2 for all $i$ and $\Sigma_k$ denotes
the symmetric group. The plus sign means forming the ideal
generated by elements in positive degrees. 
By \cite{Borel}, 30.2 one sees that this
isomorphism comes from the fibration
\[ \frac {U(n+m)} {U(n)\times U(m)} \to 
\frac {EU(n+m)} {U(n)\times U(m)} \to 
\frac {EU(n+m)} {U(n+m)}, \]
or equivalently the fibration
\[ \Gr n {n+m} \xrightarrow{j_{n,m}} BU(n) \times BU(m)
\xrightarrow{B\rho_{n,m}} BU(n+m), \]
where $\rho_{n,m}: U(n)\times U(m) \to U(n+m)$ is the usual 
block matrix inclusion. The isomorphism above appears as the   
factorization of $j_{n,m}^*$ through the positive degree part of the
image of $(B\rho_{n,m})^*$. 

One can describe $j_{n,m}$ as the composite map
\[ j_{n,m} : \Gr n {n+m} \to \Gr n {n+m} \times \Gr m {n+m} \to
\Gr n \infty \times \Gr m \infty , \]
where the first map is given by $V\mapsto (V,V^\perp )$.
Using the fact that the pullback of the canonical bundle $\gamma_m$ 
along the map $\Gr n {n+m} \to \Gr m {n+m}$, $V\mapsto V^\perp$
is the orthogonal complement $\gamma_n^\perp$, one sees that the 
Chern classes maps as follows:
\[ j_{n,m}^* (c_i (\gamma_n (\CC^\infty ))) = c_i(\gamma_n) \quad , \quad
j_{n,m}^* (c_i (\gamma_m (\CC^\infty ))) = c_i(\gamma_n^\perp ). \]

Put $c_i=c_i(\gamma_n )$ and $\bar c_j = c_j(\gamma_n^\perp )$.
By the dimension of the bundles we have that $c_i=0$ for $i>n$ and 
$\bar c_j=0$ for $j>m$. Since 
$\gamma_n \oplus \gamma_n^\perp \cong \epsilon^{n+m}$
we also have that $\sum_{i+j=k} c_i\bar c_j=0$ for $k>0$. 
This gives us a quotient map into the cohomology of the Grassman 
manifold, which by (\ref{eq:Borel}) is an isomorphism
\[
H^*(\Gr n {n+m} ;\ZZ ) \cong
\ZZ [c_i,\bar c_j | i,j>0] / 
(c_i|i>n)+(\bar c_j|j>m)+(\sum_{i+j=k}c_i\bar c_j |k>0).
\]

In the special case $n=2$ and $m=r-1$ we have that
\begin{align*}
& \ZZ[c_i,\bar c_j|i,j>0]/(c_i|i>2)+(\sum_{i+j=k} c_i\bar c_j |k>0)
\cong  \\
&\ZZ[c_1,c_2,\bar c_1,\bar c_2, \dots]/ ( \bar c_k-\phi_k(c_1,c_2)|k>0)
\end{align*}
Dividing by $(\bar c_j|j>r-1)$ we see that
\[
H^*(\Gr 2 {r+1} ;\ZZ ) \cong
\ZZ[c_1,c_2]/ (\phi_i| i\geq r).
\]
However, it follows from the inductive definition that
the $\phi_i$ is contained in the ideal generated by
$\phi_r,\phi_{r+1}$ for $i\geq r$ and this finishes the proof.
\end{proof}

We now turn to the projective bundle over the Grassmann
space. 

\begin{theorem}
\label{th:ProjectiveCohomology}
Let $\pi : \PP (\gamma_2 ) \to \Gr 2 {r+1}$ be the projective
bundle of the canonical bundle $\gamma_2$. There is an isomorphism
\[ H^*(\PP (\gamma_2 );\ZZ ) \cong 
\ZZ [x_1,x_2] / (Q_r,Q_{r+1}) , \]
where $x_1$ and $x_2$ have degree $2$ and 
\[ Q_k(x_1,x_2) = \sum_{i=0}^k x_1^ix_2^{k-i} \in \ZZ [x_1,x_2]. \]
We also have that
$\pi^*(c_1(\gamma_2 ))=x_1+x_2$ and $\pi^*(c_2(\gamma_2 ))=x_1x_2$.
\end{theorem}

\begin{proof}
We use \cite{Husemoller}, 17.2.5 and 17.2.6. 
Let $\lambda \to \PP (\gamma_2 )$ be the sub bundle of the
pullback $\pi_*(\gamma_2 )$ defined such that a point in the 
total space of $\lambda$ is a pair 
$(\{ V_1\subset V_2 \subset \CC^{r+1} \}, v)$ where the complex
dimension of $V_i$ is $i$ and $v\in V_1$. Let $\bar \lambda$ be
the conjugate bundle of $\lambda$. Then, we have
\[ 
H^*(\PP (\gamma_2 );\ZZ ) \cong
H^*(\Gr 2 {r+1} ;\ZZ )[c_1(\bar \lambda )] /
(c_1(\bar \lambda )^2 + c_1(\gamma_2 )c_1(\bar \lambda ) +
c_2(\gamma_2 )).
\]

Combining this with theorem \ref{th:GrassmannCohomology} we see that 
$H^*(\PP (\gamma_2 );\ZZ )$ is generated by the three classes 
$\pi^*(c_1(\gamma_2 ))$, $\pi^*(c_2(\gamma_2 ))$ and $c_1(\lambda )$ 
with the three relations
\begin{align*}
& \phi_r \big( \pi^* (c_1(\gamma_2 )),\pi^* (c_2(\gamma_2 ))\big) =0, \\
& \phi_{r+1} \big( \pi^* (c_1(\gamma_2 )),\pi^* (c_2(\gamma_2 ))\big) =0, \\
& \pi^* (c_2(\gamma_2 ))=
c_1(\lambda )\pi^* (c_1(\gamma_2 ))-c_1(\lambda )^2.
\end{align*}

Define $x_1=c_1(\lambda )$ and 
$x_2=\pi^* (c_1(\gamma_2 ))-c_1(\lambda )$.
Using the last of the above equations to eliminate 
$\pi^* (c_2(\gamma_2 ))$, we get that $H^*(\PP (\gamma_2 );\ZZ )$ 
is generated by the classes $x_1$ and $x_2$ subject to the relations 
we get by substituting $\pi^* (c_1(\gamma_2 ))$ and $\pi^* (c_2(\gamma_2 ))$ 
expressed by $x_1$ and $x_2$ into $\phi_r$ and $\phi_{r+1}$.
Note that $\pi^* (c_1(\gamma_2 ))=x_1+x_2$ and that
\[
\pi^* c_2(\gamma_2 )=
c_1(\lambda )\pi^* (c_1(\gamma_2 ))-c_1(\lambda )^2 =
x_1(x_1+x_2)-x_1^2 = x_1x_2.
\]

So we put $Q_r(x_1,x_2)=(-1)^r\phi_r (x_1+x_2,x_1x_2)$.
The two relations are polynomials $Q_r$, $Q_{r+1}$
in $x_1$ and $x_2$. The polynomials $Q_i$ are given inductively
by substituting into the inductive definition of $\phi_i$.
The inductive formula becomes
\[
Q_0 = 1, \quad Q_1 = x_1+x_2, \quad
Q_i = (x_1+x_2)Q_{i-1}-x_1x_2Q_{i-2}.
\]
It is easy to check that the polynomials
\[
Q_i(x_1,x_2)=\sum_{j=0}^i{x_1^jx_2^{i-j}}
=(x_1^{i+1}-x_2^{i+1})/(x_1-x_2)
\]
satisfy this inductive definition. 
\end{proof}

Because of lemma \ref{lemma:ProjectiveBundle} we have the
following result:

\begin{corollary}
\label{cor:GeoCohomology}
There is an isomorphism 
$H^*(\SG r;\ZZ)\cong\ZZ[x_1,x_2]/\{Q_r,Q_{r+1}\}$. 
The induced of the stabilization map $\SG r \hookrightarrow \SG {r+1}$ maps 
$x_1$ and $x_2$ in $H^*(\SG {r+1} ;\ZZ )$ to the classes with the 
same names in $H^*(\SG r ;\ZZ )$. 
\end{corollary}

\begin{proof}
The statement about the stability of the classes $x_1$ and $x_2$ 
follows from the fact that they are defined using Chern classes 
of the bundles $\gamma_2$ and $\lambda$. These bundles pull back 
to bundles with the same names.
\end{proof}

We note that it is very easy to check whether 
a polynomial is in the ideal generated by $Q_r$ and $Q_{r+1}$.

\begin{lemma}
\label{le:Qcheck}
Let $P=\sum_{i=0}^m p_ix_1^ix_ 2^{m-i}$ be a homogeneous polynomial of 
degree $m$. Then $P$ is contained in the ideal $I=(Q_r,Q_{r+1})$
if and only if $p_i=p_j$ for all $i,j$ such that 
$m-r\leq i\leq r$ and $m-r\leq j\leq r$.
\end{lemma}

\begin{proof}
Since $x_1^{r+1}=Q_{r+1}-x_1Q_r$, and similarly for $x_2^r$,
the monomials $x_1^ix_2^{m-i}$ are contained in $I$ if
$0\leq i\leq m-r-1$ or $r+1\leq i\leq m$. It follows that
$P\in I$ if and only if 
$\sum_{m-r \leq i\leq r}{p_ix_1^ix_ 2^{m-i}}\in I$.
If $p_{m-r}=p_{m-r+1}=\dots =p_r$, we see that
$P$ is congruent to $Q_m$ modulo $I$, and $P\in I$.

To see that this condition also is necessary, 
let $J\subset I$ be the ideal generated by $x_1^{r+1}$
and $x_2^{r+1}$. Then $x_1Q_k=x_2Q_k=Q_{k+1} \mod J$, 
so $I$ is generated as an abelian group  by $J$ 
together with $Q_k$ for $k\geq r$. So for any homogeneous
polynomial $P\in I$ of degree $m$, there is a 
$\lambda \in \ZZ$ such that $P-\lambda Q_m\in J$. 
This completes the proof.
\end{proof}

Next we consider the space of parametrized geodesics
$S(\tau(\CP^r))$. We consider two fibration, namely the
middle horizontal fibration in  (\ref{eq:fiberdiagram})
and the obvious spherical fibration: 
\[ SO(3) \to S(\tau (\CP^r )) \xrightarrow{q} \Gr 2 {r+1}, \quad  \quad
S^{2r-1} \to S(\tau (\CP^r )) \xrightarrow{\eta} \CP^r .\]
Let $\gamma_1$ be the canonical bundle over $\CP^r$, so that 
\[ 
H^*(\CP^r; \ZZ )=\ZZ [c_1(\gamma_1 )]/(c_1(\gamma_1 )^{r+1}).
\]
We have the following result (compare with the remark in the
introduction of \cite{AGMP}):

\begin{lemma}
\label{lemma:TangentCohomology}
Put $t=\eta^* (c_1(\gamma_1 ))$. There is a class 
$\bar \sigma \in H^{2r+1}(S(\tau (\CP^r ));\ZZ )$ such that
\[
H^*(S(\tau (\CP^r ));\ZZ )=\ZZ [t,\bar \sigma ]/
(t^{r+1}, \medspace \bar \sigma^2 , \medspace (r+1)t^r, 
\medspace t^r\bar\sigma ).
\]
Furthermore, $q^*(c_1(\gamma_2 ))=2t$
and $q^* (c_2(\gamma_2 ))=t^2$.
\end{lemma}

\begin{proof}
We consider the Serre spectral sequence for the spherical 
fibration $\eta$:
\[
E_2^{**}=\ZZ [c_1(\gamma_1 ),\sigma ]/
(c_1(\gamma_1 )^{r+1}, \medspace \sigma^2 )\Rightarrow 
H^*(S(\tau (\CP^r ));\ZZ ),
\]
where $\deg (\sigma )=2r-1$. For dimensional reasons, the only possibly
non-trivial differential is $d_{2r}(\sigma )$ which is given by the
Euler characteristic of $\CP^r$.
$d_{2r}(\sigma )= e(\eta ) = (r+1)c_1(\gamma_1 )^r$.
For this, cf \cite{Milnor} corollary 11.12 and theorem 12.2.
Let $\bar \sigma$ be a class representing
$c_1(\gamma_1 )\sigma$. There cannot be either further differentials 
nor extensions for dimensional reasons. This finishes the cohomology 
computation.  

\textit{Claim:} There is a bundle isomorphism
$q^*(\gamma_2 )\cong \eta^* (\gamma_1) \oplus \eta^*(\gamma_1 )$. \\
To prove the claim, we use the diffeomorphism $\phi$ from
lemma \ref{lemma:UnitTangents}. We have a commutative diagram

\[ \xymatrix@C=1.5 cm{
& \CP^{r} \\
\Fr 2 {r+1} \ar[r]^\phi_\cong \ar[ru]^{pr_1} \ar[rd]_p & 
S(\tau (\CP^r )) \ar[u]_\eta \ar[d]^q \\
& \Gr 2 {r+1},
} \]
where $p$ is the canonical projection and $pr_1$ is given by
projection on the first factor. We also have the projection on 
the second factor $pr_2$ and 
$pr_1^*(\gamma_1 ) \cong pr_2^*(\gamma_1 )$. 
So it suffices to show that 
$p^*(\gamma_2 ) \cong pr_1^*(\gamma_1 ) \oplus pr_2^*(\gamma_1 )$.
But this isomorphism is obvious.

We can now compute the total Chern class as follows:
\[
1+c_1(q^*(\gamma_2 ))+c_2(q^*(\gamma_2 ))=
(1+c_1(\eta^*(\gamma_1 )))^2=
1+2c_1(\eta^*(\gamma_1))+(c_1(\eta^*(\gamma_1 )))^2.
\]
The final statement of the lemma  follows from the  
naturality of the Chern classes.
\end{proof}

\begin{remark} 
\label{SCohomologyModP}
Similar to the integral computation we can compute 
cohomology with coefficients in $\FF_p = \ZZ/p$ for $p$ prime.
Let $x$ be the mod $p$ reduction of $t$.
If $p\mid (r+1)$, there is a class
$\sigma \in H^{2r-1}(S(\tau (\CP^r ));\FF_p )$ such that
\[
H^*(S(\tau (\CP^r ));\FF_p )\cong
\FF_p [x,\sigma ]/(x^{r+1},\medspace \sigma^2 ).
\]
Similarly, if $p\nmid (r+1)$, there is a class 
$\bar \sigma \in H^{2r+1}(S(\tau (\CP^r));\FF_p )$ such that
\[
H^*(S(\tau (\CP^r ));\FF_p )\cong
\FF_p [x,\bar \sigma ]/(x^r, \medspace \bar \sigma^2 ).
\]
\end{remark}

\begin{corollary}
\label{cor:QuotientMap}
The composition $S(\tau (\CP^r )) \xrightarrow{q^\prime} \SG r 
\xrightarrow{\psi} \PP (\gamma_2)$ of the quotient map with the
diffeomorphism satisfies
\[ (\psi \circ q^\prime )^*(x_1)= (\psi \circ q^\prime )^*(x_2)=
\eta^* (c_1(\gamma_1 )). \]
The kernel of the map  $(\psi\circ q^\prime)^*$ is a copy of the
integers, generated by $x_1-x_2$. 
\end{corollary}

\begin{proof}
Put $f=\psi \circ q^\prime$ and consider the composition
\[
S(\tau (\CP^r ))\xrightarrow{f} \PP (\gamma_2 )
\xrightarrow{\pi} \Gr 2 {r+1} ,
\]
which equals the map $q$ from (\ref{eq:fiberdiagram}).
By lemma \ref{lemma:TangentCohomology} we have that
\[ (\pi \circ f)^*(c_1(\gamma_2 ))=2t \quad , \quad
(\pi \circ f)^*(c_2(\gamma_2 ))=t^2, \] 
where $t=\eta^* (c_1(\gamma_1 ))$.
According to \ref{th:ProjectiveCohomology} we also have
that $\pi^* (c_1(\gamma_2 ))=x_1+x_2$ and 
$\pi^* (c_2(\gamma_2 ))=x_1x_2$. So we conclude that
\[
f^*(x_1)+f^*(x_2) = 2t \quad , \quad  
f^*(x_1)f^*(x_2) = t^2.
\]

If $r\geq 3$ one has 
$H^i(S(\tau(\CP^r));\ZZ)\cong \ZZ$ for $i=2,4$ generated by 
$t$ and $t^2$ respectively. So these equations have the unique 
solution $f^*(x_1)=f^*(x_2)=t$ which proves the corollary in this
case. If $i\leq 2$, we use the standard inclusion $\CC^i \subset \CC^3$
together with naturality to get the desired result.
\end{proof}

\section{Borel cohomology of spaces of geodesics}
\label{sec:Orbits}

We consider the space of parametrized geodesics $S(\tau (\CP^r ))$
with the free ${\TT}$-action. Let $S(\tau (\CP^r ))^{(n)}$ denote
$S(\tau (\CP^r ))$ where we have twisted the ${\TT}$-action 
by composing it with the $n$th power map $\lambda_n :{\TT} \to {\TT}$;
$z\mapsto z^n$. Write $C_n \subseteq {\TT}$ for the cyclic group of
order $n$. The map $\lambda_n$ passes to an isomorphism ${\TT}/C_n \to {\TT}$
with inverse $R_n$ given by sending $z$ to an $n$th root of $z$.
We write $BC_n$ for $E{\TT}/C_n$ with ${\TT}$-action through 
$R_n :{\TT}\to {\TT}/C_n$. 

The action of $C_n$ on $S(\tau (\CP^r ))^{(n)}$ is trivial. This makes it 
possible to consider $S(\tau (\CP^r ))^{(n)}$ as a ${\TT}/C_n$ space. 
We get homeomorphisms as follows:
\[ 
E{\TT}\times_{{\TT}} S(\tau (\CP^r ))^{(n)} \cong
E{\TT}/C_n \times_{{\TT}/C_n} S(\tau (\CP^r ))^{(n)} \cong
BC_n \times_{{\TT}} S(\tau (\CP^r )).
\]
We can write ${\TT}\to S(\tau (\CP^r )) \to \SG r$ as a pullback
of the universal ${\TT}$-bundle ${\TT}\to E{\TT}\to B{\TT}$ along a map 
$f:\SG r \to B{\TT}$. By this pullback diagram we get a map of
associated fibration sequences:

\begin{equation}
\label{eq:CpAction}
\xymatrix@C=2 cm{
BC_n \ar[r] \ar@{=}[d] 
& BC_n\times_{{\TT}} S(\tau (\CP^r )) \ar[r]^-{\rho} \ar[d] 
& \SG r \ar[d]^f \\
BC_n \ar[r] & B{\TT} \ar[r] & B{\TT}.
}
\end{equation}
Note that we have used the homotopy equivalence
$pr_1:BC_n\times_{{\TT}} E{\TT} \xrightarrow{\simeq} B{\TT}$ 
in the middle of the lower part of the diagram.
We can now compute the cohomology of the homotopy orbit spaces.

\begin{theorem}
\label{th:CohomologyOrbits}
For any prime $p$ one has
\[
H^*(E{\TT}\times_{{\TT}} S(\tau (\CP^r ))^{(n)})\cong
\begin{cases}
\FF_p [x_1,x_2]/(Q_r,Q_{r+1}), 
& p\nmid n, \\
\FF_p [u,x,\sigma ]/(x^{r+1},\sigma^2 ), 
& p\mid n, \medspace p \mid (r+1), \\
\FF_p [u,x,\bar \sigma ]/(x^r,\bar \sigma^2 ) ,
& p\mid n, \medspace p\nmid (r+1),
\end{cases}
\]
where $u$, $x$, $x_1$, $x_2$ have degree $2$ and $\deg (\sigma )=2r-1$,
$\deg (\bar \sigma )=2r+1$. 
\end{theorem}

\begin{proof}
If $p$ does not divide $n$, then the mod $p$ cohomology of $BC_n$
equals that of a point, and by an obvious spectral sequence argument 
the projection map $\rho$ induces an isomorphism in $\FF_p$
cohomology. The result follows by corollary \ref{cor:GeoCohomology}. 

Now assume that $p$ does divide $n$. Then one has
$H^*(BC_n)\cong \FF_p[v,w]/I_{n,p}$
where $I_{n,p}$ is the ideal defined by 
$I_{n,p}=(v^2-w)$ if $p=2$ and $4\nmid n$
and $I_{n,p}=(v^2)$ otherwise. The degrees are $\deg (v)=1$ and
$\deg (w)=2$.

The two fibrations of diagram (\ref{eq:CpAction}) each give a spectral 
sequences, and the vertical maps induce a map of spectral sequences. Let us denote the 
spectral sequence derived from the lower row of the diagram 
by $E_r(I)$, the one derived from the upper row by
$E_r(II)$ and the map by $f^*: E_r(I) \to E_r(II)$.

We have that $E_2(I)=\FF_p [u]\otimes \FF_p [v,w]/I_{n,p}$
where the only non-trivial differential is $d_2$, which 
is determined by $d_2w=0$, $d_2v=u$ and the product structure.
This follows since the inclusion of the fiber 
$BC_n \to B{\TT}$ is given by the inclusion $C_n\subseteq {\TT}$ and
the spectral sequence converges to $H^*(B{\TT})$ 
such that $E_\infty (I)= \FF_p [w]$.

We compute the $E_2$ page of the other spectral sequence:
\[
E_2(II) \cong H^*(\SG r)\otimes \FF_p [v,w]/I_{n,p}
\Rightarrow
H^*(BC_n\times_{{\TT}} S(\tau (\CP^r ))).
\]
The class $w$ is a permanent cycle, since it is the image of a 
permanent cycle in the spectral sequence $E_r(I)$.
The whole spectral sequence is generated by
the permanent cycle $w$ together with 
the classes in $E_2^{*,0}(II)$ and $E_2^{*,1}(II)$.

But this implies, for formal reasons, that the only possible
non-trivial differential is $d_2$. Using the product structure we 
also see that this non-trivial differential is determined by
$d_2(v)$, which by naturality equals the image $f^*(d_2(v))=f^*(u)$.

We claim that $f^*(u)=\lambda (x_1-x_2)$
for some $\lambda \in \FF_p \setminus \{ 0 \}$.
To see this, we consider the fibration sequence
\[ \xymatrix{
S(\tau (\CP^r ))\ar[r] &\SG r \ar[r]^-f & B{\TT}.
} \]
We have already investigated the involved spaces.
According to lemma \ref{lemma:TangentCohomology}
and corollary \ref{cor:GeoCohomology}, the corresponding 
spectral sequence for integral cohomology has the form 
\[
E_2=\ZZ [u] \otimes \ZZ [t,\bar \sigma]/J 
\Rightarrow \ZZ [x_1,x_2]/(Q_r,Q_{r+1}),
\]
where $J=(t^{r+1},\bar \sigma^2, (r+1)t^r,t^r\bar \sigma )$. 
In the notation, we do not distinguish between the classes 
$x_1,x_2\in H^2(\PP (\gamma_2);\ZZ)$ and their pull--backs under
$\psi\colon \SG r \xrightarrow{\psi}  \PP (\gamma_2)$

Consider the case $r\geq 2$ first. We have a short exact sequence
\[ \xymatrix{
0\ar[r] & H^2(B{\TT};\ZZ )\ar[r]^-{f^*} &
H^2(\SG r ;\ZZ )\ar[r]^-{(q^\prime )^* } & 
H^2(S(\tau (\CP^r ));\ZZ ) \ar[r] & 0.\\
} \]
According to corollary \ref{cor:QuotientMap} 
the kernel of  of $(q^\prime )^*$ is the free group generated 
by $x_1-x_2$.
It follows that $f^*(u)=\pm(x_1-x_2)$. Using naturality on the canonical 
inclusion $\SG 1 \subset \SG 2$, which is a ${\TT}$-map,
we see that this formula is also true for $r=1$. So we have
proved the claim.

We return to cohomology with  $\FF_p$ coefficients.
Let $K_r$ and $C_r$ be the kernel and the cokernel of 
multiplication by the element $x_1-x_2$ on 
$H^*(\SG r)$. Then 
\[
E_3(II)=(\FF_p [w] \otimes C_r)\oplus (\FF_p [w] \otimes vK_r), 
\] 
and the spectral sequence collapses from the $E_3$-page.

By theorem \ref{th:ProjectiveCohomology} and the equation 
$Q_k(x_1,x_1)=(k+1)x_1^k$ we see that
\[
C_r=
\begin{cases}
\FF_p [x_1]/x_1^{r+1}, & p \mid (r+1),\\
\FF_p [x_1]/x_1^{r}, & p \nmid (r+1).
\end{cases}
\]
The dimension of the cokernel agrees with the dimension of the kernel. 
So the kernel of multiplication by 
$x_1-x_2$ is a vector space over $\FF_p$ of dimension 
$r+1$ if $p\mid (r+1)$, and dimension $r$ if $p\nmid (r+1)$.

We need more precise information about the kernel,
since we want to compute the multiplicative structure
of the spectral sequence $E_r(II)$.
To determine the kernel, it is enough to exhibit as many linearly
independent elements in the kernel as it's known dimension.
So it suffices to find $r+1$ non-trivial elements in pairwise different 
degrees if $p\mid (r+1)$, and $r$ non trivial elements in pairwise
different degrees if $p\nmid (r+1)$. 

Consider the following homogeneous polynomial 
\[
a_k=x_1^k\sum_{i=0}^{r-1} (i+1) x_1^ix_2^{r-1-i}.
\]
A computation shows that
$(x_2-x_1)a_k=x_1^k(Q_r(x_1,x_2)-(r+1)x_1^r)$.
So if $p\mid (r+1)$, the elements $a_k$ for $k\geq 0$
are all in the kernel of multiplication by $x_1-x_2$. 
If $p\nmid (r+1)$ then $a_k$ is in the kernel if
$k\geq 1$, since $x_1^{r+1}=Q_{r+1}-x_2Q_r$.

To show that the kernel is generated by these classes, we
have to check that each $a_k$ is non-trivial in 
$\FF_p [x_1,x_2]/(Q_r,Q_{r+1})$, as long as $k\leq r$. 
We use lemma \ref{le:Qcheck} to do so. In $a_k$, the coefficient 
of $x_1^{k-1}x_2^r$ is $0$ and the coefficient of
$x_1^kx_2^{r-1}$ is $1$. Furthermore, both $k-1$ and $k$ lies
between $(k+r-1)-r=k-1$ and $r$ so by the lemma we conclude
that $a_k \notin (Q_r,Q_{r+1})$ for $k\leq r$ as desired.

It follows that if $p\mid (r+1)$ then the $r+1$ elements 
$a_0,a_1,\dots ,a_r$ form a basis for the kernel and 
if $p\nmid (r+1)$ then the $r$ elements 
$a_1,a_2,\dots ,a_r$ form a basis for the kernel.

By the explicit formula for $a_k$, we see that the kernel has
basis $a_0,x_1a_0,\dots ,x_1^ra_0$ when $p\mid (r+1)$ and
$a_1,x_1a_1,\dots ,x_1^{r-1}a_1$ when $p\nmid (r+1)$.
Let $\sigma$ represent $a_0$ when $p\mid (r+1)$ and 
and let $\bar \sigma$ represent $a_1$ when $p\nmid (r+1)$.

We obtain the $E_3$ term 
\[
E_3(II) = 
\begin{cases}
\FF_p [w,x_1,\sigma ]/(x_1^{r+1},\sigma^2 ), & p\mid (r+1), \\
\FF_p [w,x_1,\bar \sigma]/(x_1^r,\bar \sigma^2 ), & p\nmid (r+1).
\end{cases}
\]
We already noted that there can be no nontrivial differentials beyond 
the second one, so $E_3(II)=E_\infty (II)$.
\end{proof}

\begin{corollary}
\label{cor:SerreCollapse}
If $p$ divides $n$, then the mod $p$ cohomology Serre spectral 
sequence for the fibration
\[ 
\xymatrix@C=2 cm{
S(\tau (\CP^r ))^{(n)} \ar[r]^-{i} 
& E{\TT}\times_{{\TT}} S(\tau (\CP^r ))^{(n)} \ar[r]^-{pr_1} 
& B{\TT} }
\]
collapses at the $E_2$ page. 
If $p$ does not divide $n$, the inclusion $i$ of the fiber induces a 
surjection in even degrees
\[ i^*: H^{2*}(E{\TT}\times_{{\TT}} S(\tau (\CP^r ))^{(n)}) 
\twoheadrightarrow H^{2*}(S(\tau (\CP^r ))). \]
\end{corollary}

\begin{proof}
Assume that $p$ divides $n$. We can compute the $E_2$ page by remark 
\ref{SCohomologyModP} and we have just computed the cohomology of the 
total space. It turns out that the $E_2$ page is isomorphic to the 
cohomology of the total space. So there is not room for any differentials.

Assume that $p$ does not divide $n$. We only have to check that the 
classes $x^j$ are in the image of $i^*$. But by the homeomorphisms we
used for computing the cohomology of the Borel construction and
corollary \ref{cor:QuotientMap} we have that $i^*(x_1)=x$ and $i^*(x_2)=x$ 
so the result follows.
\end{proof}

\section{Equivariant vector bundles}
\label{sec:Bundles}
 
In this section we collect results about the homotopy theory of 
$\TT$-vector bundles which we will need later.
For greater clarity, the level of generality will be greater
than strictly necessary.  
 
We are interested in the cohomology of the Borel construction on 
Thom spaces. Let $G$ be a compact Lie group (we will only need the 
case $G={\TT}$). Assume that $\xi =(p:E \to X)$ is a $G$-vector-bundle 
over $X$ in the sense of \cite[I.9]{tomD}. 

\begin{lemma}
\label{lemma:Thom}
The map $EG\times_G \medspace p :EG\times_G E \to EG\times_G X$ is
the projection of a vector bundle which we denote $\xi_{hG}$. 
There is a homeomorphism
$\Th (\xi_{hG} )\cong EG_+\wedge_G \Th (\xi )$. 
Moreover, if $G$ is connected, then $\xi_{hG}$ is oriented if and 
only if $\xi$ is oriented when viewed as an ordinary vector bundle.
\end{lemma}

\begin{proof}
We have that $EG\times X$ is a free $G$-space. If this $G$-space 
is also strongly locally trivial, then proposition \cite[I.9.4]{tomD} 
gives us that $\xi_{hG}$ is a vector bundle. So we must show that 
$EG\times X$ has a tube around all of its point \cite[page 46]{tomD}.

For $(e,x)\in EG\times X$ we have that $G/G_{(e,x)}=G$. 
The universal principal $G$-bundle $\pi :EG\to BG$ is locally 
trivial, so we have a neighborhood $V$ around $\pi (e)$ and 
a trivialization $\phi :\pi^{-1} (V) \xrightarrow{\cong} V\times G$.
We get a tube around $(e,x)$ as follows:

\[ \xymatrix@C=2 cm{
\pi^{-1} (V)\times X \ar[r]^-{pr_1} & \pi^{-1} (V) \ar[r]^-{\phi}
& V\times G \ar[r]^-{pr_1} & G.
} \]
Thus $\xi_{hG}$ is a vector bundle as stated.

We may assume that there is a Riemannian metric on $\xi$ which is 
invariant under the $G$ action. We get a Riemannian metric
on $\xi_{hG}$ such that
$D(\xi_{hG} )=EG\times_G D(\xi )$ and 
$S(\xi_{hG} )=EG\times_G S(\xi )$.
So we have
\begin{align*}
\Th (\xi_{hG} ) &= EG\times_G D(\xi )/EG \times_G S(\xi ) \cong
EG_+\wedge_G (D(\xi )/S(\xi ))\\
&= EG_+\wedge_{G} \Th( \xi ).
\end{align*}

Regarding orientability, consider the fibration 
$G\to EG\times X \xrightarrow{\pi} EG\times_G X$. It is the 
pullback of $G\to EG\to BG$ along $pr_1:EG\times_G X \to BG$.
Furthermore, $BG$ is 1-connected since $G$ is connected so we have 
trivial coefficients in the $E_2$ pages of the Serre 
spectral sequences for $\pi$. Thus  
$H^1(EG\times_G X;\FF_2 )\to H^1(X;\FF_2 )$ 
is injective. The vector bundle $\xi_{hG}$ is the pullback of 
$\xi$ along $X\to EG\times_G X$.
Since orientability is equivalent to the vanishing of 
the first Stiefel Whitney class the result follows.
\end{proof}

\begin{corollary}
\label{cor:ThomIso}
Let $\eta$ be an $n$-dimensional ${\TT}$-vector bundle
over the ${\TT}$-space $X$. Assume that $\eta$ is
oriented for $H^*(-;\ZZ/p)$. Then there is an isomorphism
\[
\widetilde H^*(\Th (\eta)_{h{\TT}} ;\FF_p ) \cong 
H^{*-n}(X_{h{\TT}} ;\FF_p ).
\]
\end{corollary}

We are going to need a special case of the localization theorem for
Borel cohomology. The set up for the theorem is the following. 
\begin{definition}
 A based space $X$ homotopy equivalent to a $\TT$ CW complex 
satisfies the localization finiteness condition if
$X$ contains only finitely many ${\TT}$ orbit types, and $X$
is finite dimensional.
\end{definition}

Let $Y$ be any ${\TT}$ space with a fixed base point.
The space $E{\TT}_+\wedge Y$ is equipped with a diagonal map
$\tilde \Delta :E{\TT}_+\wedge Y \to 
E{\TT}_+\wedge E{\TT}_+ \wedge Y$
given by $\tilde \Delta (s,y)=(s,s,y)$. After taking 
quotients with group actions, we obtain
\[
\Delta_Y :E{\TT}_+ \wedge_{{\TT}} Y \to 
(B{\TT}_+)\wedge (E{\TT}_+ \wedge_{{\TT}} Y).   
\] 
This map induces a product 
\[
H^*(B{\TT})\otimes \widetilde H^*(E{\TT}_+ \wedge_{{\TT}} Y)
\to \widetilde H^*(E{\TT}_+ \wedge_{{\TT}} Y),
\]
which makes $\widetilde H^*(E{\TT}_+ \wedge_{{\TT}} Y)$ 
into a module  over the ring $H^*(B{\TT})$. A ${\TT}$-equivariant map 
$f:Y_1\to Y_2$ induces a module map by naturality. This module 
structure will pervade the whole theory.

One thing we can use it for is to localize. We invert the generator,
to form the cohomology localized away from $u$ as follows:
\[
H^*(X_{h{\TT}})\left[ \frac 1 u \right]=
H^*(X_{h{\TT}})\otimes_{\FF_p[u]} \FF_p [u,u^{-1}].
\]

\begin{theorem}
\label{th:localization}  
If $X$ satisfies the localization finiteness condition, 
the inclusion $X^{C_p}\subseteq X$
induces an isomorphism of localized cohomology
\[ \xymatrix@C=2 cm{
H^*(X_{h{\TT}})\left[ \frac 1 u \right] \ar[r]^-{\cong} 
& H^*((X^{C_p})_{h{\TT}})\left[ \frac 1 u \right].
} \]
\end{theorem}

\begin{proof}
This is a special case of the localization theorem
\cite[III.4.2]{tomD}. The parameters of the localization theorem, 
as given by tom Dieck, are chosen as follows:
$G=\TT$, the cohomology theory $H^*(-)=H^*(-;\FF_p)$
and the set $S\subseteq H^*(BG)=H^*(B{\TT};\FF_p)$ is 
$S=\{ 1,u,u^2,\dots \}$. The subset $A\subseteq X$ is the empty set.
 
The statement of the localization theorem is to be interpreted as 
follows: The family ${\cal F} (S)$ of subgroups of $\TT$ (defined in
\emph{ibid} after III.3.1) is the family of subgroups $C_n\subseteq \TT$ 
with $n$ not divisible by $p$. The set $X({\cal F})\subseteq X$
defined in \emph{ibid} I.6.1, is the complement of the
set $X^{C_p}$ in $X$. But then $FX=X^{C_p}$, so that theorem
\cite[III.4.2]{tomD} indeed specializes to our theorem.
\end{proof}

The application we want is the following: Let $X$ be a 
space, satisfying the localization finiteness condition. 
Assume that the action is not effective, but that 
the subgroup $C_p\subseteq {\TT}$ acts trivially. On the other hand, 
we assume that all isotropy groups of points in $X$ are contained in 
$C_n\subseteq {\TT}$ for some fixed $n$.
Let $\xi$ be a ${\TT}$-vector bundle over $X$. We do not assume that 
the action of $C_p$ on $\xi$ is trivial. The fixed points of $C_p$ 
forms a ${\TT}$ subbundle $\xi^{C_p}\subseteq \xi$. 

\begin{theorem}
The following map induces an isomorphism on mod $p$ cohomology localized 
away from $u$:
\[
E{\TT}_+\wedge_{\TT} \Th( \xi^{C_p} ) \to 
E{\TT}_+\wedge_{{\TT}} \Th (\xi ).
\]
\end{theorem}

\begin{proof}
The only possible isotropy groups of $\Th (\xi )$ are
${\TT}$ itself (for the base point in the Thom space), and subgroups 
of $C_n$. In particular, there are finitely many orbit types.
Since $\Th (\xi^{C_p} )=\Th (\xi )^{C_p}$, the result now follows 
by applying theorem \ref{th:localization} on 
the $\TT$-space $\Th (\xi )$. 
\end{proof}

\section{The twisted action}
\label{sec:Twisting}
Let $X$ be a $\TT$-space with action map $\mu :\TT \times X \to X$.
We can twist this action by the power map 
$\lambda_n :\TT \to \TT$; $\lambda_n(z)=z^n$ 
and obtain another $\TT$-space $X^{(n)}$. The underlying spaces of 
$X$ and $X^{(n)}$ are equal, but the action map for $X^{(n)}$ is
$\mu_n: \TT \times X^{(n)}\to X^{(n)}$;
$\mu_n (z,x)=\mu (\lambda_n (z),x)$.

\begin{lemma}
\label{le:FunctFib}
Let $X$ be an ${\TT}$-space. We have a pullback of fibration 
sequences which is natural in $X$ as follows:
\[ \xymatrix@C=2 cm{
BC_n \ar[r] \ar[d] 
& E{\TT}\times_{\TT} X^{(n)} \ar[r]^{g_n} \ar[d] 
& E{\TT}\times_{\TT} X\ar[d] \\
BC_n \ar[r] & B{\TT} \ar[r]^{B\lambda_n} & B{\TT}.
} \]
\end{lemma}

\begin{proof}
It is convenient to consider models of $E{\TT}$ and $B{\TT}$ 
which are realizations of simplicial topological abelian groups.
So let $E{\TT}=|E{\TT}\simp|$ and $B{\TT}=|B{\TT}\simp|$, where 
$E\TT_q$ is the $q+1$-fold Cartesian product of $\TT$ by itself with
\[ 
d_i(z_0,\dots ,z_q)=(z_0,\dots ,\widehat{z_i},\dots , z_q), \quad
s_i(z_0,\dots ,z_q)=(z_0,\dots ,z_i,z_i,\dots , z_q)
\]
and $B\TT_q =E\TT_q /\diag_q (\TT )$. The hat means that the element is
left out.

We have simplicial maps $E(\lambda_n )\simp$ and $B(\lambda_n )\simp$
given by rising all elements in a tuple to the $n$th power.
Since the kernel of $B(\lambda_n)_q$ is exactly $(BC_n)_q$, the bottom
row in the diagram is the realization of a short exact sequence of 
simplicial abelian groups, which means that it is a fibration.
Since $E(\lambda_n )$ is a map over $B(\lambda_n)$ the right hand
square in the diagram commutes when we put 
$g_n([e,x])=[E(\lambda_n )(e),x]$. 
Note that this makes $g_n$ well-defined since 
$E(\lambda_n )(ez)=E(\lambda_n )(e)z^n$.
The map of the vertical fibers is the identity so the right hand
square is a vertical pullback. It follows that it is also a 
horizontal pullback. 
\end{proof}

Here is our main result on Borel cohomology of twisted actions. 
\begin{theorem}
\label{th:TwistedCollapse}
For any prime $p$ and $\TT$-space $X$ there is an isomorphism, 
natural in $X$:
\[ \xymatrix@C=2 cm{
H^*(X)\otimes \FF_p [u] \ar[r]^\cong  
& H^*(E{\TT}\times_\TT X^{(p)}).
} \]
\end{theorem}

To prove the theorem, we first show that there exists a natural 
subgroup ${\cal F}_1^*(X)$ of $H^*((X^{(p)})_{h\TT})$
with the property that the statement is true if we 
replace $H^*(X)$ by ${\cal F}_1^*(X)$. After we have 
constructed this subgroup, it is easy to identify it with 
$H^*(X)$.

The spectral sequence derived from the fibration of 
lemma \ref{le:FunctFib} defines a filtration  
\[
{\cal F}_0^*(X) \subseteq {\cal F}_1^*(X) \subseteq \dots \subseteq 
H^*(E{\TT}\times_{{\TT}} X^{(p)})
\]
such that 
${\cal F}_i^*(X)/{\cal F}_{i-1}^*(X) \cong 
\bigoplus_{*\geq 0} E_\infty^{i,*}$.
Each of these filtration subgroups is a functor on the category
of ${\TT}$-spaces. It is well known how to interpret the first
group:
\[ {\cal F}_0^*(X)=g_p^*(H^*(X_{h\TT})) \]

We need to consider the
next step in the filtration. This group does not have an
equally simple description.

\begin{definition}
The group of low elements 
${\cal F}_1^*(X)\subseteq H^*(E{\TT}\times_{{\TT}} X^{(p)})$
is the subgroup consisting of those classes which in the Serre 
spectral sequence represents elements in $E_\infty^{*,0}$ or 
$E_\infty^{*,1}$. 
\end{definition}

We define a reduced version of the group in the obvious way. 
If $X$ is a pointed ${\TT}$-space, we have inclusions of groups, 
natural on such spaces
\[
g_n^*(H^*(E{\TT}\times_{{\TT}} X,B{\TT}))=
{\cal F}_0^*(X,*)\subseteq {\cal F}_1^*(X,*)\subseteq
H^*(E{\TT}\times_{{\TT}} X^{(n)},B{\TT}). 
\]

Here is a general fact about the low elements corresponding to
$p$-fold twisting.

\begin{lemma} \label{LowBasis}
$H^*(X^{(p)}_{h\TT})$ is a free module over 
$H^*(B\TT)= \FF_p[u]$. The low elements ${\cal F}_1^*(X)$ 
form a basis for it as an $\FF_p[u]$ module.
\end{lemma}

\begin{proof}
The module structure on $H^*(X^{(p)}_{h\TT})$
is given by the map $pr_1:X^{(p)}_{h\TT}\to B{\TT}$.
The Serre spectral sequence associated to the upper fibration
sequence in lemma \ref{le:FunctFib} has the form
\[
E_2^{n,m}(X)\cong H^n(X_{h\TT})\otimes H^m(BC_p)
\Rightarrow H^{n+m}(X^{(p)}).
\]
Remember that $H^*(BC_p)\cong \FF_p[v,w]/I_{p}$, 
where $I_p$ is a principal ideal, generated by $v^2-w$ if $p=2$ and
by $v^2$ otherwise. 

The Serre spectral sequence of the lower fibration sequence in
the lemma has $E_2^{*,*}=\FF_p [u]\otimes \FF_p [v,w]/I_p$ and
converges towards $H^*(B\TT )$ so we must have $d_2(v)=u$ and 
$E_3(*)=E_\infty (*)=\FF_p [w]$.
Using this together with naturality in $X$ of the spectral sequence 
we get that for any $X$ the class $w\in E_2^{0,2}(X)$ is a permanent
cycle. Furthermore,
\[
E_2^{*,*}(X)\cong \FF_p[w]\otimes (E_2^{*,0}\oplus E_2^{*,1})
\] 
so it follows formally that 
$E_3^{*,*}(X)\cong \FF_p[w]\otimes (E_3^{*,0}\oplus E_3^{*,1}).$
But then it also follows formally that higher differentials in
this spectral sequence vanish, so that
\[
E_\infty^{*,*}(X)\cong 
\FF_p[w]\otimes (E_\infty^{*,0}\oplus E_\infty^{*,1}).
\]
   
The low elements are by definition a subspace of $H^*(X^{(p)}_{h\TT})$.
We can extend this inclusion to a unique map of $\FF_p [u]$-modules
\[
f:\FF_p [u]\otimes {\cal F}_1^*(X)\to H^*(X^{(p)}_{h\TT}).
\]
This map sends $u^i\otimes {\cal F}_1^*(X)$ to
${\cal F}_i^*(X)$, and the map induces an isomorphism on
the corresponding graded rings, by the above computation of
$E_\infty^{*,j}$, and because $u\in H^2(B\TT)$ represents 
$w\in E_\infty^{0,2}$. 
It follows that $f$ itself is an isomorphism.
\end{proof}

\begin{proof}[proof of theorem \ref{th:TwistedCollapse}]
We have to exhibit a natural isomorphism  
${\cal F}_1^*(X)\to H^*(X)$. From the spherical fibration
\[ \xymatrix@C=1cm{
\TT \ar[r] & E\TT \times X^{(p)} \ar[r]^i & E\TT \times_\TT X^{(p)}
} \]
we get a long exact Gysin sequence
\[ \xymatrix@C=1cm{
\ar[r]
& H^{*-2}(X^{(p)}_{h\TT}) \ar[r]^-{\cdot u} 
& H^*(X^{(p)}_{h\TT}) \ar[r]^-{i^*}
& H^*(X^{(p)}) \ar[r]
& H^{*-1}(X^{(p)}_{h\TT}) \ar[r]
&
} \]
(One can use the map $X\to *$ to see that the Euler class is $u$.)
Since $H^*(X^{(p)}_{h{\TT}})$ is a free $\FF_p [u]$ module, 
multiplication by $u$ is injective, and the long exact sequences 
breaks up into short exact sequences. So via lemma \ref{LowBasis}
we get the short exact sequence
\[ \xymatrix@C=1cm{
0 \ar[r]
& \FF_p [u]\otimes {\cal F}_i^*(X) \ar[r]^-{\cdot u} 
& \FF_p [u]\otimes {\cal F}_i^*(X) \ar[r]^-{i^*}
& H^*(X^{(p)}) \ar[r]
& 0.
} \]
Thus, $i^*$ factors through 
$\big( \FF_p [u]\otimes {\cal F}_i^*(X)\big) /
\image (\cdot u)= {\cal F}_i^*(X)$
and gives the desired isomorphism.
\end{proof}

\section{The $n$-fold iteration map}
\label{sec:Iteration}

Let $M$ be a compact Riemannian manifold. In this chapter we are 
examining the Hilbert manifold model of the free loop space.
We denote this Hilbert manifold by $LM$. An element in $LM$ is 
a curve $\gamma :\TT \to M$ of class $H^1$. The Hilbert manifold model 
is homotopy equivalent to the usual continuous mapping space 
model. The energy integral 
\[ E:LM \to \RR; \quad 
E(\gamma )=\int_\TT |\gamma^{\medspace \prime}(z)|^2dz  \]
is a smooth function on $LM$. It's critical points are the closed 
geodesics on $M$. We review the main results of Morse theory in
this setting. For details we refer to \cite{Klingenberg},
especially to chapter 1.

The power map $\lambda_n :\TT \to \TT$; $\lambda_n (z) =z^n$ gives us an
an $n$-fold iteration map
\[
\PS_n :LM\to LM; \quad \PS_n(\gamma )= \gamma \circ \lambda_n .
\]
Just like other versions of the iteration map, 
this map is not equivariant with respect to the ${\TT}$-action,
but becomes equivariant if we twist the action. That is,
the formula defines an equivariant map $\PS_n :(LM)^{(n)}\to LM$.

The map preserves energy up to a constant factor, 
$E(\PS_n(\gamma))=n^2E(\gamma)$.
That is, it induces an equivariant map of energy filtration
$\FC (a)=E^{-1}(-\infty,a)$ as follows:
\[
\PS_n :(\FC (a))^{(n)} \to  \FC (n^2a)
\]
We also get an equivariant map of quotients
\[
\PS_n : (\FC (a)/\FC ({a-\epsilon}))^{(n)} \to  
\FC({n^2a})/\FC({n^2a-n^2\epsilon}).
\]

Let $T_\gamma LM$ denote the tangent space of $LM$ at 
$\gamma$. It consists of the vector fields along $\gamma$ of 
class $H^1$, and is a real Hilbert space with inner product
\[
{\inp \xi \eta}_1= \int_{\TT} \big( \inp {\xi (z)} {\eta (z)} +
\inp {\nabla \xi (z)} {\nabla \eta (z)} \big) dz.
\]
Here $\nabla$ denotes the covariant differentiation 
along the curve defined by the Levi--Civita connection of the 
Riemannian manifold $M$. The inner product defines a Riemannian 
metric on the Hilbert manifold $LM$.
One technical problem is that the iteration map is not an (injective)
isometry with  respect to this metric. To study the properties of
$\PS_n$ 
it will be convenient to consider a modified inner product on 
$T_\gamma LM$. Let
\[
{\inp \xi \eta }_{c,1}=\int_{\TT} \big( \inp {\xi (z)} {\eta (z)} +
c\inp {\nabla \xi (z)} {\nabla \eta (z)} \big) dz.
\]
The differential of the iteration map is given by  
\[ D_\gamma (\PS_n ): T_\gamma LM \to T_{\PS_n (\gamma )}LM ;\quad
D_\gamma (\PS_n ) (\xi ) =\xi \circ \lambda_n .\]
So we have 
$\nabla (D_\gamma(\PS_n )\xi )= n\nabla (\xi ) \circ \lambda_n =
n D_\gamma (\PS_n )(\nabla \xi )$ and
\[
{\inp {D_\gamma (\PS_n )\xi} {D_\gamma (\PS_n )\eta }}_1 
= \int_{\TT} \inp {\xi (z^n)} {\eta (z^n)} +
n^2 \inp {\nabla (\xi )(z^n)} {\nabla (\eta )(z^n)} dz= 
{\inp \xi \eta }_{n^2,1}.
\]
This means that even if the iteration map does not preserve the 
inner product, it becomes an isometry if we replace the
inner product on the target by a suitably modified
inner product.

The metric on $M$ determines an exponential map 
$\exp_p\colon T_pM\to M$ for each $p\in M$. This induces an 
exponential map 
\[
\exp^\sim_\gamma : T_\gamma LM \to LM; \quad
\exp^\sim_\gamma (\xi )(t)=\exp_{\gamma (t)}(\xi (t)).
\]
This is not the exponential map derived from the 
Hilbert metric on $LM$, but it has the advantage that 
it is compatible with the iteration map in the sense
that the following diagram commutes:
\[
\xymatrix@C=2cm{
T_\gamma LM \ar[r]^{D_\gamma (\PS_n )} \ar[d]^{\exp^\sim_\gamma}
& T_{\PS_n (\gamma )} LM \ar[d]^{\exp^\sim_{\PS_n(\gamma )}} \\
LM \ar[r]^{\PS_n} & LM. \\ 
}
\] 
The exponential map is a diffeomorphism from a neighborhood of
$0$ in $T_\gamma LM$ to a neighborhood of $\gamma$ in $LM$.

Let $\gamma$ be a closed geodesic on $M$. This is a critical point 
for the energy function $E$. The Hessian of the energy function defines 
a quadratic form $D^2E(\gamma )$ on $T_\gamma LM$. This form determines 
a self adjoint operator $A_\gamma :T_\gamma LM \to T_\gamma LM$ by the 
equation $D^2E(\gamma )(\xi ,\eta )={\inp {A_\gamma \xi } \eta }_1$. 
The operator $A_\gamma$ is the sum of the identity with a compact
operator, so there are at most a finite number of negative
eigenvalues, each corresponding to a finite dimensional vector space 
of eigenvectors of $A_\gamma$. 

Obviously the operator $A_\gamma$ depends not just on
the Hessian of the energy integral, but also on the metric on $LM$. 
Since we are considering modifications of the inner product,
we also have to consider the corresponding modifications of the self adjoint operator. 
To emphasize this dependence, for a real vector space
$V$ with inner product ${\inp \medspace \medspace }_\alpha $ and with a
quadratic form $Q$, we define the operator
$A({\inp \medspace \medspace }_\alpha ,Q)$ by the property that 
$Q(\xi ,\eta )=
{\inp {A({\inp \medspace \medspace }_\alpha ,Q)\xi } \eta }_\alpha$.

Now let $K_a$ be the space of critical points of $E$ with energy 
level $a$. The negative bundle $\mu^-$ over $K_a$ is the vector
bundle whose fiber at $\gamma$ is the vector space spanned by the 
eigenvectors belonging to negative eigenvalues of $A_\gamma$. 
Similarly, $\mu^0$ and $\mu^+$ are the vector bundles with fibers
spanned by the eigenvectors corresponding to the eigenvalue $0$ and the 
positive eigenvalues respectively.

Using the modified inner product ${\inp \medspace \medspace}_{c,1}$ 
we obtain a modified negative bundle $\mu_c^-$ etc. This bundle is 
not the same as $\mu^-$. But on the other hand, the difference 
between these two bundles is not so dramatic, since the orthogonal 
projection from $\mu_c^-$ to $\mu^-$ with respect to the standard 
inner product defines an isomorphism of bundles. 

\begin{lemma}
The iteration map induces a bijection 
$\PS_n :K_a \to (K_{an^2})^{C_n}$. 
The group $C_n$ acts on the fiber $(\mu^- )_{P_n\gamma }$, and for 
any $\gamma \in K_a$ the differential of the iteration map induces 
an isomorphism of vector spaces 
\[
\xymatrix@C=2cm{
D_\gamma (\PS_n ): (\mu_{n^2}^- )_\gamma  \ar[r]^-{\cong}
& ((\mu^- )_{\PS_n \gamma })^{C_n}. }
\] 
\end{lemma}

\begin{proof}
The property of being a geodesic is a local property of a curve, 
so $\PS_n \gamma$ is a geodesic if and only if $\gamma$ is a geodesic. 
It follows that $\PS_n (K_a)\subseteq (K_{n^2a})^{C_n}$. 
On the other hand, if $\theta \in (K_{n^2a})^{C_n}$ we can write
$\theta =\PS_n \gamma$ for some $\gamma \in LM$.
But $\gamma$ has to be a geodesic, since $\theta$ is,
and $E(\gamma )=E(\PS_n \gamma )/n^2=a$, which means that
$\gamma \in K_{a}$. Thus, $\PS_n$ induces a bijection as stated.

Now we look at the tangent map. We first claim that this map induces 
an isomorphism
\[
\xymatrix@C=2cm{
D_\gamma (\PS_n ): T_\gamma LM \ar[r]^-{\cong}
& (T_{\PS_n \gamma }LM)^{C_n}.}
\]
This is clear, since $\PS_n \gamma$ is periodic of period $n$,
and any $n$-periodic vector field along $\PS_n \gamma$ is the
image under the differential of the iteration map of a vector
field on $\gamma$. 

We must show that this isomorphism restricts to an isomorphism as 
stated. The groups $C_n$ acts on $V=T_{\PS_n \gamma}LM$,
let us say that a generator $\sigma \in C_n$ acts by
\[ (\sigma \xi )(z)=\xi (z\zeta_n )\in T_{\PS_n \gamma (z\zeta_n )}M=
T_{\PS_n \gamma (z)}M; \quad \zeta_n =e^{2\pi i/n} . \] 
Let $N\in \CC [C_n]$ be the sum $\frac 1 n \sum_{k=0}^{n-1} \sigma^k$.
This $N$ acts as an idempotent on $V$, and so does $\Id -N$.
The inner product is manifestly invariant under the action on $C_n$, 
so a simple calculation shows that $N$ is self adjoint such that
$N$ and $\Id-N$ are orthogonal idempotents. Moreover,
the idempotents commute with the action of the group $C_n$. So $V$ splits as a 
$C_n$ representation into an orthogonal sum $NV\oplus (\Id -N)V$.

The energy function is also invariant under the group action,
so by the same argument as we just used for the inner product,
the quadratic form $(V,D^2_{\PS_n \gamma}(E)$ splits as a direct sum
$(NV,D^2_{\PS_n \gamma}(E))\oplus ((\Id-N)V,D^2_{\PS_n \gamma}(E))$. 

Let $A=A({\inp \medspace \medspace}_1 ,D^2_{\PS_n \gamma}(E))$.
Since both the inner product and the quadratic form split
as direct sums, the linear endomorphism $A:V\to V$ is a direct sum
of two self adjoint endomorphisms 
\[ A|_{NV}:NV\to NV, \quad A|_{(\Id-N)V}:(\Id-N)V\to (\Id-N)V. \] 
It follows that the subspace of $V$ generated by negative eigenvectors 
equals the direct sum of the negative eigenvector spaces of 
$A|_{NV}$ and $A|_{(\Id-N)V}$. 

We claim that $NV=V^{C_n}$ and $((\Id -N)V)^{C_n}=0$. This follows  
since $\sigma N=N$ so that $NV\subseteq V^{C_n}$ and if 
$\sigma (\Id -N)\xi =(\Id -N)\xi$ then $\sigma \xi =\xi$ so that
$(\Id -N)\xi =0$. 

To finish the proof, we have to show that 
$D_\gamma (\PS_n )(\mu^-_{n^2} )_\gamma$ equals
the negative eigenvector space of $A|_{NV}$. 

But $D_\gamma (\PS_n )$ is an isometry (with respect to the modified
inner product), and it preserves $D^2E$ up to multiplication by $n^2$.
So have a commutative diagram
\[
\xymatrix@C=2cm{
T_\gamma LM \ar[r]^{D_\gamma \PS_n} 
\ar[d]_{n^2A({\inp \medspace \medspace}_{n^2,1},D^2_\gamma E)}
& T_{\PS_n \gamma }LM \ar[d]^{A({\inp \medspace \medspace}_1 ,
D^2_{\PS_n \gamma }E)} \\
T_\gamma LM \ar[r]^{D_\gamma \PS_n } & T_{\PS_n \gamma} LM 
}
\]
by which we get the desired result.
\end{proof}

Assume that $N_a\subseteq LM$ is a nondegenerate
critical manifold of energy $a$. Also assume that
$\PS_p (N_a) \subseteq N_{p^2a}$ is a non-degenerate
critical manifold, and that there are no other critical
points at energy levels in the intervals
$(a-\epsilon,a)$ and $(p^2(a-\epsilon),p^2a)$. 

The Whitney sum $\mu=\mu^- \oplus \mu^+$ is the normal
bundle of the submanifold $N_a  \subseteq LM$ because we are
assuming that $N_a$ satisfies Bott's non-degeneracy condition.  

\begin{theorem}
\label{th:IterationMap}
The $p$-fold iteration map induces an isomorphism in cohomology 
localized away from $u$ as follows:
\[
\widetilde H^*((\FC (a)/\FC({a-\epsilon}))_{h\TT})
\left[ \frac 1 u \right] \cong 
\widetilde H^*((\FC ({p^2a})/\FC ({p^2(a-\epsilon)}))^{(p)}_{h{\TT}})
\left[ \frac 1 u \right] .\]
\end{theorem}

\begin{proof}
We have a commutative diagram (possibly after decreasing $\epsilon$)
\[
\xymatrix@C=1cm{
(D(\mu_{p^2}^- (N_a)),S(\mu_{p^2}^- (N_a)))\ar[r]^-{D\PS_p} 
\ar[d]^{\exp^\sim} 
& (D(\mu^- (N_{p^2a})),S(\mu^- (N_{p^2a})))\ar[d]^{\exp^\sim} \\
(\FC (a),\FC ({a-\epsilon}))\ar[r]^-{\PS_p}
& (\FC ({p^2a}),\FC ({p^2(a-\epsilon)})).
}
\]
According to theorem \ref{th:localization} the upper horizontal map 
induces an isomorphism in mod $p$ cohomology localized away from $u$.

The exponential map in the right column induces a $\TT$ homotopy 
equivalence on quotients by theorem \cite[2.4.10]{Klingenberg}. 
This is the Hilbert manifold version of the fundamental theorem of 
Morse theory. So to prove the theorem, we need to see that the left 
vertical map induces a ${\TT}$ equivariant homotopy equivalence on 
quotients.

But this is true for the same reason. The proof of  
theorem \cite[2.4.10]{Klingenberg} does not use the explicit
form of the metric on the Hilbert manifold $LM$, but only the 
nondegeneracy of $D^2E$. The statement of the variation of the 
theorem using $LM$ with the modified metric
${\inp \medspace \medspace }_{1/p^2,1}$ is exactly the statement
that the left vertical map is a homotopy equivalence.
\end{proof}

\begin{remark}
\label{re:IterationMap}
The theorem and its proof are also valid if $a$ is not a critical value.
In this case the interpretation of the statement is that if $\lambda$ is 
a critical value of the energy function, such that $\lambda/p^2$ is not a 
critical value, then for $\epsilon >0$ sufficiently small,
\[
\widetilde H^*((\FC ({\lambda })/\FC 
({\lambda -\epsilon}))^{(p)}_{h\TT})
\left[ \frac 1 u \right] =0.
\]
\end{remark}

\begin{remark}
Since $P_p$ is a homeomorphism whose image is the $C_p$ fixed points, 
we could try to use theorem \ref{th:localization} directly on
the inclusion $\FC ({p^2a})^{C_p}\subseteq \FC ({p^2a})$, cleverly
bypassing this entire section. The problem with this approach is that
you need a finite dimensionality condition for the localization statement 
\ref{th:localization} to be true. The only role of Morse theory in the 
proof of theorem \ref{th:IterationMap} is to reduce the 
infinite dimensional situation to a finite dimensional one.

It seems likely that this reduction can be done in greater
generality, and that the non-degeneracy condition in \ref{th:IterationMap} 
is far stronger than necessary. 
\end{remark}

\section{The Morse spectral sequences}
\label{sec:MSS}

Let $M$ be a compact Riemannian manifold as before. But in this 
section we are going to assume that the critical points of the 
energy function on $LM$ are collected on compact critical manifolds.
We also assume that each of these critical manifolds satisfy the 
Bott non-degeneracy condition. This is a strong assumption
on the metric of $M$, but for instance the symmetric spaces
satisfy this, according to \cite{Ziller}. 

We are considering the filtration induced by the levels of
the energy function. Let the critical values of the energy function be
$0=\lambda_0 <\lambda_1 <\dots $. There is a filtration of
$LM$ by $\FC (\lambda_i )=E^{-1}(-\infty ,\lambda_i )$.
This filtration is equivariant with respect to the action of
the circle, and induces spectral sequences of various forms
of cohomology. We call these spectral sequences Morse spectral 
sequences. The condition we impose on the metric of $M$ means that for
each $\lambda$ we have a finite dimensional critical manifold
$N(\lambda )$. The tangent bundle of $L M$ restricted to
$N(\lambda )$ splits ${\TT}$-equivariantly into a sum of three bundles.
\[
TLM|_{N(\lambda )}\cong
\mu^- (\lambda )\oplus \mu^0 (\lambda )\oplus \mu^+ (\lambda ).
\]
The standard Morse theory argument can be carried through equivariantly 
on the Hilbert manifold $LM$. This was done by Klingenberg. For an account
of this work see section \cite[2.4]{Klingenberg}, especially 
theorem 2.4.10. The statement of this theorem implies that we have 
an equivariant homotopy equivalence
\[
\FC (\lambda_n )/\FC (\lambda_{n-1} ) \simeq \Th (\mu^- (\lambda_n )).
\]

We will consider cohomology of the homotopy 
orbits $LM_{h\TT}$. We are also going to consider the twisted action. 
This is not because we have a particular interest in 
$H^*(LM_{h\TT}^{(p)})$ in itself. As we will see, 
these groups are easy to compute anyhow. The purpose of considering
them is rather to be able to use a comparison argument to
obtain information about the Morse spectral sequence
of $H^*(LM_{h\TT})$. The abstract situation is as follows:

\begin{theorem}
\label{th:ThreeSS}  
There are three spectral sequences 
\begin{align*}
& \{ E_r^{n,m}(\Morse)(LM) \} \Rightarrow 
H^{n+m}(LM), \\
& \{ E_r^{n,m}(\Morse)(LM_{h\TT}) \} \Rightarrow 
H^{n+m}(LM_{h\TT}), \\
& \{ E_r^{n,m}(\Morse)(LM^{(p)}_{h\TT}) \} \Rightarrow  
H^{n+m}(LM^{(p)}_{h\TT}). 
\end{align*}
The $E_1$ pages are given by
\[ \widetilde H^m(\Th (\mu^- (\lambda_n )))), \quad 
\widetilde H^m(\Th (\mu^- (\lambda_n)_{h\TT}))), \quad  
\widetilde H^m(\Th (\mu^- (\lambda_n )_{h\TT}^{(p)} ))) \]
respectively. We have a natural isomorphism of spectral sequences 
\[
E_*(\Morse)(LM_{h\TT }^{(p)} )\cong 
\FF_p [u]\otimes_{\FF_p} E_*(\Morse)(LM) 
\] 
as a module over $H^*(B{\TT })\cong \FF_p [u]$.
\end{theorem}

\begin{proof}
Consider the energy filtration 
$\FC (\lambda_0 )\subseteq \FC (\lambda_1 )\subseteq \dots$ 
with union $LM$. Because the filtration is equivariant, it induces 
filtrations of the spaces $LM$, $LM_{h\TT}$ and $LM^{(p)}_{h{\TT}}$. 
These filtrations give rise to three exact couples which gives
the three spectral sequences. The spectral sequences all converge 
(strongly). One can use \cite[Ch. 9,1.6]{Sp} or \cite[theorem 12.6]{Bo}
and the remark following \cite[theorem 7.1]{Bo} to see this.
 
The functorial isomorphism of $H(B{\TT} ;\FF_p )$-modules from 
theorem \ref{th:TwistedCollapse} imply the natural isomorphism of 
spectral sequences.
\end{proof}

We next have to discuss the $p$-fold iteration map.
If there exists a periodic geodesic of energy $\lambda$, its
$p$-fold iterate is a geodesic of energy $p^2\lambda$. So the
sequence $p^2\lambda_0 <p^2\lambda_1 <p^2\lambda_2 <\dots$
is a subsequence of the sequence
$\lambda_0 <\lambda_1 <\lambda_2 <\dots$.
The $p$-fold iteration map induces an equivariant map of filtrations:
\[ \xymatrix@C=1cm{
\FC (\lambda_0 )^{(p)} \ar[r] \ar[d] & 
\FC (\lambda_1 )^{(p)} \ar[r] \ar[d] & \dots \ar[r] & LM^{(p)} \ar[d] \\
\FC (p^2\lambda_0 ) \ar[r] & 
\FC (p^2\lambda_1 ) \ar[r] & \dots \ar[r] & LM. \\
} \]

\begin{lemma}
\label{le:ItLocalization}
Let $p^2\lambda_{i-1} =\lambda_j <\lambda_{j+1} <\dots 
<\lambda_k =p^2\lambda_i$ be the critical values of the energy function 
in the interval $[ p^2\lambda_{i-1} , p^2\lambda_i ]$. Then
\[ 
H^*(\FC (\lambda_\ell )_{h\TT} , \FC (\lambda_{\ell-1})_{h\TT} ;\FF_p )
\left[ \frac 1 u \right] = 
\begin{cases}
0 & , \ell \neq k, \\
H^*(\FC (\lambda_i ),\FC (\lambda_{i-1} );\FF_p )\otimes \FF_p [u,u^{-1}]
& , \ell=k.
\end{cases}
\]
The isomorphism is induced by the iteration map followed by the
isomorphism of theorem \ref{th:TwistedCollapse}.
\end{lemma}

\begin{proof}
This follows from theorem \ref{th:IterationMap} and remark 
\ref{re:IterationMap}. 
\end{proof}

\begin{theorem}
\label{th:relabel}
Up to a re-indexing of the columns in the spectral sequence
the $p$-fold iteration map induces a natural isomorphisms 
of spectral sequences 
\[
E_*(\Morse )(LM_{h\TT}) \left[ \frac 1 u \right] \cong
E_*(\Morse )(LM) \otimes \FF_p [u,u^ {-1}].
\]
\end{theorem}

\begin{proof}
We have to explain the re-indexing.
Let $j(i)$ be the non-negative number such that
$p^2\lambda_i=\lambda_{j(i)}$, and $k(i)$ the largest
number $s\leq i$ such that $s=j(t)$ for some $t$. So by
definition, $k(i)\leq i$. 

Besides the filtration $\FC_i:=\FC (\lambda_i)$ we consider two 
derived filtrations of $LM$;
\[ \FC_i^\prime := \FC_{j(i)} \quad , \quad
\FC_i^{\prime \prime} := \FC_{k(i)} . \]
The filtrations define spectral sequences
\[
E_*(\Morse )^\prime (LM_{h\TT}) \Rightarrow
H^*(LM_{h\TT}) \quad , \quad 
E_*(\Morse )^{\prime \prime}(LM_{h\TT}) \Rightarrow 
H^*(LM_{h{\TT}}). 
\]

The $p$ fold iteration map $\PS :=\PS_p$ induces a map of filtrations
$\PS :\FC_i \to \FC_{j(i)}$, so we have a ladder
\[ \xymatrix@C=1cm{
(\FC_0 )^{(p)}_{h\TT} \ar[r] \ar[d]^\PS 
& (\FC_1)^{(p)}_{h\TT} \ar[r] \ar[d]^\PS
& \dots \ar[r] & LM^{(p)}_{h\TT} \ar[d]^\PS \\
(\FC_0^\prime )_{h\TT} \ar[r] 
& (\FC_1^\prime )_{h\TT} \ar[r]
&\dots \ar[r] & LM_{h\TT}. \\
} \]
According to theorem \ref{th:IterationMap} we get an isomorphism
of spectral sequences 
\[ \xymatrix@C=1cm{
\PS^* :E_*(\Morse )^\prime (LM_{h\TT}) \left[ \frac 1 u \right]
\ar[r]^-{\cong} 
& E_*(\Morse )(LM^{(p)}_{h\TT} )\left[ \frac 1 u \right] .
} \]
Note that the spectral sequence on the right hand side can be
rewritten by the isomorphism in theorem \ref{th:ThreeSS}.

The filtration $\FC_i^{\prime \prime}$ is just a trivial reindexing
of the filtration $\FC_i^\prime$. 
We have that
\[
\FC_0^\prime =\FC_0^{\prime \prime}=\FC_1^{\prime \prime} =\dots =
\FC_{j(1)-1}^{\prime \prime} \subseteq 
\FC_1^\prime = \FC_{j(1)}^{\prime \prime} = 
\FC_{j(1)+1}^{\prime \prime} = \dots
\] 
So the two spectral sequences 
\[
E_*(\Morse )^{\prime \prime}(LM_{h\TT}) \left[ \frac 1 u \right]
\quad , \quad
E_*(\Morse)^\prime (LM_{h\TT}) \left[ \frac 1 u \right]
\]
are the same up to a re-indexing of the columns. 

Finally, since 
$\FC_i^{\prime \prime} = \FC_{k(i)} \subseteq \FC_i$ 
we have an inclusion of filtrations, which induces a map of spectral
sequences 
\[ \xymatrix@C=1cm{
E_*(\Morse )^{\prime \prime} (LM_{h\TT}) \left[ \frac 1 u \right]
\ar[r] 
& E_*(\Morse )(LM_{h\TT}) \left[ \frac 1 u \right] .
} \]
According to lemma \ref{le:ItLocalization} this map is an
isomorphism of $E_1$ pages, and thus of spectral sequences.
This proves the theorem.
\end{proof}

\section{The Morse spectral sequences for $\CP^r$}
\label{sec:MorseCPr}

We now turn to the special case $M=\CP^r$ with the usual symmetric space 
(Fubini--Study) metric. We know that this space satisfies the conditions 
of section \ref{sec:MSS}, in particular theorem \ref{th:relabel} is 
valid. We are going to study the Morse spectral sequence for 
$H^*(L\CP^r_ {h{\TT}})$. Using theorem \ref{th:relabel} 
we can pin down the possible structure of the differentials, but we 
cannot completely determine them. In particular, using only differential 
geometry, we can not prove that the differentials are non-trivial.
But we will make an effort to obtain as much information as possible 
about the Morse spectral sequence with the information we already have 
available at this point.

It will turn out in section \ref{se:MorseSerreSS} that there are 
non-trivial differentials in the Morse spectral sequence. In that 
section we are going to compare the Morse spectral sequence to a 
purely homotopy theoretical spectral sequence. This will not quite 
solve the differentials, but it will at least determine the dimension of
$H^*(L\CP^r_{h{\TT}})$ as a vector space over $\FF_p$.
 
\begin{theorem}
\label{th:MorsePoincare}
The Morse spectral sequence 
$E_*^{*,*}(\Morse )(L\CP^r_{h\TT})$ 
is a spectral sequence of 
$H(B{\TT})=\FF_p [u]$ modules converging towards
$H^*(L\CP^r_{h\TT})$. 
In case $p\mid (r+1)$ the $E_1$ page is given by
\[
\begin{aligned}
& E_1^{0,*} &
&=& &\FF_p [u,x]/(x^{r+1}),& 
& & \\
& E_1^{pm+k,*}& 
&=& &\alpha_{pm+k} \FF_p [x_1,x_2]/(Q_r,Q_{r+1}),\quad 
0\leq m, \quad 1\leq k\leq p-1, & \\
& E_1^{pm,*} &
&=& &\alpha_{pm} \FF_p [x,u] /(x^{r+1})
\oplus \zeta_{pm} \FF_p [x,u]/(x^{r+1}),\quad 1\leq m.& 
\end{aligned}
\]
In case $p\nmid (r+1)$ the $E_1$ page is given by
\[
\begin{aligned}
& E_1^{0,*} &
&=& & \FF_p [u,x]/(x^{r+1}), & \\
& E_1^{pm+k,*} &
&=& & \alpha_{pm+k} \FF_p [x_1,x_2]/(Q_r,Q_{r+1}), 
\quad 0\leq m,\quad 1\leq k \leq p-1, & \\
& E_1^{pm,*} &
&=& &\alpha_{pm} \FF_p [x,u]/(x^r)\oplus 
\bar \zeta_{pm} \FF_p [x,u]/(x^r),\quad 1\leq m.& 
\end{aligned}
\]
The total degree of the element $\alpha_{pm+k}x_1^ix_2^j$ is 
$2r(pm+k-1)+2i+2j+1$ and the filtration degree is $pm+k$, where
$1\leq k\leq p-1$.
The generators which are free $\FF_p[u]$-module generators have total 
degree and filtration degree as follows:

\begin{tabular}{l !{\quad} l !{\quad} l !{\quad} l}
\toprule
class & case & total degree & filtration\\
\midrule
$\alpha_{pm} x^i$ & all $r$, $0\leq i\leq r-1$ & $2r(pm-1)+2i+1$ & $pm$\\
$\alpha_{pm} x^r$ & $p\mid (r+1)$ & $2r(pm-1)+2r+1$ & $pm$\\
$\zeta_{pm} x^i$ & $p\mid (r+1)$ and $0\leq i\leq r$ & $2rpm+2i$ & $pm$\\
$\bar \zeta_{pm} x^i$ & $p\nmid (r+1)$ and $0\leq i\leq r-1$ & $2rpm+2+2i$ 
& $pm$ \\
\bottomrule
\end{tabular}
\end{theorem}

\begin{proof}
The sequence of critical values for the energy function is
given by $\lambda_n=n^2$. So in order to determine the $E_1$ page
we have to compute $\widetilde H^*(\Th (\mu^- (n^2))_{h\TT})$ 
for each $n$.

The case $n=0$ is special. The critical manifold $N(0)$
consists of the constant curves, so it is diffeomorphic
to $\CP^r$ itself, with trivial action of $\TT$. Since the
constant curves are the absolute minima of the
energy function, the negative bundle is the trivial bundle.
We find
\[
E_1^{0,*}(\Morse )(L\CP^r_{h\TT} )=
H^*(B\TT \times \CP^r)=
H^*(B\TT)\otimes H^*(\CP^r).
\]
This is what the theorem states for $E_1^{0,*}$.

Now consider $n>0$. The negative bundle $\mu^-(n^2)$ is a
$\TT$-vector bundle over $N(n^2)\cong S(\tau (\CP^r ))^{(n)}$. 
We know from theorem \ref{th:ThreeSS} that 
\[ E_1^{n,*}(\Morse )(L\CP^r_{h\TT} )\cong \widetilde 
H^*(\Th (\mu^-(n^2))_{h\TT} ).\] 
Recall from \cite{SySp} that if we forget the $\TT$-action,  
then $\mu^-(n^2)$ is the vector bundle 
$(n-1)\eta^* (\tau (\CP^r ))\oplus \epsilon$ over $S(\tau(\CP^r))$.
In particular, this means that $\mu^-(n^2)$ is orientable.
So by lemma \ref{lemma:Thom}, the vector bundle 
$\mu^-(n^2)_{h\TT}$ is orientable and 
$\Th (\mu^-(n^2))_{h\TT} \cong \Th ( \mu^-(n^2)_{h\TT}).$
The vector bundle $\mu^- (n^2)_{h\TT}$ is of the same dimension as
$\mu^- (n^2)$, that is of dimension $2r(n-1)+1$. So by the Thom
isomorphism we find that
\[
E_1^{n,*} (L\CP^r_{h\TT} )\cong \widetilde 
H^*(\Th (\mu^-(n^2)_{h\TT} )) \cong 
H^{*-(2r(n-1)+1)}(N(n^2)_{h\TT} ).
\]

The theorem now follows from the computation of the Borel cohomology 
of the space of geodesics in theorem \ref{th:CohomologyOrbits}.
\end{proof}

\begin{remark}
The symbol $\alpha_n x^i$ refers to the cup product in a critical manifold. 
The precise meaning is that it is the Thom isomorphism of the bundle
$\mu^{-}_n$ applied to the class $x^i$. In particular, the class 
$\alpha_n$ is the Thom class itself. The product is not defined 
on the cohomology of the Thom space of the negative bundle. 
It is also not defined in the spectral sequence. 
So strictly speaking, the product notation is improper. 
Products with $u$ on the other hand are genuine products,
which are defined in the spectral sequence.
\end{remark}

We also consider the non-equivariant case.
\begin{lemma}
\label{le:nonequivariant}
Consider the Morse spectral sequence 
\[
\{ E_r^{*,*}(\Morse )(L\CP^r ) \} \Rightarrow H^* (L\CP^r ).
\]  
We have the following formulas for it's $E_1$ page:
\[
E_1^{n,*}(\Morse )(L\CP^r )=
\begin{cases}
\FF_p [x]/(x^{r+1}) 
& , n=0, \\
\alpha_n \FF_p [x,\sigma ]/(x^{r+1},\sigma^2) 
& , n\geq 1, \quad p\mid (r+1), \\
\alpha_n \FF_p [x,\bar \sigma ]/(x^r,\bar \sigma^2 ) 
& , n\geq 1, \quad p\nmid (r+1). \\
\end{cases}
\]
The total degree of $\alpha_n$ is $2r(n-1)+1$, the total degree of
$\sigma$ is $2r-1$ and the total degree of $\bar\sigma$ is $2r+1$.
If $i\geq rn+1$ and $p\mid (r+1)$ or if $i\geq rn$ and $p\nmid (r+1)$
then $E_1^{n,2i+1-n}(\Morse )(L\CP^r )=0$ for all $n$. 
Finally, the canonical map 
\[ E_1^{n,2j+1-n}(\Morse )(L\CP^r_{h\TT} )\to
E_1^{n,2j+1-n}(\Morse )(L\CP^r )\] 
is surjective for all $n$ and all $j$.
\end{lemma}

\begin{proof}
Most of this follows directly from remark \ref{SCohomologyModP}. 
The class of highest even degree in  
$\FF_p [x,\sigma ]/(x^{r+1},\sigma^2 )$ is $x^r$ of degree $2r$
and in $\FF_p [x,\bar \sigma ]/(x^r,\bar \sigma^2 )$
it is $x^{r-1}$ of degree $2r-2$. The stated vanishing result follows
from this observation. The final surjectivity statement  
follows from the surjectivity result in corollary \ref{cor:SerreCollapse}. 
\end{proof}

A basic fact, which we will use again and again, is the 
following:
\begin{theorem}
\label{th:CPcollapse}
The Morse spectral sequence for non-equivariant cohomology collapses.
\end{theorem}

This theorem follows for instance from the stable splitting
of $L\CP^r$ that we constructed in \cite{SySp}. But the
theorem definitely does not belong to us. To prove 
the basic fact we only need a splitting on homology level,
and such a splitting was known to Ziller \cite{Ziller}.

\begin{corollary}
\label{cor:oddcounting}
Let $\FC_0 \subseteq \FC_1 \subseteq \dots \subseteq L\CP^r$ be the 
energy filtration and let $m\geq 0$. If $p\mid (r+1)$, the odd degree 
cohomology $H^{\text{odd} }(\FC_m)$ is a vector space of dimension 
$(r+1)m$. If $p\nmid (r+1)$, its dimension is $rm$.
\end{corollary}

\begin{proof}
According to theorem \ref{th:CPcollapse} the spectral sequence
$\{ E_r^{*,*} \}= \{ E_r^{*,*} (\Morse )(L\CP^r )\}$ collapses, so that
\[
H^{\text{odd}}(\FC_{m} ) \cong 
\bigoplus_{i,n} E^{n,2i+1-n}_\infty
\cong \bigoplus_{i,n} E^{n,2i+1-n}_1,
\]
where the direct sums are taken over pairs $n,i$ such that 
$0\leq n\leq m$ and $n\leq 2i+1$. By lemma \ref{le:nonequivariant} 
we have that $E_1^{0,2i+1}=0$, and for a fixed $n$ with $n\geq 1$ we 
have that 
\[
\sum_{n \leq 2i+1} \dim (E_1^{n,2i+1-n})=
\begin{cases}
r+1 & \text{if } p\mid (r+1), \\
r &\text{if } p\nmid (r+1), \\
\end{cases}
\] 
independently of $n$. It follows that
\[
\dim H^{\text{odd}} (\FC_{m}) =
\sum_{1\leq n\leq m} \medspace \sum_{n\leq 2i+1} \dim(E_1^{n,2i+1-n})=
\begin{cases}
(r+1)m & \text{if } p\mid (r+1),\\
rm & \text{if } p\nmid (r+1).\\
\end{cases}
\]
\end{proof}

Later it will be convenient to be able to estimate 
the sheer size of the Morse spectral sequence using a coarse 
dimension counting. So let's get that over with.

\begin{lemma}
\label{le:MorsePoincare}
In case $p\mid (r+1)$, the Poincar\'e series of 
$E_1^{*,*}(\Morse )(L\CP^r_{h\TT} )$ is
\[
\frac {1-t^{2r+2}} {(1-t)(1-t^2)(1-t^{2pr})}.
\]
In case $p\nmid (r+1)$, the Poincar\'e series of 
$E_1^{*,*}(\Morse )(L\CP^r_{h\TT} )$ is
\[
\frac {1-t^{2r+2}-t^{2pr}+t^{2pr+2}} {(1-t)(1-t^2)(1-t^{2pr})}.
\]
\end{lemma}

\begin{proof}
In theorem \ref{th:MorsePoincare} the $E_1$ term is written as a direct 
sum of three terms. We compute the Poincar\'e series of each term, and add.
First note that the Poincar\'e series of the algebra
$\FF_p [x_1,x_2]/(Q_r,Q_{r+1})$ is
\[ \frac {1-t^{2r}} {1-t^2} \medspace \frac {1-t^{2r+2}} {1-t^2}\]
as one sees by writing down a basis of monomials.

Here is the case $p\mid (r+1)$:
\begin{align*}
& \frac{1} {1-t^2} \medspace \frac {1-t^{2r+2}} {1-t^2} 
+\frac {t} {1-t^{2pr}} \medspace 
\frac {1-t^{2r(p-1)}} {1-t^{2r}} \medspace 
\frac {1-t^{2r}} {1-t^2} \medspace 
\frac {1-t^{2r+2}} {1-t^2} \medspace + \\
& \frac {1} {1-t^2} \medspace \frac {1-t^{2r+2}} {1-t^{2}} \medspace 
\frac {t^{2r(p-1)+1}+t^{2rp}} {1-t^{2pr}} ,
\end{align*}
which equals
\begin{align*}
& \frac{1-t^{2r+2}} {(1-t^2)^2 (1-t^{2pr})}
[(1-t^{2pr})+(t-t^{2r(p-1)+1})+(t^{2pr}+t^{2r(p-1)+1})] = \\
& \frac {(1-t^{2r+2})(1+t)} {(1-t^2)^2 (1-t^{2pr})} 
=\frac {1-t^{2r+2}} {(1-t)(1-t^2)(1-t^{2pr})}.
\end{align*}
In case $p\nmid (r+1)$, the last term in the direct sum is changed
and thus the Poincar\'e series equals the sum
\begin{align*}
& \frac{1} {1-t^2} \medspace \frac {1-t^{2r+2}} {1-t^2} 
+\frac {t} {1-t^{2pr}} \medspace 
\frac {1-t^{2r(p-1)}} {1-t^{2r}} \medspace 
\frac {1-t^{2r}} {1-t^2} \medspace 
\frac {1-t^{2r+2}} {1-t^2} \medspace + \\
& \frac {1} {1-t^2} \medspace \frac {1-t^{2r}} {1-t^2} \medspace
\frac {t^{2r(p-1)+1}+t^{2pr+2}} {1-t^{2pr}} ,
\end{align*}
which equals 
\begin{align*}
& \frac {1} {(1-t^2)^2 (1-t^{2pr})}
[(1-t^{2r+2}-t^{2pr}+t^{2pr+2r+2}) \\
& +(t-t^{2(p-1)r+1}-t^{2r+3}+t^{2pr+3}) 
+(t^{2(p-1)r+1}+t^{2pr+2}-t^{2pr+1}-t^{2pr+2r+2})] = \\
& \frac {1-t^{2r+2}-t^{2pr}+t^{2pr+2}} {(1-t)(1-t^2)(1-t^{2pr})} .
\end{align*}
\end{proof}

Now that we have completed the census, we will  investigate what we can 
say about the algebraic properties of the Morse spectral sequence 
for $\CP^r$.

\begin{lemma}
\label{le:naturality}
The classes $\zeta_{pm} x^i$, $\bar \zeta_{pm} x^i$, $\alpha_{pm} x^i$
and $\alpha_nx_1^i$
are not in the image of any differential. 
\end{lemma}

\begin{proof}
The inclusion map $i: L\CP^r \to L\CP^r_{h{\TT}}$
induces a map of filtrations, and a map of cohomology 
Morse spectral sequences. We know from theorem \ref{th:CPcollapse} 
that the target spectral sequence collapses.

But on the $E^1$ level, the classes $\zeta_{pm} x^i$ and 
$\bar \zeta_{pm}x^i$ survive to the classes represented by the Thom 
isomorphism applied to $\sigma x^i$ respectively $\bar \sigma x^i$, 
according to corollary \ref{cor:SerreCollapse}. These classes are not
in the image of a differential (since all differentials vanish). 
So by naturality, the classes $\zeta_{pm} x^i$ and $\bar \zeta_{pm}x^i$ 
cannot be in the image of a differential either.

Similarly, the classes $\alpha_{pm} x^i$ and $\alpha_nx_1^i$ 
are mapped non-trivially 
according to corollary \ref{cor:SerreCollapse} and the result follows.
\end{proof}

\begin{lemma} 
\label{le:Preliminary}
In the Morse spectral sequence $E_*^{*,*}(\Morse )(L\CP^r_{h\TT} )$
every non-trivial differential starts in an even total degree.
\end{lemma}
 
\begin{proof}
We have to do some counting, but the strategy of the proof is simple. 
First we note that it is enough to consider the generators of $E_1$, 
secondly we see that except for a very few special cases these 
generators cannot map non-trivially for dimensional reasons.
Thirdly we use lemma \ref{le:naturality} to dispose of the 
remaining cases.

The elements of odd degree in the spectral sequence are of the form
$\alpha_{pm} x^iu^j$ or $\alpha_nx_1^ix_2^j$. Because the spectral
sequence is a spectral sequence of $\FF_p[u]$ modules, to prove the
lemma it suffices to show the special case that all differentials 
vanish on the generators $\alpha_{pm} x^i$ and $\alpha_nx_1^i$.

We consider a class $d_s(\alpha_{pm} x^i)$. We wish to prove that it 
is trivial. This class has filtration $pm+s$ and total degree 
$2r(pm-1)+2i+2$. We ask, when does there exist a non trivial class 
of this filtration and total degree? To figure this out, we inspect 
the table of theorem \ref{th:MorsePoincare}, and realize the 
following facts.

\begin{itemize}
\item There exist non-trivial classes in 
$E_1 (\Morse )(L\CP^r_{h\TT} ;\FF_p )$
of total even degree, and of filtration $n$ if and only
if $n$ is divisible by $p$.
\item If $n$ is divisible by $p$ there is a unique such class 
of filtration $n$ with lowest possible even dimension. In case 
$p$ divides $r+1$ it is $\zeta_n$ of dimension $2rn$, 
if $p$ does not divide $r+1$ it is $\bar \zeta_n$ of dimension $2rn+2$. 
\end{itemize}

The class $d_s(\alpha_{pm} x^i)$  has total even degree, so if it 
were non-trivial it would to have degree at least
equal to the lowest possible degree of such a class. That is
\[
2r(pm-1)+2i+2 \geq 
\begin{cases}
2r(pm+s), &  p\mid (r+1), \\
2r(pm+s)+2, &  p\nmid (r+1).
\end{cases}
\]
If $p$ does not divide $r+1$, we have the two inequalities
$i\geq r(s+1)$ and $r-1 \geq i$. This system has no solution
with $s \geq 1$.
If $p$ divides $r+1$, we have the two inequalities
$i+1\geq r(s+1)$ and $r\geq i$. Now there is the unique solution
$s=r=i=1$. That is, we would have that
$d_1(\alpha_{pm-1}x)=\lambda \zeta_{pm}$ for a scalar $\lambda\not=0$. 
But this contradicts lemma \ref{le:naturality}.

A similar argument shows that $d_s(\alpha_nx_1^i)=0$.
\end{proof}

At this point, to a certain extent we can make a shortcut in the theory.
Using the elementary method as in the proof of 
lemma \ref{le:Preliminary} we can prove the following statement 
(theorem \ref{th:TransferClass} ). The term ``elementary'' should be 
read as ``without  using theorem  \ref{th:relabel}''. 
Together with the purely homotopy theoretical analysis of the
Serre spectral sequence (which we discuss in detail in section 
\ref{sec:SerreSS}), theorem \ref{th:TransferClass} is
sufficient to compute the cohomology $H^*(L\CP^r_{h\TT} )$ 
in the special case $p\mid (r+1)$. It does not seem possible to
cover the case $p\nmid (r+1)$ using a similar method.
It is precisely because of this that we need theorem \ref{th:relabel}.

\begin{theorem}
\label{th:TransferClass}
If $p$ divides $r+1$ then there is a class
$\hat \zeta_{pm} \in H^{2pmr} (L\CP^r_{h\TT} )$
such that the restriction 
$i^* (\hat \zeta_{pm}) \in H^{2pmr}(L\CP^r )$ 
is non-trivial.
\end{theorem}

\begin{proof}
The argument is very similar to the proof of lemma \ref{le:Preliminary}.
We show that the class 
$\zeta_{pm} \in E_1^{pm,(2r-1)pm}(\Morse)(L\CP^r_{h\TT} )$
of theorem \ref{th:MorsePoincare} survives to $E_\infty$. 

Suppose that $d_s(\zeta_{pm})=\alpha_{pm+s}y$
for some nontrivial scalar $y$. 
Equality of total degrees gives the relation
\[
2rpm+1=\deg (\zeta_{pm} )+1 =
\deg (\alpha_{pm+s} y)=2r(pm+s-1)+\deg(y)+1\geq 2rpm+1.   
\]
So we conclude that $s=1$, $\deg(y)=0$, and
$d_1(\zeta_{pm})=\lambda \alpha_{pm+1}$ for a non-trivial scalar $\lambda$.
But this contradicts lemma \ref{le:naturality}. That is, 
$\zeta_{pm}$ is a permanent cycle. Furthermore, $\zeta_{pm}$ is not
in the image of any differential by lemma \ref{le:naturality} ,
so it survives to $E_\infty$.

We claim that any class 
$\hat \zeta_{pn} \in H^{2pnr}(L\CP^r_{h\TT} ;\FF_p )$
representing $\zeta_{pm}$ will satisfy the property of the theorem.
By naturality, $\zeta_{pm}$ maps nontrivially by $i^*$ to the 
$E_1$ page of the Morse spectral sequence converging to
$H^*(L\CP^r ;\FF_p )$. The theorem follows, again since this spectral 
sequence collapses.
\end{proof}

In order to prove this sections final result we need a lemma
on non-negatively graded modules over the graded ring $\FF_p [u]$ 
(graded by $\deg(u)=2$). We say that the graded module $M$ is trivial 
in degrees strictly greater than $n$ if $M^i=0$ for all $i>n$.

Let us also say that a graded $\FF_p [u]$-module is generated in 
degrees less or equal to $m$ if there is a set of generators which 
have degrees less or equal to $m$. Because $M$ is bounded from below, 
a set of elements $\{ x_\alpha \}$ generate $M$ exactly if the 
reductions $[x_\alpha ]$ span the $\FF_p$-vector space $M/uM$. 

\begin{lemma}
\label{le:filterdegree}
Let $f:M \to N$ be a degree preserving map of non-negatively  
graded $\FF_p [u]$ modules. 
Assume that $M$ is generated in degrees less or equal to 
$m$, and that $N$ is trivial in degrees strictly greater than $n$. 
Then the kernel of $f$ is generated in degrees 
less or equal to $\max (m,n)$.
\end{lemma}

\begin{proof}
Assume without loss of generality that $f$ is surjective.
Multiplication by $u$ defines a map of the short
exact sequence $0\to \ker (f) \to M \to N \to 0$ to itself.
The snake lemma produces an exact sequence
\[
\xymatrix@C=1cm{
\ker (u:N \to N) \ar[r] & \ker(f) /u\ker(f) \ar[r] 
& M/uM \ar[r] & N/uN \ar[r] & 0. 
}
\]
That $M$ is generated in degrees less or equal to $m$
implies that the graded vector space $M/uM$ is zero in degrees
strictly greater than $m$. The exact sequence proves
that  $\ker(f) /u\ker(f)$ is trivial in degrees greater than 
$\max(m,n)$. Pick a set of spanning classes $\{ [x_\alpha]\}$ in 
$\ker(f) /u\ker(f)$, and lift them arbitrarily to classes 
$\{ x_\alpha\}$ in $\ker(f)$ of degree less or equal to $\max (m,n)$. These 
classes generate $\ker (f)$.
\end{proof}

We are now ready to apply our general theorem on the 
localized Morse spectral sequence to our particular specimen.

\begin{theorem}
\label{Epcollapse}
The spectral sequence $E_*(\Morse )(L\CP^r_{h\TT} )$ collapses
from the $E_p$ page.
\end{theorem}

\begin{proof}
Recall from lemma \ref{le:Preliminary} that all non-trivial 
differentials are defined on classes in even total degree. 
It is easy to see from theorem \ref{th:MorsePoincare} that the 
only classes of even total degree sit in filtrations divisible by $p$.
Let $X^{pm,*}_s \subseteq E_s^{pm,*}(\Morse )(L\CP^r_{h\TT} )$
be the subgroup of elements of even total degree. This is a module 
over the ring $\FF_p[u]$, and for all $s$ one has 
\[
X_{s+1}^{pm,*}=\ker \left( d_s :X_s^{pm,*}\to E_s^{pm+s,*} \right) .
\]

{\emph{Claim.}} If $1\leq s\leq p-1$, then the $\FF_p [u]$-module
$X^{pm,*}_{s+1}$ is generated by elements in degrees 
less than $2r(pm+s+1)-2$.

The claim follows from lemma \ref{le:filterdegree} and induction over
$s$. By theorem \ref{th:MorsePoincare} we see that $X_1^{pm,*}$ 
in all cases is generated in degrees less or equal to $2pmr+2r$.
$E_s^ {pm+s,*}$ is a quotient of the $\FF_p[u]$-module $E_1^{pm+s,*}$. 
Using theorem \ref{th:MorsePoincare} we see that for $1 \leq s \leq p-1$
this last module is trivial in dimensions strictly greater than 
$2r(pm+s)+2r-1$. Now we use  lemma \ref{le:filterdegree} on the 
differential
\[
d_1: X_1^{pm,*} \to E_1^{pm+1,*}.
\]

The differential raises total  degree by 1. Taking this
into account, lemma \ref{le:filterdegree} shows 
that $X^{pm}_1$ is generated in degree less than
$\max(2pmr+2r,2r(pm+1)+2r-2)=2r(pm+2)-2$, which proves
the claim for $s=1$.
Using  lemma \ref{le:filterdegree} inductively on $d_s$
for $2\leq s\leq p-1$ proves the claim.

In particular, $X^{mp,*}_{p}$ is generated in degrees
less or equal to $2r(pm+p)-2$. The basic fact \ref{th:CPcollapse} 
says that the Morse spectral sequence $E_*(\Morse )(L\CP^r )$ collapses
from the $E_1$ page. Because of theorem \ref{th:relabel} this implies 
that the localized spectral sequence 
$E_* (\Morse )(L\CP^r_{h\TT} )[1/u]$ 
collapses from the $E_1$ page. 

From our computation in theorem \ref{th:MorsePoincare} we can 
read off that the localization map
\[ \xymatrix@C=1cm{
E_1^{pm,*}(\Morse )(L\CP^r_{h\TT} )
\ar[r] 
& E_1^{pm,*}(\Morse )(L\CP^r_{h\TT} ) \left[ \frac 1 u \right]
} \]
is injective. By naturality, this means that no non-trivial 
differentials are arriving at
$E_1^{pm,*}(\Morse )(L\CP^r_{h\TT} )$. So in particular, the
$d_p$ differential on $X^{mp,*}_{p}$ is trivial, and
$X^{mp,*}_{p+1}=X^{mp,*}_{p}$ is generated in degrees
less than $2r(pm+p)-2$. 

We claim that all higher differentials $d_s$
vanish on $X_{p+1}^{pm,*}$. Assume to the contrary that
$s$ is the smallest $s\geq p+1$ which does not vanish on 
$X_s^{pm,*}$. Since $X_s^{pm,*}=X_{p+1}^{pm,*}$, we know that 
$X_s^{pm,*}$ is generated in degrees less or equal to $2r(pm+p)-2$.

It is enough to prove that the differential $d_s$ vanishes on the 
generators. A generator of $X_{p+1}^{pm,*}$ will be mapped by $d_s$
to a class in some $E_{p+1}^{pm+s,*}$ of total degree less or equal 
to $2r(pm+p)-1$. We have to show that there is no such non-trivial class.
 
But the non-zero class of lowest degree in $E_{p+1}^{pm+s,*}$ is 
$\alpha_{pm+s}$ which has degree $2r(pm+s-1)+1$.
Since we are assuming that $s\geq p+1$, this degree is larger than or
equal to $2r(pm+p)+1$. This finishes the proof.  
\end{proof}

\section{Derived functors at odd primes}
\label{sec:Derived}

In this section and in the following two sections we are studying 
the cohomology of the free loop space by homotopy theoretical methods. 
The main result we are aiming for is a computation of the action of
the circle on the cohomology. It seems hard to obtain this using Morse 
theory methods, since the action will cause non-trivial interaction 
between the layers in the Morse filtration. In section 
\ref{se:MorseSerreSS} we will show that this computation has consequences 
for the Morse spectral sequence $E_*(\Morse )(L\CP^r_{h\TT} )$, and 
essentially solves its differentials.

In this section we make a preliminary algebraic computation,
which we will need for the homotopy theory calculation.

Let as before $p$ be a prime number. Write $\trunc r x$ for
the truncated polynomial algebra $\FF_p [x] /(x^{r+1})$,
where $r\geq 1$ and the degree of the generator $|x|$ is
a positive even number. 
In \cite{BO} we defined the derived functors 
$H_*(\trunc r x ; \drf )$, and in the case $p=2$ we computed them.
In this section we extend this calculation to compute these 
derived functors for odd $p$. We will use definitions and notation 
from \cite{BO}. 

Assume that $p$ is an odd prime. Recall the functor 
$\drf :\F \to \Alg$ from \cite{SpSe}. The cohomology with 
$\FF_p$ coefficients of a space defines an object in $\F$, and 
$\Alg$ is the category of non-negatively graded algebras $A$
over $\FF_p$ such that $a^p=a$ for $a\in A^0$.  
We view $\trunc n x$ as an object in $\F$ with $\lambda =0$ and
$\beta =0$. Note that by definition of $\F$, any polynomial algebra 
$\FF_p [ z_i|i\in I]$ on even dimensional generators is a free object
in $\F$. This special type of free object has trivial $\lambda$ and 
$\beta$. By a similar argument as in \cite{BO} Theorem 2.1 we find the 
following result:

\begin{theorem} 
\label{res}
For odd primes $p$, there is an almost free simplicial resolution 
$R\simp \in s\F$ of the object $\trunc r x \in \F$ as follows: 
$R_q = \FF_p [x,y_1,\dots , y_q]$ for $q\geq 0$ 
where $|y_i|=(r+1)|x|$. The face and degeneracy maps
are given by $s_i(x)=x$, $d_i(x)=x$, 
\begin{align*}
s_i(y_j) &= \begin{cases}
y_j &, i\geq j ,\\
y_{j+1} &, i<j ,
\end{cases}
\\
d_i(y_j) &= \begin{cases}
x^{r+1} &, i=0, j=1, \\
y_{j-1} &, i<j, j>1, \\
y_j &, i\geq j , j<q, \\
0  &, i=q, j=q.
\end{cases}
\end{align*}
\end{theorem}

Using this resolution, the derived functors can be computed as the 
homotopy groups
$H_i(\trunc r x ; \drf ) = \pi_i \drf (R\simp )$.
It is convenient to use the normalized chain complex
$N_* (\drf R\simp )$ with $N_i (\drf R\simp )= \cap_{j=1}^i \ker (d_j)$
and differential $d_0$ for this purpose. 

The de Rham differential on $\drf R\simp$ is not the only 
simplicial derivation. There is is another one which turns out 
to be useful.

\begin{lemma} \label{std}
There is a well defined derivation $\theta : \drf R_q \to \drf R_q$
of degree 1 for each $q\geq 0$ which satisfies
$\theta (ab)=\theta (a) b +(-1)^{|a|}a \theta (b)$ 
for all $a,b\in \drf R_q$ and is defined by the following for
$1\leq j \leq q$: 
\[ \theta (x)=0,\quad \theta (y_j)=0 ,\quad \theta (\ddr x)=x,\quad
\theta (\ddr y_j)= (r+1)y_j. \]
One has that $\theta \circ \theta =0$. Furthermore,
$\theta$ commutes with the simplicial face and
degeneracy maps and hence defines a simplicial derivation
$\theta : \drf R\simp \to \drf R\simp .$
\end{lemma}

\begin{proof}
We have that $\drf R_q = \FF_p [x,y_1,\dots ,y_q] \otimes
\Lambda (\ddr x, \ddr y_1 ,\dots ,\ddr y_q )$. 
By the derivation property we see that $\theta ((\ddr x)^2)=0$
and $\theta ((\ddr y_j)^2)=0$ so $\theta$ is well defined.
The derivation property for $\theta$ also implies that 
$\theta \circ \theta$ is a derivation of degree two.
So $\theta \circ \theta$ is zero since it maps 
the algebra generators to zero.

One checks that  
$\theta (s_i z) =s_i\theta (z)$ and $\theta (d_iz) = d_i\theta (z)$
for each algebra generator $z$ by direct computations. 
The most interesting case goes as follows:
\begin{align*}
\theta (d_0 (\ddr y_1)) &= \theta (\ddr (x^{r+1})) =
\theta ((r+1)x^r\ddr x)=(r+1)x^{r+1}=d_0((r+1)y_1) \\
&= d_0\theta (\ddr y_1). 
\end{align*}
\end{proof}

Here is a complete computation of the derived functors.
\begin{theorem} 
\label{der}
If $p$ is a prime such that $p \mid (r+1)$, 
Then there is an isomorphism of bigraded $\FF_p$-algebras 
\[ H_*(\trunc r x ; \drf )\cong \trunc r x \otimes \Lambda (\ddr x )
\otimes \Gamma [\omega ], \]
where $\Vert x \Vert = (0,|x|)$, $\Vert \ddr x \Vert =(0, |x|-1)$, 
$\Vert \gamma_i (\omega ) \Vert =(i, i((r+1)|x|-1))$.
The algebra generators are represented by cycles in the normalized
chain complex $N_* (\drf R\simp )$ as follows:
$x=[x]$, $\ddr x= [\ddr x]$, 
$\gamma_i (\omega ) = [\ddr y_1 \dots \ddr y_i ]$. 
The de Rham differential induces the map
\[ \ddr_*: H_*(\trunc r x ;\drf ) \to H_*(\trunc r x ;\drf ); \quad
x \mapsto \ddr x, \quad \ddr x \mapsto 0, \quad 
\gamma_i (\omega )\mapsto 0. \]

If $p$ is a prime such that $p \nmid (r+1)$, then there is an 
isomorphism of bigraded $\FF_p$-algebras
\[ H_*(\trunc r x ; \drf )\cong \FF_p [a_i,b_i | i\geq 0]/I_r, \]
where $I_r$ is the ideal generated by the following elements
for $i,j \geq 0$:
\[ a_ia_j, \quad  
b_ib_j-\binom {i+j} i b_0 b_{i+j}, \quad
a_ib_j-\binom {i+j} i b_0 a_{i+j}, \quad
b_0^rb_i, \quad b_0^ra_i. \]
Here
$\Vert a_i \Vert = (i, i((r+1)|x|-1)+|x|-1)$ and
$\Vert b_i \Vert = \Vert a_i \Vert + (0,1)$. 
The algebra generators are represented by cycles in the normalized
chain complex $N_* (\drf R\simp )$ as follows:
$a_i=[\ddr x \ddr y_1 \dots \ddr y_i]$, 
$b_i=[\theta (\ddr x \ddr y_1 \dots \ddr y_i)]$.
The de Rham differential induces the map
\[ \ddr_*: H_*(\trunc r x ;\drf ) \to H_*(\trunc r x ;\drf ); \quad
a_i \mapsto 0, \quad b_i \mapsto (1+(r+1)i)a_i. \]
\end{theorem}

\begin{remark} \label{expl}
Explicitly, the cycle that represents $b_i$ is 
\[ \theta (\ddr x \ddr y_1 \dots \ddr y_i)=
x\ddr y_1 \dots \ddr y_i +(r+1)\ddr x \sum_{k=1}^i (-1)^k 
y_k \ddr y_1 \dots \widehat{\ddr y_k} \dots \ddr y_i. \]
\end{remark}

\begin{proof}
For $p=2$ these results were proved in \cite{BO} so assume that
$p$ is an odd prime. We first compute the derived functors as 
$\FF_p$-vector spaces. 

We have that $\trunc r x$ is the pushout of the diagram
$\FF_p \leftarrow \FF_p [y] \rightarrow \FF_p [x]$ in $\F$ where
$y\mapsto x^{r+1}$. Since $\drf$ commutes with colimits (appendix of
\cite{O}) and $\FF_p [x]$ is a free module over $\FF_p [y]$,  
proposition 6.3 of \cite{SpSe} gives us a Quillen spectral sequence
\[ E_{i,j}^2 = \Tor_i^{H_*(\FF_p [y];\drf )} 
(\FF_p , H_*(\FF_p [x];\drf ))_j \Rightarrow H_{i+j}(\trunc r x; \drf ). \]
Since polynomial algebras on even dimensional generators
are free objects in $\F$, we see that
$E_{i,j}^2=0$ for $j>0$ and that
\[ H_i(\trunc r x; \drf ) \cong E_{i,0}^2 \cong 
\Tor_i^{\drf (\FF_p [y])} (\FF_p , \drf (\FF_p [x])). \]

Note that $\drf (\FF_p [y])= \FF_p [y] \otimes \Lambda (\ddr y)$ is the
free graded commutative algebra on $\{ y,\ddr y\}$ so we have the
Koszul resolution $(K_*, \partial )$ of $\FF_p$ by free 
$\drf (\FF_p [y])$-modules:
\[ K_*=\Lambda (v)\otimes \Gamma [w]\otimes \drf (\FF_p [y]); \quad
\partial v =y,\quad \partial \gamma_i (w) = \gamma_{i-1}(w) \ddr y ,\]
where $v\in K_1$ and $\gamma_i (w)\in K_i$. In order to compute
the group $\Tor_i^{\drf (\FF_p [y])} (\FF_p , \drf (\FF_p [x]))$
we tensor this Koszul resolution with
$\drf (\FF_p [x] )$ over $\drf (\FF_p [y])$ and get a chain complex 
$(C_*,\partial )$ with
\[ C_* = \Lambda (v)\otimes \Gamma [w]\otimes \drf (\FF_p [x]), \quad
\partial v = x^{r+1}, \quad 
\partial \gamma_i (w) = (r+1)\gamma_{i-1} (w)x^r\ddr x. \]
Computing the homology of this chain complex, we find that if
$p \mid (r+1)$ then
\[ H_*(\trunc r x;\drf ) \cong 
\trunc r x \otimes \Lambda (\ddr x) \otimes \Gamma [w]. \]
If $p \nmid (r+1)$ we compute that 
$H_0(\trunc r x ;\drf ) \cong \trunc r x \otimes \Lambda (\ddr x) /
(x^r\ddr x)$. For $i>0$, 
$H_i(\trunc r x;\drf )$ is a sum of two copies of $\trunc r x/(x^r)$,
the generators being $\ddr x \gamma_i (w)$ and
$x\gamma_i (w)-(r+1)\ddr x v \gamma_{i-1}(w)$.
This completes the computation of the derived functors
as $\FF_p$-vector spaces.

We now check that the listed representatives are indeed cycles
in the chain complex $N_*(\drf R\simp )$. Define 
elements in $\drf R_i$ using the derivation $\theta$ as follows:
\[
\omega_i = \ddr y_1 \dots \ddr y_i ,\quad   
\alpha_i = \ddr x \omega_i ,\quad  
\beta_i = \theta (\alpha_i ).
\]
Here by convention $\omega_0=1$ such that 
$\alpha_0=\ddr x$ and $\beta_0= x$. 
We have that $d_j \omega_i =0$ for $0<j\leq i$ since $(\ddr y_j)^2=0$
and $d_iy_i=0$. Furthermore, 
$d_0\omega_i = (r+1)x^r\ddr x\omega_{i-1}$.
Thus $\omega_i$ is a cycle when $p \mid (r+1)$ and 
$\alpha_i$, $\beta_i$ are cycles when $p \nmid (r+1)$ as
stated. 

By a similar argument as the one presented in $\cite{BO}$
theorem 2.5 one checks that
\[ [x^j (\ddr x)^\epsilon \omega_i ], \quad 0\leq j \leq n ,\quad
\epsilon \in \{ 0 , 1 \} \]
are linearly independent in $H_i(N_*\drf R\simp )$ when
$p \mid (r+1)$ and that
\[ [\beta_0^j \alpha_i ] , [\beta_0^j \beta_i ] , \quad 
0\leq j \leq r-1 \]
are linearly independent in $H_i(N_*\drf R\simp )$ when
$p \nmid (r+1)$. By a dimension count, these two sets are then 
vector spaces bases. 

We now prove that the algebra structure of $H_*(\trunc r x ; \drf )$
is as stated. The algebra structure comes from the chain map
\[ \rho : C_*(\drf R) \otimes C_*(\drf R) \xrightarrow{g}
C_*(\drf R\otimes \drf R) \xrightarrow{C_*(m)} C_*(\drf R)\]
where $g$ is the Eilenberg--MacLane shuffle map \cite{ML} and
$m$ the multiplication map on $\drf R\simp$. By $C_*(V)$ for a 
simplicial $\FF_p$-vector space $V$, we mean the chain complex with
$C_i (V)=V_i$ and differential $\sum_{j=0}^i (-1)^j d_j$.
So we have the formula
\[ \rho (v_i \otimes w_j) = \sum_{(\mu, \nu)} (-1)^{\epsilon (\mu )}
s_\nu (v_i) s_\mu (w_j) \]
where the sum is over all $(i,j)$ shuffles $(\mu , \nu )$ of the set
$\{ 0, 1, \dots ,i+j-1\}$ and 
$\epsilon (\mu )= \sum_{k=1}^i (\mu_k -(k-1))$.

Lemma 2.4 of \cite{BO} still holds for $\FF_p$-coefficients.
(There is a small misprint in the lemma: There should be a hat over 
$\nu_j$ in the index set of the last formula.) So we find that
\[ \rho (\omega_i \otimes \omega_j ) = \binom {i+j} i \omega_{i+j} .\]
By this formula and remark \ref{expl} it follows directly that
\[ \rho (\alpha_i \otimes \alpha_j )=0 \quad , \quad 
\rho (\alpha_i \otimes \beta_j ) = 
\binom {i+j} i \beta_0 \alpha_{i+j} .\]

Since the degeneracy maps commutes with $\theta$, there is
a commutative diagram as follows where $A\simp = \drf R\simp$ :
\[
\diagram
A_i\otimes A_j \rrto^-{g} 
\dto_{\theta \otimes 1 +(-1)^{i+j} 1 \otimes \theta}
& &A_{i+j} \otimes A_{i+j} \rto^-{m} 
\dto_{\theta \otimes 1 + (-1)^{i+j} 1 \otimes \theta}
& A_{i+j} \dto^\theta \\
A_i\otimes A_j \rrto^-{g} & & A_{i+j}\otimes A_{i+j} \rto^-{m}  & A_{i+j}
\enddiagram
\]
By mapping $\alpha_i \otimes \theta (\alpha _j)$ both ways around
one finds that 
\[ \rho (\beta_i \otimes \beta_j) = 
\binom {i+j} i \beta_0 \beta_{i+j}. \]
Thus the algebra structure is as stated.
It follows directly form the formulas for the representing cycles,
that the de Rham map is as stated.
\end{proof}

\section{Cohomology of the free loop space}
\label{sec:FreeCohomology}

\begin{theorem}
\label{th:action}
Let $p$ be a prime. Assume that $X$ is a 1-connected space with 
mod $p$ homology of finite type. Assume also that
\[ H^*(X)=\FF_p [x]/(x^{r+1}), \]
where the degree $|x|$ of $x$ is even and $r\geq 1$. 
Put $\alpha=|x|$ and $\rho=(r+1)\alpha-2$. \\
1) If $p\mid (r+1)$ then there is an algebra isomorphism
\[ H^* (LX) \cong \FF_p [x]/(x^{r+1}) \otimes
\Lambda (\ddr x) \otimes \Gamma [\omega ], \]
where $|x|=\alpha$, $|\ddr x|=\alpha-1$, $|\gamma_i (\omega )|=\rho i$. 
The action differential is given by
\[ d :H^*(LX) \to H^*(LX);\quad
d(x)=\ddr x, \quad d(\ddr x)=0, \quad d(\gamma_i (\omega ))=0. \]
2) If $p\nmid (r+1)$ then there is an algebra isomorphism
\[ H^*(LX)\cong \FF_p [a_i,b_i|i\geq 0]/I, \]
where $I$ is the ideal generated by the following elements
for $i,j \geq 0$:
\[a_ia_j, \quad
b_ib_j-\binom {i+j} i b_0b_{i+j}, \quad 
b_ia_j-\binom {i+j} i b_0a_{i+j}, \quad
b_0^rb_i, \quad 
b_0^ra_i. \]
The degrees of the generators are $|a_i|=\rho i+\alpha-1$ and
$|b_i|=\rho i+\alpha$. In particular, the dimensions 
of $H^{2k}(LX)$ and $H^{2k-1}(LX)$ are the same. 
The action differential is given by
\[ d:H^*(LX) \to H^*(LX);\quad
d(a_i)=\kappa_i b_0^{r-1}b_i, \quad d(b_i)=((r+1)i+1)a_i. \]
where $\kappa_i=0$ unless $\alpha = 2$, $p=2$, $r$ is even and 
$i$ is odd.
\end{theorem}

\begin{remark} \label{rem:action}
When $p\nmid (r+1)$ we have the following formula for
$0\leq i$, $1\leq j \leq r$: 
\[ d(b_0^{j-1}b_i)=((r+1)i+j)b_0^{j-1}a_i. \]
We will show later (in corollary \ref{cor:actionproof}) that
$\kappa_i=0$. So we have $d( b_0^{j-1}a_i)=0$.
\end{remark}

\begin{proof}
According to \cite{SpSe} there is a strongly convergent second 
quadrant spectral sequence of cohomology type
\[ E_2^{-i,j} = H_i(H^*(X);\drf )^j \Rightarrow 
H^*(LX ), \]
where the $E_2$ page is the derived functors which we computed
in theorem \ref{der}. The spectral sequence is a spectral sequence
of algebras, and the induced of the de Rham differential $\ddr_*$
on $E_2$ corresponds to the action differential $d$ on 
$H^*(LX)$.

There are two distinct cases: $p\mid (r+1)$ and
$p\nmid (r+1)$. In both cases, the $E^2$ page does not admit
any differentials because of its distribution of zeros,
so $E_\infty \cong E_2$ as a vector space over $\FF_p$. 
The theorem can be paraphrased as that the $E_2$ page
is isomorphic to $H^*(LX)$ as an algebra, and
that also the action differential agrees. 
We have to look close to exclude the possibility of multiplicative 
extensions as well as extension problems concerning the action 
differential. Write $|\cdot |$ for the total degree in the spectral 
sequence.

1) Assume that $p\mid (r+1)$. Since $x$ and $\ddr x$ lie
in $H_0(\trunc r x ;\drf )$ they have unique representatives
in $H^*(LX)$ which satisfy $x^{r+1}=0$, $(\ddr x)^2=0$ and
$d(x)=\ddr x$. So the possible extension questions are if 
the relations $\gamma_i^p=0$ and $d\gamma_i=0$ are valid.
Let us denote a class in $H^*(LX)$ representing 
$\gamma_i \in E_\infty$ by the symbol $\overline \gamma_i$. 

We look at $\overline \gamma_i^p$. We know that $\overline \gamma_i^p=0$
up to classes of the same total degree and strictly higher
filtration. The class $\gamma_{ip}$ has the same total degree, but
also the same filtration degree as $\gamma_i^p$. 
Using that $ |\gamma_{j+1} |=|x^r \gamma_j |+\alpha -2$, we see that if
$\alpha \geq 4$ this is the only class of the same total degree
as $\gamma_i ^p$, so that independently of the choice of 
$\overline \gamma_i$, we indeed have that $\overline \gamma_i^p=0$.

Similarly $d\gamma_i=0$ up to classes of the same total degree and
of strictly higher filtration. But if $\alpha \geq 4$, using that 
$|d \gamma_{j+1} |=|x^{r-1}\ddr x \gamma_j |+\alpha -2$
we see that there is no such class. So $d\overline \gamma_i=0$. 

Now we consider the slightly more complicated case $\alpha =2$. 
In this case the class $x^r\gamma_j$ has the same total degree as 
$\gamma_{j+1}$ and strictly higher filtration. It is the only
such class. Similarly, $x^{r-1}\ddr x \gamma_j$ is the 
only class of same total degree as $d\gamma_{j+1}$ and of strictly 
higher filtration. 

The filtration of $H^*=H^*(LX)$ has the form
\[ H^* \supseteq \dots \supseteq F^{-i}H^* \supseteq F^{-i+1}H^*
\supseteq \dots \supseteq F^0H^* \supseteq F^1H^*=0. \]

We now show that we can choose $\overline \gamma_i \in F^{-i}H^*$
such that $[\overline \gamma_i] = \gamma_i \in E_\infty^{-i,*}$ and
$d\overline \gamma_i =0$. For $i=0$ the unit $\overline \gamma_0 =1$
has these properties. Assume that $\overline \gamma_j$ has been chosen
with these properties for $1\leq j < i$. Choose 
$\gamma_i^\prime \in F^{-i}H^*$ such that 
$[\gamma_i^\prime ] =\gamma_i \in E_\infty^{-i,*}$.
We have $[d\gamma_i^\prime ]=0$ so 
$d\gamma_i^\prime =kx^{r-1}dx \overline \gamma_{i-1}$ for
some $k\in \FF_p$. Put $\overline \gamma_i = 
\gamma_i^\prime-\frac k r x^r\overline \gamma_{i-1}$.
Then, $[\overline \gamma_i] = \gamma_i \in E_\infty^{-i,*}$ and
$d\overline \gamma_i= 
d\gamma_i^\prime -k x^{r-1}dx\overline \gamma_{i-1} =0$.

We claim that these representatives satisfy 
$(\overline \gamma_{p^i} )^p=0$ for $i\geq 0$.
From the spectral sequence we see that
$(\overline \gamma_{p^i} )^p=c x^r\overline \gamma_{p^{i+1}-1}$
for some $c \in \FF_p$. We apply the action differential 
on this equation and find
\[ 0 = d((\overline \gamma_{p^i})^p) = 
crx^{r-1}dx \overline \gamma_{p^{i+1}-1}, \]
and since $r\equiv -1$ mod $p$ it follows that $c=0$.
 
So there is a well defined algebra map
\[ \phi :\FF_p [x]/(x^{r+1}) \otimes \Lambda (\ddr x) \otimes 
\Gamma [\omega ] \to H^*(LX), \]
which maps $\gamma_i (\omega )$ to $\overline \gamma_i$.
We can filter the domain such that $\phi$ becomes a map of filtered
algebras. Since $\phi$ is an isomorphism on associated graded objects, it is
itself an isomorphism.

2) Assume that $p\nmid (r+1)$.
The classes and multiplicative relations in $E_2$ are as follows:

\begin{tabular}{l !{\quad} l !{\quad} l !{\quad} l}
\toprule
class  or relation & type & total degree/parity & filtration \\
\midrule
$b_0^ia_j$,\quad \medspace $0\leq i\leq r-1$ & class 
& $\rho j+\alpha(i+1)-1$, odd & $-j$ \\
$b_0^ib_j$, \quad $0\leq i\leq r-1$ & class & $\rho j+\alpha(i+1)$, even 
& $-j$ \\
$b_ia_j-\binom {i+j} i b_0a_{i+j}$ & relation & $\rho(i+j)+2\alpha-1$, odd
& $-i-j$ \\
$b_0^ra_j$ & relation & $\rho j+\alpha (r+1)-1$, odd & $-j$ \\
$a_ia_j$ & relation & $\rho (i+j)+2\alpha-2$, even & $-i-j$ \\
$b_ib_j-\binom {i+j} i b_0b_{i+j}$ & relation & $\rho(i+j)+2\alpha$, even
& $-i-j$ \\
$b_0^rb_j$ & relation & $\rho j+\alpha (r+1)$, even &$-j$\\
\bottomrule
\end{tabular}

\noindent We note that if $0\leq i\leq r-1$, then 
\[
0 < \alpha \leq \alpha(i+1) \leq \alpha r 
\leq \rho.
\]
From this and the list of classes above it follows easily 
that there is at most one class in each total degree. 
A differential in the spectral sequence raises total
degree by one, and it strictly increases the
filtration.
The unique class of total degree one higher than $a_j$
is $b_j$, which has the same filtration as $a_j$, so that $a_j$ is a 
permanent cycle. If $\alpha\geq 4$, there is no non-trivial class
of degree one higher than $b_j$. 
If $\alpha=2$, the unique class of degree one higher than $b_j$ is 
$b_0a_j$ which again has the same filtration so $b_j$ is also a
permanent cycle.

Each of the relations given above is true for any lifting of the 
generators in $E_\infty$ to generators of $H^*(LX)$ up to 
classes of the same total degree and strictly higher filtration.
We claim that in each case, there are no non-trivial such classes.

If $\alpha \geq 4$, there are no nontrivial classes of the
same dimension as $b_0^ra_j$ or $b_0^rb_j^r$. In case
$\alpha=2$, there are the unique classes $a_{j+1}$
respectively $b_{j+1}$. But these have lower 
filtration than the relation, so they cannot contribute to 
extensions.

If $\alpha \geq 4$ there is no class which has the same
total degree as the relation $a_ia_j$. In case
$\alpha =2$, there is the unique possibility $b_{i+j}$,
which does not matter anyway since its filtration 
is too low.

Also, the relations 
$b_ia_j-\binom {i+j} i b_0a_{i+j}$ and
$b_ib_j-\binom {i+j} i b_0b_{i+j}$ have the same
total degree as $b_0a_{i+j}$
respectively $b_0b_{i+j}$. By filtration check, 
there can be no extension problems.

Finally we have to consider the action differential.
The filtration argument says that  
$db_i$ is as stated. We have to argue that
$da_i=0$. If $\alpha\geq 4$, there is no 
non-trivial possibility for $da_i$. If
$\alpha=2$, the class $da_i$ has the same
total degree as $b_{i-1}b_0^{r-1}$. 
Because $d$ is a differential, $d(a_i)=0$
unless 
\[
0=db_i=(1+(r+1)i)a_i,
\]
 and 
\[
0=d(b_{i-1}b_0^{r-1})=
(1+(r+1)(i-1))b_0^{r-1}a_{i-1}+
(r-1)b_{i-1}b_0^{r-2}a_0=
((r+1)i-1)b_0^{r-1}a_{i-1}.
\]
Solving this, we get that $2=0$ (that is $p=2$),
that $r$ is even and that $i$ is odd.
\end{proof}

We will later want to know how big the cokernel of the 
action differential is, so we do the counting now.
We write $\NN = \{ 1,2,3,\dots \}$ for the set of natural numbers.

\begin{definition}
\label{def:TF}
Let $p$ be a prime and let $r,\alpha \in \NN$ with $\alpha$ even.
Put $\rho=(r+1)\alpha-2$ and let   
\[
\chi_p (s) = \begin{cases}
0 & \text{if } p \mid s, \\
1 & \text{if } p \nmid s.
\end{cases} 
\]
We define two subsets of $\NN$ as follows:
\begin{align*}
& \IF (r,p,\alpha )=
\{ \rho i+\alpha j \mid \chi_p (r+1)\leq j\leq r, 0\leq i 
\text{ and } p\mid ((r+1)i+j) \} \setminus \{ 0 \} , \\
& \IT (r,p,\alpha )=
\{ \rho i+\alpha j \mid \chi_p (r+1)\leq j\leq r, 0\leq i 
\text{ and } p\nmid ((r+1)i+j) \} .
\end{align*}
\end{definition}

The notation $\IF$ refers to an {\em index} set for {\em free} 
generators and $\IT$ refers to an {\em index } set for
{\em transfer} generators. This choice of notation will make sense
later on.
 
\begin{lemma}
\label{le:twosets} 
Whether a natural number $k$ is contained in 
$\IF (r,p,\alpha )$ respectively in $\IT (r,p,\alpha )$ only 
depends on the congruence class of $k$ modulo $\rho p$.
Furthermore, 
\begin{align*}
& \IF (r,p,\alpha )\cap \IT (r,p,\alpha ) =
\begin{cases}
2r \NN & \text{if $p\mid (r+1)$ and $\alpha=2$,}\\
\emptyset& \text{otherwise,}  
\end{cases}\\    
& \IF(r,p,2)\cup \IT(r,p,2) = 2\NN.
\end{align*}
If a number is in $\IF (r,p,\alpha )$ or in $\IT (r,p,\alpha )$, 
there is a unique choice of numbers $i,j$ that displays it as such.
If $p\mid (r+1)$, the set
$\{ 0< 2k\leq \rho p m \mid 2 k\in \IF (p,r,\alpha ) \}$
has $m(r+1)$ elements. If $p\nmid (r+1)$, the set has $mr$ elements. 
\end{lemma}

\begin{proof}
Note that $0\not\in \NN$. We warm up by first ignoring the 
congruence conditions. Assume that $i,i^\prime\in \ZZ$ and that
$\chi_p (r+1) \leq j, j^\prime \leq r$.
We make two claims about this situation.
\begin{enumerate}
\item If $i^\prime >i$ and 
$\rho i+\alpha j=\rho i^\prime+\alpha j^\prime$, then
$\alpha =2$, $p\mid (r+1)$, $i^\prime = i+1$, $j^\prime=0$ and $j=r$.
\item If $\rho i+\alpha j > 0$, then $i\geq 0$.
\end{enumerate}
To prove the first claim, note that
$\alpha (j-j^\prime)=\rho (i^\prime-i)\geq \rho$. So
\[
j-j^\prime \geq \frac \rho \alpha = r+\frac {\alpha-2} {\alpha}.
\]
Since $0\leq j,j^\prime\leq r$ and $\alpha\geq 2$, 
this is only possible if
$j=r$, $j^\prime=0$, $\alpha=2$. But this implies that 
$i^\prime =i+1$ and by assumption, if $j^\prime =0$, 
then $p\mid (r+1)$. 

To prove the second claim, note that
\[
i> -\frac {\alpha j} {\rho} \geq 
-\frac {\alpha r} {\rho} =
-\frac {\alpha r}{\alpha (r+1)-2} \geq 
-\frac {\alpha r} {\alpha r}= -1.
\]

We now prove the lemma. First assume that 
$k$ and $k+m\rho p$ are natural numbers. Assume also that
$k\in \IF (r,p,\alpha )$. We can write $k=\rho i+\alpha j$, where 
$i,j$ satisfy the appropriate conditions.
Then $k+m\rho p=\rho(i+mp)+\alpha j$, and the pair $(i+mp,j)$
satisfies the same congruence conditions, and conditions on $j$. 
This proves that $k+m\rho p\in \IF (r,p,\alpha )$, if we know that 
$i+mp\geq 0$. But this follows from claim 2. together with our 
assumption that $k+m\rho p$ is a natural number.

The same argument shows that $\IT (r,p,\alpha )$ is also a union
of congruence classes of natural numbers.

If $x\in \IF (r,p,\alpha )\cap \IT(r,p,\alpha )$, it must be 
possible to write $x$ in two different ways in the form 
$x=\rho i+\alpha j$.

By the first claim, we get that the only possible way
this can happen is that $\alpha=2$, $p\mid (r+1)$ and
$x=2ri + 2r=2r(i+1)+0$. This proves that 
$\IF (r,p,\alpha )\cap \IT (r,p,\alpha )$ is empty unless
$p\mid (r+1)$ and $\alpha =2$. It also shows that 
$\IF (r,p,2)\cap \IT(r,p,2) \subseteq 2r\NN .$
We have to show that if $p\mid (r+1)$
then $2r\NN \subseteq \IF (r,p,2)\cap \IT (r,p,2)$. 
In this case we write
$2rm=\rho i+\alpha j=\rho i^\prime +\alpha j^\prime$
for $(i,j)=(m,0)$ and $(i^\prime ,j^\prime )=(m-1,r)$.
This proves the claim, since $p\mid j$ but $p\nmid j^\prime$.

We have $\IF (r,p,2)\cup \IT (r,p,2) =
\{ 2(ri+j) \mid \chi_p (r+1) \leq j \leq r, 0 \leq i \}$,
which equals $2\NN$ as stated. The uniqueness statement on
$i,j$ follows directly from claim 1.

To prove the final statement about the number of elements, 
it is enough to show that the number of congruence
classes in $\IF (r,p,\alpha )$ modulo $\rho p$
is $r+1$ respectively $r$. In case $p\mid (r+1)$
the congruence classes of $2k$ are the classes
of form $\rho i+\alpha pj^\prime$ for $0\leq i < p$
and $0\leq j^\prime <(r+1)/p$, and there are clearly
$p((r+1)/p)=r+1$ of those. In case $p\nmid (r+1)$,
each $j$ uniquely determines a congruence class 
$i(j)$ modulo $p$ such that 
$(r+1)i(j)+j\equiv 0\mod p$. That is, each $j$
with $1\leq j \leq r$ uniquely determines an $i$,
$0\leq i<p$ such that $(r+1)i+j\equiv 0\mod p$. 
So there are exactly $r$ pairs $(i,j)$ qualifying,
and there are $r$ congruence classes in $\IF (r,p,\alpha )$.
\end{proof}

\begin{example}
Let $\alpha=2$, $r=2$. Then $\rho = 4$ and 
\begin{align*}
& \IF (2,2,2) = \{ 4+8m, 6+8m | \medspace m\geq 0 \} , \\ 
& \IF (2,3,2) = 4\NN , \\ 
& \IF (2,5,2) = \{ 8+20m, 14+20m | \medspace m\geq 0 \} , \\
& \IF (2,7,2) = \{ 10+28m, 20+28m | \medspace m\geq 0 \} .
\end{align*}
\end{example}

\begin{lemma}
\label{le:countingtransfer}
Let $X$ be as in theorem \ref{th:action}, $k\in \NN$ and put $H^*=H^*(LX)$. 
\begin{enumerate}
\item \label{action:kernel}
The kernel of the action differential
$d: H^{2k} \to H^{2k-1}$ is either a trivial or a one dimensional 
vector space. It is non-trivial if and only if 
$2k\in \IF(r,p, \alpha )$.
\item \label{action:image}
The image of $d:H^{2k}\to H^{2k-1}$ is either a trivial or a one 
dimensional vector space. It is non-trivial if and only if 
$2k\in \IT(r,p, \alpha )$.
\item \label{action:cokernel}
The cokernel of $d:H^{2k}\to H^{2k-1}$ is either a trivial or a one 
dimensional vector space. In case  $p\nmid (r+1)$, it is non-trivial 
if and only if $2k\in \IF (r,p,\alpha )$. 
In case $p\mid (r+1)$ it is non-trivial 
if and only if either $2k\in \IF(r,p,\alpha )$
and $\rho \nmid 2k$ or if $k>1$ and $2k\equiv 2\mod \rho$.
\item \label{action:kernelsum}
The kernel of the map
$d:\oplus_{2\leq 2k\leq \rho pm} H^{2k}\to
\oplus_{1\leq 2k-1\leq \rho pm-1} H^{2k-1}$
is a vector space of dimension $mr$ if $p\nmid (r+1)$, and
of dimension $m(r+1)$ if $p\mid (r+1)$.
\item \label{action:cokernelsum1}
The cokernel of the map
$d:\oplus_{2\leq 2k\leq \rho pm} H^{2k}\to
\oplus_{1\leq 2k-1\leq \rho pm-1} H^{2k-1}$ is a vector space of 
dimension $rm$ when $p\nmid (r+1)$.
\item \label{action:cokernelsum2}
The cokernel of the map
$d:\oplus_{2\leq 2k\leq \rho pm+2} H^{2k}\to
\oplus_{1\leq 2k-1\leq \rho pm+1} H^{2k-1}$ 
is a vector space of dimension $(r+1)m$ when $p\mid (r+1)$.
\end{enumerate}
\end{lemma}

\begin{proof}
For the action differential $d:H^{2k} \to H^{2k-1}$ we have the equation
\begin{gather}
\begin{split}
\label{eq:cokerdimension}
\dim \coker (d) &=\dim H^{2k-1}-\dim \im (d)\\
&= \dim H^{2k-1} - (\dim H^{2k}-\dim \ker (d))\\
&= -(\dim H^{2k}-\dim H^{2k-1})+\dim \ker (d).
\end{split} 
\end{gather}

We first consider the case $p\nmid (r+1)$. 
The even part $H^{\text{even}}$ has basis $\{ b_0^{j-1}b_i\}$ and
the odd part $H^{\text{odd}}$ has basis 
$\{ b_0^{j-1}a_i\}$ where $0\leq i$, $1\leq j\leq r$. The basis
elements sit in degrees $\rho i +\alpha j$ and $\rho i + \alpha j -1$
respectively. The kernel of the action differential
$d:H^{\text{even}}\to H^{\text{odd}}$ is generated by 
those $b_0^{j-1}b_i$ for which $p\mid ((r+1)i+j)$ and
its image of those $b_0^{j-1}a_i$ for which $p\nmid ((r+1)i+j)$.
Combining this with the uniqueness statement for $i,j$ 
in lemma \ref{le:twosets} we see that \ref{action:kernel}. and
\ref{action:image}. are valid. Furthermore we see that 
$\dim H^{2k}=\dim H^{2k-1}$ which combined with 
(\ref{eq:cokerdimension}) gives us \ref{action:cokernel}.

Next we consider the case $p\mid (r+1)$. Here $H^{\text{even}}$ has
basis $\{ x^j\gamma_i (\omega ) \}$ and $H^{\text{odd}}$ has basis
$\{ x^jdx\gamma_i (\omega )\}$ where $0\leq i$, $0\leq j\leq r$.
The basis elements sit in degrees $\rho i+\alpha j$ and 
$\rho i+\alpha j+\alpha -1$ respectively. The kernel for the action
differential $d:H^{\text{even}}\to H^{\text{odd}}$ is generated by those
$x^j\gamma_i (\omega )$ for which $p\mid j$ and the image is
generated by those $x^{j-1}dx\gamma_i (\omega )$ for which
$p\nmid j$. As above lemma \ref{le:twosets} give us that
\ref{action:kernel}. and \ref{action:image}. are valid. 

We now prove \ref{action:cokernel}. One checks it directly 
for $(r,p,\alpha )=(1,2,2)$. Assume that $(r,p,\alpha )\neq (1,2,2)$ 
which implies $\rho >2$. By a counting argument we find the following:
\begin{equation}
\label{eq:homdimension}
\dim H^{2k}-\dim H^{2k-1}=
\begin{cases}
1 & \text{if $\rho \mid 2k$,}\\
-1 & \text{if $\rho \mid (2k-2)$ and $k>1$,}\\
0 & \text{otherwise.}
\end{cases}
\end{equation}
We combine (\ref{eq:homdimension}) and (\ref{eq:cokerdimension}) 
in order to prove statement \ref{action:cokernel}. 

If $\rho \mid 2k$ we have $2k\in \IF (r,p,\alpha )$ so the
dimension of the cokernel becomes $0$. Assume that
$\rho \mid (2k-2)$ and $k>1$. We claim that this implies that
$2k\notin \IF (r,p,\alpha )$. For if $2k\in \IF (r,p,\alpha )$
we have that $2k=\rho i+\alpha j$ where $0\leq i$, $0\leq j \leq r$
and $p\mid j$. It follows that
$2\equiv 2k \equiv \alpha j \mod \rho$. Since 
$\rho >2$, we cannot have $j=0$. We conclude that
$1 \leq p\leq j$. On the other hand,
$0\equiv 2k-2=\rho i+ \alpha j -2\equiv \alpha j-2\mod \rho$,
so $\rho \leq \alpha j-2$. Now we have our contradiction, 
finishing the proof of \ref{action:cokernel}.
since $\alpha j-2 \leq \alpha r-2 <\rho$.

We have reduced the last three statements of the lemma
to statements  about the sets $\IT (r,p,\alpha )$ and 
$\IF (r,p,\alpha )$. We see that \ref{action:kernelsum}.
is equivalent to the statement that the set
$\{ 0< 2k \leq \rho pm | 2k \in \IF (r,p,\alpha) \}$
has $rm$ elements if $p\nmid (r+1)$ and 
$(r+1)m$ elements if $p\mid (r+1)$. But this is exactly 
the content of the last statement of lemma \ref{le:twosets}.
Statement \ref{action:cokernelsum1}. follows from statement
\ref{action:kernelsum}. and (\ref{eq:cokerdimension}). 

Finally, we prove statement \ref{action:cokernelsum2}.
One verifies it directly for
$(r,p,\alpha )=(1,2,2)$. Assume that $(r,p,\alpha )\neq (1,2,2)$. 
Formula (\ref{eq:homdimension}) gives us that
\[
\sum_{2\leq 2k\leq \rho pm+2}{(\dim H^{2k}-\dim H^{2k-1})}=0.
\]
So by (\ref{eq:cokerdimension}), \ref{action:kernelsum}. and 
\ref{action:kernel}. it suffices to check that 
$\rho p m +2\notin \IF (r,p,\alpha )$. This follows from
lemma \ref{le:twosets} since $\rho pm+2\equiv 2$ mod $\rho p$.
\end{proof}

We define two formal power series in $\ZZ [[t]]$ by
\[ 
P_{\IT (r,p,\alpha )}(t)= \sum_{n \in \IT (r,p,\alpha )}t^n,
\quad
P_{\IF (r,p,\alpha )}(t)= \sum_{n\in \IF (r,p,\alpha )}t^n.
\]
Note that the constant terms are zero in both series. Furthermore,
if the numbers $\kappa_i$ in theorem \ref{th:action} are
zero for all $i$, then $P_{\IT (r,p,\alpha )}(t)=tP_{\im (d)}(t)$ where 
$P_{\im (d)}(t)$ is the Poincar\'e series for the image of the 
action differential. By lemma \ref{le:twosets} we find the following 
result:
\begin{lemma} \label{le:Poincaretwosets}
When $\alpha =2$ such that $\rho = 2r$ we have that
\[  
P_{\IT (r,p,2)}(t)+P_{\IF (r,p,2)}(t)=
\frac {t^2} {1-t^2}+(1-\chi_p (r+1)) \frac {t^{2r}} {1-t^{2r}} .
\]  
\end{lemma}

\section{The $E_3$-term of the Serre spectral sequence}
\label{sec:SerreSS}

Let $Y$ be a $\TT$-space with 
$H_*(Y)$ of finite type. Write $q:E\TT \times Y \to E\TT \times_\TT Y$
for the quotient map. As described in appendix
\ref{Appendix:s1transfer} there is a $\TT$-transfer map
$\tau$ such that the composite $q^*\circ \tau$ equals the action 
differential:
\[ d: H^{*+1}(Y) \xrightarrow{\tau } H^*(Y_{h\TT}) \xrightarrow{q^*} 
H^*(Y) .\] 
Since $B\TT$ is 1-connected and $H^*(B\TT) =\FF_p [u]$ where $|u|=2$,
the Serre spectral sequence for the fibration $Y\to Y_{h\TT} \to B\TT$
has the following form:
\[ E_2^{*,*} = \FF_p [u] \otimes H^*(Y) \Rightarrow H^*(Y_{h\TT}). \]
The $d_2$ differential is determined by the action differential since 
$d_2y=udy$ for all $y \in H^*(Y) $. Thus, the $E_3$ term has the
following form:
\[ E_3^{*,*} = \image (d) \oplus (\FF_p [u]\otimes H(d)), \]
where $\image (d)$ and $H(d)$ denotes the image and the homology of the 
action differential respectively.

\begin{proposition} \label{imd}
The subspace $\image (d) \subseteq E_3^{*,*}$ survives to 
$E_\infty^{*,*}$. For any $a\in H^*(Y)$ one has that 
$\tau (a) \in H^*(Y_{h\TT })$ represents 
$da \in E_{\infty }^{0,*}$ and that $pr_1^*(u)\tau (a)=0$ 
in $H^*(Y_{h\TT })$, where $pr_1 : Y_{h\TT} \to B\TT$ denotes the 
projection on the first factor.
\end{proposition}

\begin{proof}
There are two commutative diagrams 
\[ \xymatrix@C=1cm{ 
H^*(Y) \ar[ddr]_-{d} \ar[r]^-{\tau} 
& H^*(Y_{h\TT }) \ar[dd]_-{q^*} \ar@{>>}[r] 
& E_\infty^{0,*} \ar@{^{(}->}[d] 
& & H^*(Y_{h\TT }) \ar[r]^-{q^*}  
& H^*(E\TT \times Y) \\
& & E_3^{0,*} \ar@{^{(}->}[d] 
& & H^*(B\TT ) \ar[r] \ar[u]^-{pr_1^*} & H^*(E\TT ) \ar[u]_-{pr_1^*}\\
& H^*(Y) \ar@{=}[r] & E_2^{0,*}
} \]
Assume that $a\in H^*(Y)$ has $da\neq 0$. The diagram to the left 
shows that $\tau (a)\neq 0$ and that $da$ survives to $E_\infty$ and
is represented by $\tau (a) \in H^*(Y_{h\TT })$.
The diagram to the right shows that $q^*\circ pr_1^*(u)=0$ since 
$E\TT$ is contractible. By Frobenius reciprocity
$pr_1^*(u)\tau (b) = \tau (q^*(pr_1^* (u))b)=0$.   
\end{proof}

We now take $Y=LX$. By the result in the previous section we 
can compute the $E_3$ term of the Serre spectral sequence.

\begin{proposition} \label{th:E3Serre}
Let $p$ be a prime. Assume that $X$ is a 1-connected space with 
mod $p$ homology of finite type. Assume also that
\[ H^*(X;\FF_p )=\FF_p [x]/(x^{r+1}), \]
where $\alpha=|x|$ is even and $r\geq 1$. Put $\rho=(r+1)\alpha-2$.\\
1) If $p\mid (r+1)$ then 
\[ E_3^{*,*} \cong 
\big( \FF_p [u,\phi , q, \delta_0 , \delta_1 , \dots ,\delta_{p-2}] 
/I \big) \otimes \Gamma [\omega ], \]
where I is the ideal
\[ I=(\phi^{(r+1)/p},\medspace q^2,\medspace   
\delta_j u,\medspace \delta_j q ,\medspace \delta_j \delta_k
\medspace | 0\leq j\leq p-2, \medspace 0\leq k \leq p-2 ). \]
The bidegrees are 
$\Vert u \Vert = (2,0)$,
$\Vert \phi \Vert = (0,p\alpha )$, $\Vert q \Vert = (0,p\alpha -1)$,
$\Vert \delta_j \Vert = (0, j\alpha +\alpha -1)$ and
$\Vert \gamma_i (\omega ) \Vert =(0,\rho i)$.
The generators are represented by elements in the $E_2$ term as
follows:
\[ \quad u = [u], \quad \phi = [x^p], \quad q = [x^{p-1} dx], \quad 
\delta_j = [x^jdx], \quad
\gamma_i (\omega ) = [\gamma_i (\omega )]. \]
2) If $p\nmid (r+1)$ and the numbers $\kappa_t$ from 
theorem \ref{th:action} are zero for all $t$, then 
\begin{align*}
E_3^{*,*} \cong \FF_p [v_i^{(k)} ,w_i^{(k)} ,\trcl i h , u |
& 1\leq k \leq r, \medspace p\mid ((r+1)i+k), \\
& 1\leq h \leq r, \medspace p\nmid ((r+1)i+h),  
\medspace 0 \leq i] /I ,
\end{align*}
where $I$ is the ideal generated by the elements
\begin{align*}
& \trcl i h u,\quad \trcl i h w_j^{(\ell )},\quad 
\trcl i h \trcl j m ,\quad w_i^{(k)}w_j^{(\ell)}, \\
& v_i^{(k)} v_j^{(\ell )} - 
\epsilon_{r} (k+\ell ) \binom {i+j} i v_{i+j}^{(k+\ell)}, \\
& v_i^{(k)} w_j^{(\ell )} - 
\epsilon_{r} (k+\ell ) \binom {i+j} i w_{i+j}^{(k+\ell )}, \\
& v_i^{(k)} \trcl j h - 
\epsilon_{r} (k+h) \binom {i+j} i \trcl {i+j} {k+h}. 
\end{align*}
Here the number $\epsilon_{r}(s)$ equals $1$ if $1\leq s \leq r$ and
$0$ otherwise. The bidegrees of the generators are 
\begin{align*}
& \Vert v_i^{(k)}\Vert = (0, \rho i +\alpha k), \quad
\Vert w_i^{(k)}\Vert = (0, \rho i +\alpha k-1), \\
& \Vert \trcl i h\Vert = (0, \rho i +\alpha h-1), \quad
\Vert u \Vert = (2,0),
\end{align*}
and the generators are represented by elements in the $E_2$ term
as follows:
\[ v_i^{(k)}=[b_0^{k-1}b_i], \quad w_i^{(k)}=[b_0^{k-1}a_i], \quad
\trcl i h = [b_0^{h-1}a_i], \quad u=[u]. \]
\end{proposition}

\begin{proof}
1) Assume that $p\mid (r+1)$ such that $r+1 =mp$ for some $m\geq 1$.
By the K\" unneth formula we have that
\[ H(d) = H\big( \FF_p [x]/(x^{mp}) \otimes \Lambda (dx);d\big) 
\otimes \Gamma [\omega ]. \]
Since $d(x^j) = jx^{j-1}dx$ and
$d(x^jdx)=0$, the kernel and image of the differential on
$\FF_p [x]/(x^{mp})\otimes \Lambda (dx)$ has the 
following $\FF_p$-bases:
\begin{align*}
& \{ x^{kp}, x^j dx |  
0\leq k \leq m-1, 0 \leq j \leq mp-1 \} , \\
& \{ x^{j-1}dx |1 \leq j \leq mp-1, p\nmid j \} .
\end{align*}
It follows that 
$H(d)=\FF_p [\phi ]/(\phi^m )\otimes \Lambda (q) \otimes 
\Gamma [\omega ],$ where $\phi = x^p$ and $q=x^{p-1}dx$.
We can get the basis elements in $\image (d)$ by multiplying the 
elements in the set $\{ x^kdx | 0\leq k \leq p-2 \}$
by the elements $\gamma_i (\omega )$ and powers of $\phi$. 
The result follows.

2) Assume that $p\nmid (r+1)$. By remark \ref{rem:action} we see that
$\ker (d)$, $\image (d)$, $H(d)$ has respective $\FF_p$-bases
as follows:
\begin{align*}
& b_0^{j-1} a_i, \medspace b_0^{k-1} b_i & \text{ for } &
p\mid ((r+1)i+k), \\
& b_0^{h-1} a_i & \text{ for } &
p\nmid ((r+1)i+h), \\
& [b_0^{k-1}b_i], \medspace [b_0^{k-1}a_i] & \text{ for } & 
p\mid ((r+1)i+k),
\end{align*}
where $0\leq i$ and $1 \leq h,k,j \leq r$. The result follows.
\end{proof}

\begin{corollary} \label{cor:Poincare}
If $X$ satisfies the hypothesis of proposition \ref{th:E3Serre},
then the Poincar\'e series for the $E_3$ term of the
Serre spectral sequence is given by the following when $p\mid (r+1)$:
\[
\frac {1-t^{(r+1)\alpha}} {(1-t^{p\alpha})(1-t^{(r+1)\alpha -2})} \cdot
\big( \frac {1+t^{p\alpha -1}} {1-t^2}+
\frac {t^{\alpha -1}-t^{p\alpha -1}} {1-t^{\alpha}} \big) ,
\]
and by the following when $p\nmid (r+1)$:
\[
\frac 1 {1-t^2} \big( 1+P_{\IF (r,p,\alpha )}(t)+
\frac 1 t P_{\IF (r,p,\alpha )}(t) \big)+
\frac 1 t P_{\IT (r,p,\alpha )}(t).
\]
\end{corollary}

By proposition \ref{imd} and proposition \ref{th:E3Serre} we have

\begin{corollary}
If $X$ satisfies the hypothesis of proposition \ref{th:E3Serre} then,
\begin{enumerate}
\item If $p\mid (r+1)$ then the element
$\gamma_i(\omega ) \phi^j \delta_k\in E_3$ survives
to $E_\infty$ and is represented by 
$\tau (x^{pj+k+1}\gamma_i (\omega ))$ for 
$0\leq i$, $0\leq j \leq (r+1)/p$ and $0\leq k \leq p-2$. \\
\item If $p\nmid (r+1)$ then the generator $\trcl i h \in E_3$ survives
to $E_\infty$ and is represented by
$\tau (b_0^{h-1}b_i)$ for $1\leq h\leq r$ with $p\nmid ((r+1)i+h)$
and $0\leq i$.
\end{enumerate}
\end{corollary}

\section{Comparing the two spectral sequences.}
\label{se:MorseSerreSS}

In this section we will complete the investigation of the Morse
spectral sequence $E_*(\Morse )(L\CP^r_{h{\TT}})$. 
We write $\FC_0 \subseteq \FC_1 \subseteq \dots \subseteq L\CP^r$
for the energy filtration. So far, the main structural facts which 
we proved in section \ref{sec:MorseCPr} are the following:

\begin{enumerate}[{SF}(1)]
\item \label{Geometry}
The classes of even total degree are concentrated in
$\bigoplus_m E_*^{pm,*}(\Morse )(L\CP^r_{h{\TT}})$ (theorem
\ref{th:MorsePoincare}).
\item \label{SSmodule}
$E_*^{pm,*}(\Morse )(L\CP^r_{h{\TT}} )$ is a free $\FF_p [u]$ module. 
If $p\nmid n$, then  $E_*^{n,*}(\Morse )(L\CP^r_{h{\TT}} )$ is a finite
dimensional vector space (theorem \ref{th:MorsePoincare}).
\item \label{Oddeven}
Every non-trivial differential goes from even total degree to
odd total degree (lemma \ref{le:Preliminary}).
\item \label{Epcoll}
The spectral sequence collapses from the $E_p$ page (theorem
\ref{Epcollapse}).
\end{enumerate}

\begin{remark}
\label{re:filtsurjective}
It follows from SF(\ref{Oddeven}) that the inclusion 
$j:(\FC_n)_{h\TT} \hookrightarrow L\CP^r_{h\TT}$ induces a surjective
map $j^*:H^{\text{odd}}(L\CP^r_{h\TT} )\to H^{\text{odd}}((\FC_n)_{h\TT})$ 
for all $n\geq 0$.
\end{remark}

By SF(\ref{Geometry}) and SF(\ref{Epcoll}) we see that the 
only possibly non-trivial differentials in the spectral sequence are
\[ \xymatrix@C=1cm{
E_s^{pm,*}(\Morse )(L\CP^r_{h{\TT}} ) \ar[r]^-{d_s} &
E_s^{pm+s,*}(\Morse )(L\CP^r_{h{\TT}} )\cong 
E_1^{pm+s,*}(\Morse )(L\CP^r_{h{\TT}} )
} \]
for $1\leq s\leq p-1$.

We are also going to use the non-equivariant spectral sequence
$E_*(\Morse )(L\CP^r )$. From theorem \ref{le:nonequivariant}
we have the following structural facts:

\begin{enumerate}[{SF}(1)]
\setcounter{enumi}{4}
\item \label{Oddvanishing}
$E_1^{n,2i+1-n}(\Morse )(L\CP^r )=0$ if $p\mid (r+1)$ and $i\geq rn+1$ or if
$p\nmid (r+1)$ and $i\geq rn$.
\item \label{Morsesurject}
The map 
$E^{\text{odd}}_1(\Morse )(L\CP^r_{h\TT} )
\to E^{\text{odd}}_1(\Morse )(L\CP^r )$ is surjective. 
\end{enumerate}

\begin{remark}
\label{re:Homologysurject}
By SF(\ref{Morsesurject}) and SF(\ref{Oddeven}) the map
$E^{\text{odd}}_\infty (\Morse )(L\CP^r_{h\TT} )
\to E^{\text{odd}}_\infty (\Morse )(L\CP^r )$
is a surjection. A filtration argument shows that
if a map in some degree induces a surjective map on the $E_\infty$ pages,
then it also induces a surjective map on the cohomology
of the targets of the spectral sequences. So   
$H^{\text{odd}}(L\CP^r_{h\TT} ) \to H^{\text{odd}}(L\CP^r )$ is also 
surjective.
\end{remark}

Our plan for computing $H^*(L\CP^r_{h\TT})$ goes as follows.
First, we concentrate on the odd part of $H^*(L\CP^r_{h\TT})$.
The sum of the odd dimensional cohomology groups
$H^{\text{odd}}(L\CP^r_{h\TT} )$ is a submodule of 
$H^*(L\CP^r_{h\TT} )$ over 
$H^*(B\TT )$, since this ring is concentrated in even degrees.
We will list a set of elements in this module, and use the above 
properties of the spectral sequences to show that these elements 
are generators. We also give the relations satisfied by these elements,
thus computing the odd part of the Borel cohomology. 

To determine the odd part of the cohomology is not quite
the same as to determine the differentials in the Morse
spectral sequence, but in our situation it is close enough. 
We can use the knowledge of the odd cohomology together with 
the spectral sequences to determine the even dimensional cohomology. 

The first step in this program is to find elements in
$H^\text{odd}(L\CP^r_{h\TT})$ which can serve as generators. 
We use the results on the transfer map to find the right elements.  

Consider the ${\TT}$ transfer map $\tau$ in the context
of the Morse filtration of $L \CP^r$. Let 
$i: L\CP^r \to L\CP^r_{h{\TT}}$ be the inclusion.
It follows from theorem \ref{th:s1transfer} that the composite  
\[ \xymatrix@C=1cm{
H^{*+1}(L\CP^r )\ar[r]^-{\tau} & H^*(L\CP^r_{h{\TT}} )\ar[r]^-{i^*} 
& H^* (L\CP^r )
} \]
equals the action differential $d$.

We can now chose one half of our generators. This bunch of generators 
come with the relation that the generators are annihilated by 
multiplication by $u$.

\begin{lemma}
\label{le:torsionGen}
There is a graded subgrup 
$\Trans^* \subseteq H^{\text{odd}}(L\CP^r_{h{\TT}} )$
such that
\begin{enumerate}
\item $u\Trans^*=0$.
\item The map $i^*: H^*(L\CP^r_{h{\TT}} )\to H^*(L\CP^r)$ 
restricts to an injective map on $\Trans^*$. 
\item The subgroup $i^*(\Trans^* )\subseteq H^*(L\CP^r )$ 
agrees with the image of the composite map   
$i^*\circ \tau :H^{*+1}(L\CP^r )\to H^*(L\CP^r )$. 
\end{enumerate}
\end{lemma}

\begin{proof}
We chose a graded subgroup 
$\overline \Trans^{*+1} \subseteq \tilde H^{*+1}(L\CP^r )$
which maps isomorphically to the image of $d$. That is,
we chose (arbitrarily) a splitting of the surjective map
$d: H^{*+1}(L\CP^r )\to \im(d)^*$. 
Then $\Trans^* =\tau (\overline \Trans^{*+1} )\subseteq 
H^*(L\CP^r_{h{\TT}} )$ is a subgroup which by its definition satisfies 
the second and third property. It also satisfies the first property, 
since $u\tau=0$ by theorem \ref{th:s1transfer}. 
\end{proof}

Our second bunch of generators is not involved with any relations.
\newcommand{\Free}{{\mathcal U}}

\begin{lemma} \label{le:freeGen}
There is a graded subgrup 
$\Free^* \subseteq H^{\text{odd}}(L\CP^r_{h{\TT}} )$
such that the composite 
\[ \xymatrix@C=1cm{
\Trans^* \oplus \Free^* \ar[r] &
H^{\text{odd}}(L\CP^r_{h{\TT}} )\ar[r]^-{i^*} &
H^{\text{odd}}(L\CP^r)
} \]
is  an isomorphism. In addition to this, the restriction
\[ \xymatrix@C=1cm{
\Free^{2i+1} \ar[r] &
H^{2i+1}(L\CP^r_{h{\TT}} ) \ar[r]^-{j^*} &
H^{2i+1}((\FC_{pm})_{h{\TT}} )
} \]
is trivial if either $p\nmid (r+1)$ and $i\geq rpm$,
or $p\mid (r+1)$ and $i\geq rpm+1$.
\end{lemma}

\begin{proof}
To construct $\Free^*$, we  make a choice 
$\overline \Free^* \subseteq H^{\text{odd}}(L\CP^r )$ of 
a complementary subgroup of $i^*(\Trans )$, so that we have a 
direct sum decomposition of vector spaces
$H^{\text{odd}}(L\CP^r )\cong i^*(\Trans^* )\oplus \overline \Free^*$.
We intend to find $\Free^* \subseteq H^{\text{odd}}(L\CP^r_{h{\TT}} )$
such that $i^*$ maps this subgroup isomorphically to
$\overline \Free^*$ and such that the statement for the
restriction map holds.

According to the long exact sequence of
theorem \ref{th:s1transfer}, the following diagram has exact rows.
Because of remark \ref{re:filtsurjective} the left and middle 
vertical maps are surjections, because of remark \ref{re:Homologysurject}
the upper right horizontal map is a surjection. 
\[ \xymatrix@C=2cm{
H^{2i-1}(L\CP^r_{h{\TT}} ) \ar@{>>}[d] \ar[r]^-{u} &
H^{2i+1}(L\CP^r_{h{\TT}} ) \ar@{>>}[d] \ar@{>>}[r]^-{i^*} &
H^{2i+1}(L\CP^r ) \ar[d] \\
H^{2i-1}((\FC_{pm} )_{h{\TT}} )\ar[r]^-{u} &
H^{2i+1}((\FC_{pm} )_{h{\TT}} )\ar[r] &
H^{2i+1}(\FC_{pm} ). \\
} \]

Assume that $p\nmid (r+1)$ and $i\geq rpm$, or
$p\mid (r+1)$ and $i\geq rpm+1$. By SF(\ref{Oddvanishing}) this
implies that $H^{2i+1}(\FC_n , \FC_{n-1})=0$ for $0\leq n \leq pm$ 
such that $H^{2i+1}(\FC_{pm})=0$. Thus  
$\overline \Free^{2i+1}$ is contained in the kernel of the 
right vertical map.

The rest is a diagram chase. By the surjectivity of the upper 
right map $\overline \Free^*$ is the isomorphic image of a 
subgroup of $H^*(L\CP^r_{h{\TT}})$. The degree $2i+1$ part of this 
subgroup might not itself be in the kernel of the middle vertical map, 
but using the surjectivity of the left vertical map, we can replace 
it with a subgroup $\Free^{2i+1} \subseteq H^{2i+1}(L\CP^r_{h{\TT}} )$
which also maps isomorphically to $\overline \Free^*$, and 
such that $\Free^{2i+1}$ is in the kernel of the middle vertical map.
\end{proof}

\begin{remark}
\label{re:dimension}
It follows by \ref{action:cokernelsum1}. and
\ref{action:cokernelsum2}. of lemma \ref{le:countingtransfer}, 
that if $p\nmid (r+1)$, the dimension of the group 
$\bigoplus_{1\leq 2k-1\leq 2rpm-1} \Free^{2k-1}$ is $rm$.
If $p\mid (r+1)$ the dimension of the group 
$\bigoplus_{1\leq 2k-1\leq 2rpm+1} \Free^{2k-1}$ is 
$(r+1)m$. 
\end{remark}

\begin{theorem}
\label{th:mainodd}
There is a  map of $\FF_p [u]$-modules
\[
h_1 \oplus h_2:(\FF_p [u]\otimes \Free^* )\oplus \Trans^*
\to H^{\text{odd}}(L\CP^r_{h{\TT}} )
\]
which is an isomorphism of $\FF_p [u]$-modules.
\end{theorem}

\begin{proof}
We can extend the inclusion of $\Free^*$ in a unique way to an 
$\FF_p [u]$-linear map 
$h_1:\FF_p [u]\otimes \Free^* \to H^{\text{odd}}(L\CP^r_{h{\TT}} )$.  
Because of lemma  \ref{le:torsionGen} the inclusion of $\Trans^*$
is already an $\FF_p[u]$-linear map 
$h_2:\Trans^* \hookrightarrow H^{\text{odd}}(L\CP^r_{h{\TT}} )$. 

Lemma \ref{le:freeGen} and the exact sequence
\[ \xymatrix@C=2cm{
H^{2i-1}(L\CP^r_{h{\TT}} ) \ar[r]^-{u} &
H^{2i+1}(L\CP^r_{h{\TT}} ) \ar[r]^-{i^*} &
H^{2i+1}(L\CP^r )
} \]
shows that $h_1\oplus h_2$ is surjective on indecomposables.
But then it is a surjective $\FF_p [u]$-linear map, and it
suffices to show that it is also injective.

The main step is to prove injectivity of the 
localized map $(h_1\oplus h_2)[1/u]$.
Because of the vanishing statement of lemma \ref{le:freeGen} 
we have a commutative diagram as follows, where $\chi_p (r+1)$ 
is the number defined in \ref{def:TF}:
\[ \xymatrix@C=2cm{
\FF_p[u] \otimes \bigoplus_{0 \leq i} \Free^{2i+1}
\ar[d]^-{id \otimes pr} \ar[r]^-{h_1} &
H^{\text{odd}}(L\CP^r_{h{\TT}} ) \ar@{>>}[d]^-{j^*} \\
\FF_p[u]\otimes \bigoplus_{0 \leq i\leq rpm-\chi_p (r+1)} \Free^{2i+1} 
\ar[r]^-{\overline h_1} &
H^{\text{odd}}((\FC_{pm})_{h{\TT}} ).
} \]
The map $j^*$ is surjective by remark \ref{re:filtsurjective}. 

We localize by inverting $u$. Because all modules are of finite type, 
this localization agrees with tensoring with $\FF_p [u,u^{-1}]$ over
$\FF_p [u]$. Since $h_1\oplus h_2$ is surjective,
and the localization of $h_2$ vanishes, we know that
\[ \xymatrix@C=2cm{
h_1 [ \frac 1 u ] : \FF_p [u,u^{-1}] \otimes \Free^* \ar[r] &
H^{\text{odd}}(L\CP^r_{h{\TT}} ) [ \frac 1 u ]
} \]
is surjective. It follows by the diagram that the map 
\[ \xymatrix@C=2cm{
\overline h_1[\frac 1 u ] : \FF_p [u,u^{-1}] \otimes
\bigoplus_{0 \leq i\leq rpm-\chi_p (r+1)} \Free^{\text{2i+1}} \ar[r] &
H^{\text{odd}}((\FC_{pm})_{h{\TT}}) [\frac 1 u ]
} \]
is also surjective. We know from remark \ref{re:dimension} 
that as abstract module, the domain space is given by
\[
\FF_p [u,u^{-1}] \otimes \bigoplus_{0 \leq i\leq rpm-\chi_p (r+1)}
\Free ^{\text{2i+1}} \cong
\begin{cases}
\FF_p [u,u^{-1}]^{\oplus (r+1)m} & \text{if } p\mid (r+1) \\
\FF_p [u,u^{-1}]^{\oplus rm}   & \text{if } p\nmid (r+1).  
\end{cases}
\]
The target space is determined by theorem \ref{th:relabel}
and corollary \ref{cor:oddcounting}. By the proof of theorem
\ref{th:relabel} the re-indexing is given by multiplying the 
filtration degree by $p$. We find that 
\[
 H^{\text{odd}}((\FC_{pm})_{h{\TT}} ) \left[ \frac 1 u \right] \cong
\begin{cases}
\FF_p [u,u^{-1}]^{\oplus (r+1)m} & \text{if } p\mid (r+1) \\
\FF_p [u,u^{-1}]^{\oplus rm}   & \text{if } p\nmid (r+1).  
\end{cases}
\]

The punch line is that $h_1[1/u]$ is a surjective map between two 
finitely generated, free $\FF_p [u,u^{-1}]$-modules of the same rank. 
But then the map has to be an isomorphism.

This proves the injectivity of the localization $h_1$. To get
the injectivity of $h_1\oplus h_2$, we consider an element 
$(c,t)\in \ker (h_1\oplus h_2)$. We have to prove that this element 
is trivial. Since the localization of $t$ vanishes, the localization 
of $c$ is in the kernel of the localization of $h_1$. This localized 
map is injective, so the localization of $c$ vanishes. But the
canonical map
$\FF_p [u] \otimes \Free^* \to \FF_p [u,u^{-1}] \otimes \Free^*$
is injective, so $c$ itself vanishes, and $t$ is in the kernel of 
$h_2$. But $h_2$ is injective on $\Trans$, which proves that $t=0$.
\end{proof}

We can finally complete the main calculation of the paper.
Recall the index sets $\IT (r,p,2)$ and $\IF (r,p,2)$ 
from definition \ref{def:TF}. We need a small perturbation
as follows:
\[ \IF^\prime (r,p,2)=
\begin{cases}
(\IF (r,p,2) \setminus 2r\NN) \cup (2+2r\NN) & 
\text{if } p \mid (r+1), \\
\IF (r,p,2) & \text{if } p\nmid (r+1).
\end{cases}
\]
This makes sense, since $2r\NN \subseteq \IF(r,p,2)$ by 
lemma \ref{le:twosets} when $p\mid (r+1)$.

\begin{theorem} 
\label{th:main}
As a graded $\FF_p [u]$-module,
$H^*(L\CP^r_{h\TT} ;\FF_p )$ is isomorphic to the direct sum
\[
\FF_p [u] \oplus \bigoplus_{2k\in \IF (r,p,2)} \FF_p [u] f_{2k}
\oplus \bigoplus_{2k\in \IF^{\prime} (r,p,2)} \FF_p [u]f_{2k-1}
\oplus \bigoplus_{2k\in \IT(r,p,2)} (\FF_p [u]/(u)) t_{2k-1}.  
\] 
Here the lower index denotes the degree of the generator.
\end{theorem}

\begin{proof}
In theorem \ref{th:mainodd} we proved the formula for
the group of odd degree elements. We have to show that
$H^\text{even} (L\CP^r_{h\TT} )$ is a free $\FF_p [u]$-module, 
with generators in the stated degrees. 
Put $E_r^{**} = E_r^{**}(\Morse )(L\CP^r_{h\TT} )$.
From SF(\ref{Geometry}) and SF(\ref{SSmodule}) we see that
$E_1^{\text{even}}$ is a free $\FF_p [u]$-module. The degrees of the 
generators can be read off theorem \ref{th:MorsePoincare}. 
We see that 
\begin{align*}
E_1^{(0,*)(\text{even})}
& \cong \oplus_{0\leq 2i\leq 2r}{\FF_p [u] x_{2i}}, 
\\
E_1^{(pm,*)(\text{even})}
&\cong
\begin{cases}
\oplus_{2rpm\leq 2i\leq 2r(pm+1)} \FF_p [u]x_{2i} & 
\text{if } p\mid (r+1), \\
\oplus_{2rpm+2\leq 2i\leq 2r(pm+1)} \FF_p[u]x_{2i} & 
\text{if } p\nmid (r+1).
\end{cases}
\end{align*}
Because of SF(\ref{Epcoll}) together with SF(\ref{SSmodule})
we see that 
$E_\infty^{(pm,*)(\text{even})}$ is a
submodule of finite index in $E_1^{(pm,*)(\text{even})}$. By abstract
structure theory of $\FF_p [u]$-modules, it follows that
$E_\infty^{(pm,*)(\text{even})}$ is a free module on certain
generators. If we can filter a graded module with quotient which are 
free modules, the original module was also a free module. The 
generators are in the same degrees as the generators for the direct 
sum of the quotients.

What is left is to figure out the degrees of the generators of the 
free module $E_\infty^{(pm,*)(\text{even})}$. It would be best if we 
could compute the differentials. Unfortunately, we cannot do this. 
What we can do, is to compute the dimension of
$E_{\infty}^{\text{even}}$ in each degree, and recover the degrees of 
the generators from this.

To do this, we consider Poincar\'e series. Let 
$P(t)$, $P^\text{even}(t)$ and $P^\text{odd}(t)$
be the Poincar\'e series of 
$H^*(L\CP^r_{h\TT})$, $H^{\text{even}}(L\CP^r_{h\TT})$ and
$H^\text{odd}(L \CP^r_{h\TT})$ respectively.
Similarly, let $P_r (t)$, $P_r^{\text{even}}(t)$ and
$P_r^{\text{odd}}(t)$ be the Poincar\'e series of $E_r$, the 
even total degree part of $E_r$ and the odd total degree part $E_r$
respectively. 

The series $P_1^{\text{even}}(t)$ and $P_1^{\text{odd}}(t)$ can be 
recovered from $P_1(t)$ as the sum of the even respectively the 
sum of the odd monomials occurring in it.
From the computation in Lemma \ref{le:MorsePoincare},
we recall that if $p\mid (r+1)$, then 
\[
P_1(t)=\frac 1 {1-t} \cdot \frac {1-t^{2r+2}} {(1-t^2)(1-t^{2pr})}
\]
If follows that 
\[
P_1^{\text{even}}(t)= 
\frac 1 {1-t^2} \cdot \frac {1-t^{2r+2}} {(1-t^2)(1-t^{2pr})}, \quad
P_1^{\text{odd}}(t)=
\frac t {1-t^2} \cdot \frac {1-t^{2r+2}} {(1-t^2)(1-t^{2pr})}.
\]

The important thing here is that $P_1^\text{odd}(t)=tP_1^\text{even}(t)$.
This means that the dimensions of the group of classes of total
degree $2i$ in $E_1$ agrees with the dimension of the classes  
of total degree $2i+1$. Now recall that by SF(\ref{Oddeven}), all
differentials in the spectral sequence go from even total degree to
odd total degree. It follows that for any $r$, including
$r=\infty$, we have that $P_r^{\text{odd}}(t)=tP_r^{\text{even}}(t)$,
and also that $P^{\text{odd}}(t)=tP^{\text{even}}(t)$.

In case $p\nmid (r+1)$ the formula for $P_1(t)$ in 
lemma \ref{le:MorsePoincare} is different, but the argument 
above still applies, so that $P^{\text{odd}}(t)=tP^{\text{even}}(t)$ 
in both cases.

We can use theorem \ref{th:mainodd} to give an expression 
for $P^{\text{odd}}(t)=P_\infty^{\text{odd}}(t)$. Using the notation of
\ref{def:TF}, we see that
\[
tP^{\text{odd}}(t)=P_{\IT (r,p,2)}(t)+\frac 1 {1-t^2} 
P_{\IF^\prime (r,p,2)}(t). 
\]
By definition of $\IF^\prime (r,p,2)$, we have
\begin{equation} \label{eq:IFprime}
 P_{\IF^\prime (r,p,2)}(t)=P_{\IF (r,p,2)}(t)
-(1-\chi_p(r+1))\frac {t^{2r}(1-t^2)} {1-t^{2r}},
\end{equation}
which we insert above and find
\[ 
tP^{\text{odd}}(t)=P_{\IT (r,p,2)}(t)+\frac 1 {1-t^2} 
P_{\IF (r,p,2)}(t)-(1-\chi_p(r+1))\frac {t^{2r}} {1-t^{2r}}.
\]
Rewriting the first two terms, and using the result 
$P^{\text{odd}}(t)=tP^{\text{even}}(t)$, we see that
\begin{align*}
t^2P^{\text{even}}(t) &= tP^{\text{odd}}(t) \\
&= P_{\IT (r,p,2)}(t)+ P_{\IF (r,p,2)}(t) +  
\frac {t^2} {1-t^2} P_{\IF(r,p,2)}(t) 
-(1-\chi_p(r+1))\frac {t^{2r}} {1-t^{2r}}.
\end{align*}
We rewrite the sum of the first two terms 
using lemma \ref{le:Poincaretwosets} and obtain
\[
t^2P^{\text{even}} (t)=\frac {t^2} {1-t^2} + 
\frac {t^2} {1-t^2} P_{\IF (r,p,2)}(t),
\]
which completes the proof.
\end{proof}

\begin{corollary}
\label{cor:actionproof}
The numbers $\kappa_i$ in theorem \ref{th:action} are always zero.
\end{corollary}

\begin{proof}
In theorem \ref{th:action} this is proved for $\alpha \geq 4$, and
in some cases also when $\alpha =2$. So assume that $\alpha =2$.
Then an obstruction argument shows that $X$ is homotopy
equivalent to $\CP^r$. So we can without loss of generality assume 
that $X=\CP^r$. Since the action differential factors over the
transfer map, it sufficient to show that the transfer map
\[
\tau :H^{\text{odd}}(L\CP^r ) \to H^{\text{even}}(L\CP^r_{h\TT} )
\]
is the zero map. The image of the transfer map is annihilated
by multiplication by $u$ (theorem \ref{th:s1transfer}). 
According to theorem \ref{th:main} we also have that
$H^{\text{even}}(L\CP^r_{h\TT} )$ is a free module over $\FF_p [u]$, so
multiplication by $u$ is injective. The result follows.
\end{proof}

We can now give an other version of our main result.
Recall that
\[ IF (r,p,2)= \{ 2(ri+j) \mid 0\leq i,\medspace 
\chi_p (r+1) \leq j \leq r,\medspace p\mid ((r+1)i+j) \}, \]
where $\chi_p (s)$ equals $0$ when $p$ divides $s$ and
$1$ when $p$ does not divide $s$.
\begin{theorem} \label{th:main2}
Let $\{ E_*\}$ be the mod $p$ Serre spectral sequence for the
fibration sequence $L\CP^r \to (L\CP^r )_{h\TT} \to B\TT$. That is
\[ E_2^{*,*}=H^*(B\TT; \FF_p )\otimes  H^*(L\CP^r ;\FF_p )
\Rightarrow H^*((L\CP^r )_{h\TT}; \FF_p ). \]
For any positive integer $r$ and any prime $p$ one has that
$E_3=E_\infty$. Furthermore, the Poincar\' e 
series $P_{r,p}(t)$ for $H^*((L\CP^r )_{h\TT};\FF_p )$ is given by
\[ 
P_{r,p}(t)=\frac 1 {1-t} \big( 1+\sum_{k\in \IF (r,p,2)} t^k \big) .
\]
If $p$ divides $r+1$ we can rewrite this as 
\[
P_{r,p}(t)= \frac {1-t^{2(r+1)}} {(1-t)(1-t^{2r})(1-t^{2p})}. 
\] 
\end{theorem}

\begin{remark}
We have described the $\FF_p$-algebra structure of the 
$E_\infty$ page in proposition \ref{th:E3Serre} with $\alpha =2$.
\end{remark}

\begin{proof}
We first use theorem \ref{th:main} to prove that the Poincar\' e 
series is as stated. By this theorem we have that
\begin{equation} \label{eq:mainpoincare}
P_{r,p}(t) = \frac 1 {1-t^2} 
\big( 1+P_{\IF (r,p,2)}(t)+\frac 1 t P_{\IF^\prime (r,p,2)}(t)\big)
+ \frac 1 t P_{\IT (r,p,2)}(t).
\end{equation}
By using equation (\ref{eq:IFprime}) and lemma 
\ref{le:Poincaretwosets} we can rewrite this as
\[ P_{r,p}(t)= \frac 1 {1-t^2} 
\big( 1+P_{\IF (r,p,2)}(t)+\frac 1 t P_{\IF (r,p,2)}(t) \big)
-\frac 1 t P_{\IF (r,p,2)}(t)+\frac t {1-t^2}. \]
The desired result follows by a small reduction.

When $p$ divides $r+1$, we can write the index set as  
\[ \IF (r,p,2)=
\{ 2(ri+pm)\mid 0\leq i, \medspace 0\leq m \leq \frac {r+1} p -1\}
\setminus \{ 0\} .\] 
The last formula for $P_{r,p}(t)$ follows. 

Since corollary \ref{cor:actionproof} holds, we have a formula for 
the Poincar\' e series of the $E_3$ term in corollary
\ref{cor:Poincare} with $\alpha =2$. This series is the same
as $P_{r,t}(t)$ when $p$ does not divide $r+1$ by 
(\ref{eq:mainpoincare}). When $p$ divides $r+1$ the two series
also agree by our last formula for $P_{r,p}(t)$. Thus, the
Serre spectral sequence collapses from the $E_3$ page.
\end{proof}

The proof of the main theorem involves the existence of
non-trivial differentials in the Morse spectral sequence. 
This fact has a certain geometrical content which we describe
below for $r\geq 2$. By similar methods one can get the
same result for $r=1$, but we will not go into this here.

\begin{corollary}
\label{cor:elastomania}
For any $r\geq 2$ and any $n\geq 2$ there is a
trajectory of loops on $\CP^r$ which converges in positive time 
towards a geodesic with period $n$ and in negative time
towards a non-constant geodesic with period $n+1$. 
\end{corollary}

\begin{proof}
Assume that there are no such trajectories. Then the geometric map
\[
(\FC_{n+1} /\FC_{n} )_{h\TT} \to \Sigma (\FC_n )_{h\TT} \to
\Sigma (\FC_n /\FC_{n-1} )_{h\TT} ,
\]
which induces the $d_1$ differential in the Morse spectral 
sequence, is nullhomotopic. So in the mod $p$ Morse spectral 
sequence $E_*=E_* (\Morse) (L\CP^r_{h\TT} )$, we have  
$d_1=0 : E_1^{n,*}\to E_1^{n+1,*}$ for any prime $p$. We intend 
to show that you can always chose a prime $p$, such that this 
$d_1$ must be non trivial, leading us to the desired contradiction.

Since $n\geq 2$, there is a prime $p$ that divides $n$, say $n=pm$. 
We actually claim that any such $p$ will do.
Assume to the contrary, that $d^1$ is trivial on
$E_1^{pm,*}$. The generators of the $\FF_p [u]$-module $E_1^{pm,*}$ 
have degrees less or equal to $2r(pm+1)$. Since the lowest class 
in $E_1^{pm+s,*}$ is in degree $(2r)(pm+s-1)+1$, the only possible 
non-trivial differential originating in $E_*^{pm,*}$ would be a
$d_2$ hitting the non-trivial class in $E_2^{pm+2,2r(pm+1)-pm-1}$.
But comparison to the non-equivariant spectral sequence shows that
this class is a permanent cycle, so every class in $E_1^{pm,*}$ is 
a permanent cycle. (cf. the proof of theorem \ref{Epcollapse}).

But this means that every generator of the free 
$\FF_p[u]$-module $E_1^{(pm,*)\text{even}}$ corresponds to a 
generator in the free $\FF_p [u]$-module 
$H^{\text{even}}(L\CP^r_{h\TT})$. So this modules has a generator 
in every degree $2k$ where $2rpm+2 \leq 2k \leq 2r(pm+1)$. 
In case $p\mid (r+1)$, in addition to this, it has a generator 
in degree $2k=2rpm$. We have listed the
generators of this module in theorem \ref{th:main},
so this says that all these numbers $2k$ are contained
in $\IF^\prime (r,p,2)$.     

The rest of the proof is just very elementary number theory.
As usual, there are two possible cases. 
We first deal with the case $p\mid (r+1)$. According to
the definition of $\IF^\prime (r,p,2)$ just prior to
theorem \ref{th:main}, this set does not contain
any numbers divisible by $2r$. So we cannot have that
$2rpm \in \IF^\prime (r,p,2)$, contradicting the assumption.

So now assume that $p\nmid (r+1)$. We first show that there 
are not many pairs of consecutive even numbers in 
$\IF (r,p,2)= \IF^\prime(r,p,2)$. Assume that both $2k$ and $2k+2$ 
are contained in $\IF (r,p,2)$. As in definition
\ref{def:TF} we find four numbers
$i_1,i_2,j_1,j_2$  pairwise satisfying the 
conditions mentioned in that definition, and such that
$2k=2ri_1+2j_1$ and $2k+2=2ri_2+2j_2$.
Then $2=2r(i_2-i_1)+2(j_2-j_1)$. Since
$\mid j_1-j_2\mid \leq 2(r-1)$ it follows that
either $i_1=i_2$ and $j_2=j_1+1$, or
$i_2=i_1+1$, $j_2=0$ and $j_1=r-1$. The first
possibility together with the congruence condition
leads to $p\mid ((r+1)i_1+j_1)$ and
$p\mid ((r+1)i_1+j_1+1)$, which is a contradiction.

The second possibility leads to $p\mid ((r+1)i_1+r-1)$
and $p\mid ((r+1)(i_1+1))$. Subtracting, we get that
$p\mid 2$, that is $p=2$. So the only possibility for
that $\{ 2k,2k+2\}\in \IF (r,p,2)$ is that
$p=2$, $r$ is even, $j_1=r-1$ is odd and  
$k=ri_1+j_1$ is odd. In particular, we can never have that 
$\{ 2k,2k+2,2k+4\} \in \IF (r,p,2)$. Since the set of $k$ 
such that $2rpm+2 \leq 2k \leq 2rp(m+1)$ is a set of
$r$ consecutive numbers, we have that $r=2$.

We are now reduced to showing that $\IF (2,2,2)$ does not
contain two numbers of the form $8m+2,8m+4$. We already noted that
if $2k, 2k+2\in \IF (r,p,2)$, then $k$ is odd.
So if $2k=8m+4$, this particular $k$ does not qualify.
The proof is complete.
\end{proof}

\section{Appendix: The circle transfer map}
\label{Appendix:s1transfer}

In this appendix we describe an elementary construction of a $\TT$ 
transfer map $\tau$. There are several discussions of transfer maps in the 
literature, and also of extensive refinements and generalizations. See 
for example \cite{MS} or \cite{MMM}. For our purposes, we need only a 
rather coarse version. We give a simple, self contained construction
close to the one we gave in \cite{Fund}. 

\begin{theorem}
\label{th:s1transfer}
Let $X$ be a based $\TT$-CW complex such that the action of $\TT$ 
is free away from the base point. Write $q:X\to X/\TT$ for the 
canonical projection. Let $R$ be a principal ideal domain and 
$M$ an $R$-module. There is a linear map which is natural in 
$X$ and $M$ as follows:
\[
\tau : \tilde H^n (X;M) \to \tilde H^{n-1}(X/\TT ;M).
\]
It is the connecting homomorphism in a long exact Gysin sequence
\[
\xymatrix@C=0.6cm{
\dots \ar[r] 
& \tilde H^{n-2}(X/\TT ;M) \ar[r]^-{\gamma} 
& \tilde H^{n}(X/\TT ;M) \ar[r]^-{q^*} 
& \tilde H^{n}(X;M) \ar[r]^-{\tau} 
& \tilde H^{n-1}(X/\TT ;M) \ar[r]
& \dots 
}\]
in particular, $\tau \circ q^*=0$. Assume that $M=R$. 
Then Frobenius reciprocity holds 
\[
\tau (a q^* (b))=\tau(a) b,
\]
and $q^*\circ \tau=d$ where $d$ denotes the action 
differential. 
\end{theorem}

\begin{remark}
We have $H^*(B\TT ;R)=R[u]$ where 
$\deg u=2$ which gives us a class $pr_1^*(u) \in H^2(X_{h\TT};R)$. 
The map $\gamma$ in the Gysin sequence is given by 
multiplication by this class:
\[
\xymatrix@C=1.5cm{
H^{n-2}(X_{h\TT},B\TT ;M) \ar[r]^-{pr_1^*(u)\cdot } 
& H^{n}(X_{h\TT},B\TT ;M) \\ 
\tilde H^{n-2}(X/\TT ;M) \ar[r]^-{\gamma} \ar[u]^-{\cong}  
& \tilde H^{n}(X/\TT ;M). \ar[u]^-{\cong}
}\]
\end{remark}

\begin{proof}
The key point is to compare $X$ to $E\TT \times X/E\TT$.
So we start by considering the spherical fibration
\[ 
\xymatrix@C=1cm{
\TT \ar[r] & E\TT \times X \ar[r]^-Q & E\TT \times_{\TT} X, } 
\]
together with the two subspaces $E\TT= E\TT\times *$ of $E\TT \times X$
and $B\TT = E\TT \times_{\TT} *$ of $E\TT \times_{\TT} X$ with
$Q^{-1} (B\TT ) = E\TT$.

The fibration $Q$ is a pullback of the fibration $\TT \to E\TT \to B\TT$
along the projection map $pr_2$. Since $B\TT$ is $1$-connected 
it follows that $Q$ is orientable. Thus there is a relative Gysin sequence
for $Q$ \cite[VII.5.12]{Whitehead}. We let $\tilde \tau$ be the 
connecting homomorphism in this sequence.

Let $X_0$ denote $X$ with trivial $\TT$ action. Pick a point 
$e\in E\TT$ and define the $\TT$ map 
$\theta :\TT \times X_0 \to E\TT \times X$ 
by $(z,x)\mapsto (ez^{-1},zx)$.
We have a commutative diagram
\[ 
\xymatrix@C=1cm{
(\TT \times X_0, \emptyset ) \ar[r]^{pr_2} \ar[d]^{\theta} 
& (X, \emptyset ) \ar[d] \\
(E\TT \times X, E\TT ) \ar[r]^{Q} 
& (E\TT \times_\TT X,B\TT ).}
\]
By naturality of the Gysin sequence we have the following diagram,
where the cohomology groups have coefficients in $M$:
\[ \xymatrix@C=0.6cm{
\ar[r] 
& H^{n-2}(X_{h\TT}, B\TT ) \ar[r]^-{u\cdot } \ar[d]^-{Q^*}
& H^{n}(X_{h\TT}, B\TT ) \ar[r]^-{Q^*} \ar[d]^-{Q^*}
& \tilde H^{n}(X) \ar[r]^-{\tilde \tau} \ar[d]^-{\theta^*}
& H^{n-1}(X_{h\TT}, B\TT ) \ar[r]^-{u\cdot} \ar[d]^-{Q^*}
& \\
\ar[r] 
& H^{n-2}(X) \ar[r]^-{0}
& H^{n}(X) \ar[r]^-{pr_2^*}
& H^{n}(\TT \times X) \ar[r]^-{\overline \tau} 
& H^{n-1}(X) \ar[r]^-{0}
& 
} \]

By the upper sequence, $\tilde \tau \circ Q^*=0$. Write 
$\eta : \TT \times X \to X$ for the action map. $\eta*$ and
$\theta^*$ are the same in positive degrees. Assume that $M=R$
Then $\eta^* (y)=1\otimes y+v\otimes dy$ where $v$ has degree $1$. 
By the diagram $\overline \tau (1\otimes a)=0$ and 
$\overline \tau (v\otimes b)=b$ such that 
$Q^*\circ \tilde \tau =\overline \tau \circ \eta^*=d$.
Finally, we have Frobenius reciprocity for $\tilde \tau$ by 
\cite[VII.5.16]{Whitehead}.

The map $pr_2: E\TT \times X/E\TT \to X$ is a homotopy equivalence
if we restrict it to the fixed points for any closed subgroup of 
$\TT$. So it is a $\TT$-equivariant homotopy equivalence,
see for instance \cite[II.2.7]{tomD}. Thus, we have a natural
isomorphism 
\[ pr_2^* : H^*(X/\TT ,*;M) \to H^*(X_{h\TT},B\TT ;M). \]
The result follows. 
\end{proof}

\section{Appendix: Rational coefficients}

In this appendix we use theorem \ref{th:action} to obtain results
for rational coefficients.

\begin{proposition} \label{pr:actionQ}
Let $X$ be a 1-connected space with $H_*(X;\ZZ )$ of finite type
and assume that
\[ H^*(X;\ZZ )=\ZZ [x]/(x^{r+1}), \]
where $\alpha = \deg (x)$ is even and $r\geq 1$. Put 
$\rho =(r+1)\alpha -2$. Then,
\[
H^k(LX;\QQ )= 
\begin{cases}
\QQ , & k\in \{ 0\} \cup 
\{\rho i+\alpha j, \rho i+\alpha j-1|0\leq i, 1\leq j \leq r \} ,\\
0 , & \text{otherwise.} 
\end{cases}
\]
The action differential 
$d: H^k(LX;\QQ )\to H^{k-1}(LX;\QQ )$
is an isomorphism when 
$k\in \{ \rho i +\alpha j|0\leq i, 1\leq j \leq r \}$
and zero otherwise. 
\end{proposition}

\begin{proof}
By the Serre spectral sequence for the fibration 
$\Omega X \to PX \to X$ we see that $H_*(\Omega X;\ZZ )$ is of finite
type. The Serre spectral sequence for $\Omega X \to LX \to X$ then
gives us that $H_*(LX;\ZZ )$ is of finite type.

Consider the universal coefficient sequence where
$A$ is an abelian group:
\[ \xymatrix{
0\ar[r] & \Ext (H_{k-1}(LX;\ZZ ), A) \ar[r] & H^k(LX;A) \ar[r]
& \Hom (H_k(LX;\ZZ ), A) \ar[r] & 0. 
} \]
By choosing $A=\ZZ /p$ for $p$ sufficiently large, we obtain 
that the $\Ext$ group is zero and that we can apply theorem
\ref{th:action} part 2). Thus,
\[
H_k(LX;\ZZ )/T_k(LX) \cong 
\begin{cases}
\ZZ , & k \in \{ 0\} \cup
\{\rho i+\alpha j, \rho i+\alpha j-1|0\leq i, 1\leq j \leq r\} ,\\
0 , & \text{otherwise,} 
\end{cases}
\]
where $T_k(LX)$ denotes the torsion subgroup of $H_k(LX;\ZZ )$. 
We then choose $A=\QQ$ and obtain the stated result for 
$H^k(LX;\QQ )$.

Let $\eta :\TT \times LX \to LX$ denote the action map. By definition
of the action differential $d$ we have that 
$\eta^* (y)=1\otimes y+v\otimes dy$ where $\deg (v)=1$ for
$\ZZ /p$ or $\QQ$ coefficients. There is also the projection
$pr_2 : \TT \times LX \to LX$ with $pr_2^*(y)=1\otimes y$.

By theorem \ref{th:action} we have $b_0^{j-1}b_i$ in $H^*(LX;\ZZ /p)$
with $\deg (b_0^{j-1}b_i)=\rho i +\alpha j$ and
\[ d(b_0^{j-1}b_i)=(j+(r+1)i)b_0^{j-1}a_i \]
for $0\leq i$, $1\leq j\leq r$. Fix such $i$ and $j$ and put
$k=\rho i+ \alpha j$. 

We use the universal coefficients sequence with abelian group 
$A$ again. By naturality we get commutative diagrams for 
$\eta$ and for $pr_2$. Choose $A=\ZZ /p$ where $p$ is so 
large that the $\Ext$ groups vanish, $p>r+1$ and $p>j+(r+1)i$. 
Then,
\[ 0\neq \eta_*-(pr_2)_*: H_k(\TT \times LX;\ZZ )/T_k(\TT \times LX)
\to H_k (LX;\ZZ )/T_k(LX).\]
We then take $A=\QQ$ and find that $0\neq \eta^*-pr_2^*$ in 
degree $k$ for $\QQ$ coefficients. Thus, the action differential 
$d:H^k(LX;\QQ)\to H^{k-1}(LX;\QQ )$ is an isomorphism as stated. 
Since $d\circ d=0$ the vanishing statement for $d$ follows.
\end{proof}

The ring structure of $H^*(LX;\QQ )$ was first computed by \v Svarc
in \cite{Svarc}. Combining his computation with the proposition 
above we obtain the following result:

\begin{proposition}
Let $X$ be as in proposition \ref{pr:actionQ}. Then,
\[ H^*(LX;\QQ )= \QQ [a_i,b_i|i\geq 0]/I, \]
where $I$ is the ideal generated by the following elements 
for $i,j\geq 0$:
\[ a_ia_j, \quad b_ib_j-b_0b_{i+j}, \quad b_ia_j-b_0a_{i+j}, \quad
b_0^rb_i, \quad b_0^ra_i. \]
The degrees of the generators are $|a_i|=\rho i+\alpha -1$ and
$|b_i|=\rho i +\alpha$. Furthermore, 
\[H^*(LX_{h\TT};\QQ )=\QQ [u]\otimes 
\QQ [w_i^{(j)}|0\leq i, 1\leq j \leq r]/J, \]
where $J$ is the ideal generated by all the products
$w_i^{(j)}w_k^{(\ell )}$.
The degrees of the generators are $|u|=2$ and
$|w_i^{(j)}|=\rho i + \alpha j -1$.
\end{proposition}

\begin{proof}
The result of \v Svarc's computation can be found in 
\cite[\S 6]{Klein}. Performing the substitutions $a_i=g_{i+1}$ for
$i\geq 0$ and $b_0=x$, $b_i=h_i$ for $i\geq 1$ in Theorem 6.2 and 6.3
of \cite{Klein}, we obtain the stated description of the cohomology of
$LX$. 

Consider the Serre spectral sequence for the homotopy orbit space.
Note that the elements $u^k$ for $k\geq 0$ are not hit by any
differentials since we may factor $id_{B\TT}$ as 
$B\TT \to LX_{h\TT} \to B\TT$, where we use a constant loop to define 
the first map and the second map is the projection $pr_1$.

The $d_2$ differential is given by the action differential $d$.
By proposition \ref{pr:actionQ} we see that only the elements
of the form $b_0^{j-1}a_i$ in $E_2^{0,*}$ and $u^k$ in $E_2^{*,0}$
survive to the $E_3$ page and that $E_3=E_\infty$.
The unique class $w_i^{(j)}$ in $H^*(LX_{h\TT};\QQ )$ 
representing $b_0^{j-1}a_i$ satisfies $(w_i^{(j)})^2=0$ since it has
odd degree. Thus any product $w_i^{(j)}w_k^{(\ell )}$ has
finite multiplicative order and hence cannot equal a nonzero 
constant times a power of $u$. Hence the multiplicative structure
is as stated.
\end{proof}

\end{document}